\documentclass[10pt]{article} 
\usepackage{geometry}                       
\usepackage{graphicx}
\usepackage[table]{xcolor}             
\usepackage{amssymb}
\usepackage{amsmath}
\usepackage{mathtools} 
\usepackage{svg} 
\usepackage[utf8]{inputenc}
\usepackage{bm} 
\usepackage{xcolor}

\usepackage{booktabs}
\newcommand{\E}{\mathrm{e}}
\usepackage{amsthm}
\usepackage[normalem]{ulem}
\usepackage{makecell}

\newtheorem{exmp}{Example}

\newtheorem{rmrk}{Remark}

\newtheorem{algo}{Algorithm}
\newtheorem*{wh}{Working Hypothesis}

\usepackage{MnSymbol}
\usepackage{algorithm}
\usepackage[noend]{algpseudocode}
\usepackage{cite}
\usepackage{listings}
\usepackage[]{hyperref}
\hypersetup{
    colorlinks = true,
    linkbordercolor = {white},
    linkcolor=black,
    citecolor=black,
    urlcolor=black
}
\makeatletter
\newcommand{\doubletilde}[1]{{%
  \mathpalette\double@tilde{#1}%
}}
\newcommand{\double@tilde}[2]{%
  \sbox\z@{$\m@th#1\tilde{#2}$}%
  \ht\z@=.9\ht\z@
  \tilde{\box\z@}%
}
\makeatother
\usepackage{amsthm}
\usepackage[export]{adjustbox}
\usepackage{mathdots}

\definecolor{jlbase}{HTML}{444444}            
\definecolor{jlkeyword}{HTML}{444444}         
\definecolor{jlliteral}{HTML}{78A960}         
\definecolor{jlbuiltin}{HTML}{397300}         
\definecolor{jlcomment}{HTML}{888888}         
\definecolor{jlstring}{HTML}{880000}          
\definecolor{jlbackground}{HTML}{F5F5F5} 
\lstdefinelanguage{Julia}%
  {morekeywords=[1]{abstract,break,case,catch,const,continue,do,else,elseif,%
      end,export,false,for,function,immutable,import,importall,if,in,%
      macro,module,otherwise,quote,return,switch,true,try,type,typealias,%
      using,while},%
   morekeywords=[2]{LinearAlgebra,GenericLinearAlgebra},%
    morekeywords=[3]{diagm,zeros,ones,cos,range,LinRange,sort,length,setprecision,big,real,convert,isreal,eigvals},
        morekeywords=[4]{Float64,BigFloat,Bool,Integer,Array},
   sensitive=true,%
   alsoother={$},%
   morecomment=[l]\#,%
   morecomment=[n]{\#=}{=\#},%
   morestring=[s]{"}{"},%
   morestring=[m]{'}{'},%
}[keywords,comments,strings]%

\begin{document}
\title{A Matrix-Less Method to Approximate the Spectrum and the Spectral Function of Toeplitz Matrices with Complex Eigenvalues}
\author{Sven-Erik Ekstr\"om$^\dag$\\ {\small\texttt{see@2pi.se}} \and Paris Vassalos$^\dag$\\ {\small\texttt{pvassal@aueb.gr}}}
\date{
    $^\dag$Athens University of Economics and Business\\
    {\small 76 Patission Str., GR-10434 Athens, Greece}\\[2ex]
    \today
}
\maketitle
\begin{abstract}
It is known that the generating function $f$ of a sequence of Toeplitz matrices $\{T_n(f)\}_n$ may not describe the asymptotic distribution of the eigenvalues of $T_n(f)$ if $f$ is not real. 
In a recent paper, we assume as a working hypothesis that, if the eigenvalues of $T_n(f)$ are real for all $n$, then they admit an asymptotic expansion where the first function $g$ appearing in this expansion is real and describes the asymptotic distribution of the eigenvalues of $T_n(f)$. 
In this paper we extend this idea to Toeplitz matrices with complex eigenvalues. The paper is predominantly a numerical exploration of different typical cases, and presents several avenues of possible future research.
\end{abstract}

\section{Introduction}
\label{sec:introduction}
Given a function $f\in L^1([-\pi,\pi])$, we have the Fourier coefficients $\hat{f}_k$, $k\in \mathbb{Z}$, and the Fourier series of $f$ by
\begin{align}
\hat{f}_k=\frac{1}{2\pi}\int_{-\pi}^{\pi}\!\!f(\theta)\E^{-k\mathbf{i}\theta}\mathrm{d}\theta,\qquad f(\theta)=\!\!\sum_{k=-\infty}^{\infty}\!\!\hat{f}_k\E^{k\mathbf{i}\theta}.\label{eq:introduction:fourier}
\end{align}
The $n\times n$ Toeplitz matrix $T_n(f)$ is said to be generated by $f$, if 
\begin{align}
T_n(f)=\left[\hat{f}_{i-j}\right]_{i,j=1}^n=\left[
\begin{array}{ccccccccc}
\hat{f}_0&\hat{f}_{-1}&\cdots&\hat{f}_{2-n}&\hat{f}_{1-n}\vphantom{\ddots}\\
\hat{f}_1&\hat{f}_0&\ddots&&\hat{f}_{2-n}\\
\vdots&\ddots&\ddots&\ddots&\vdots\\
\hat{f}_{n-2}&&\ddots&\hat{f}_{0}&\hat{f}_{-1}\\
\hat{f}_{n-1}&\hat{f}_{n-2}&\cdots&\hat{f}_1&\hat{f}_0\vphantom{\ddots}\\
\\[-0.6em]
\end{array}
\right],\nonumber
\end{align}
where $\hat{f}_k$ are the Fourier coefficients given in \eqref{eq:introduction:fourier}.

It is known that the generating function $f$, also known as the symbol of the matrix sequence $\{T_n(f)\}_n$,
describes the asymptotic distribution of the singular values of $T_n(f)$; if $f$ is real or if $f\in L^\infty([-\pi,\pi])$ and its essential range has empty interior and does not disconnect the complex plane, then $f$ also describes the asymptotic distribution of the eigenvalues of $f$; see \cite{bottcher991,garoni171,tilli991} for details and~\cite[Section~3.1]{garoni171} for the notion of asymptotic singular value and eigenvalue distribution of a sequence of matrices. We write $\{T_n(f)\}_n\sim_\sigma f$ to indicate that $\{T_n(f)\}_n$ has an asymptotic singular value distribution described by $f$ and $\{T_n(f)\}_n\sim_\lambda f$ to indicate that $\{T_n(f)\}_n$ has an asymptotic eigenvalue distribution described by $f$.
In \cite{ekstrom193} the cases of interest were those in which $\{T_n(f)\}_n\not\sim_\lambda f$ and the eigenvalues of $T_n(f)$ are real for all $n$. We believe that in these cases there exist a real function $g$ such that $\{T_n(f)\}_n\sim_\lambda g$ and the eigenvalues of $T_n(f)$ admit an asymptotic expansion of the same type as considered in previous works; e.g., \cite{ekstrom171,ahmad171,ekstrom181,ekstrom183,ekstrom193,ekstrom184,ekstrom185,batalshchikov192,barrera181,bogoya151,bogoya171,bottcher101}. 

In this paper we extend this notion to the case when the eigenvalues of $T_n(f)$ are complex-valued for all $n$. We then assume that there exist a function $g(\theta)=g^\Re(\theta)+\mathbf{i}g^\Im(\theta)$ that describes the eigenvalue distribution. In the case $\{T_n(f)\}_n\sim_\lambda f$ we have $f=g$. For $\{T_n(f)\}_n\not\sim_\lambda f$ we have $f\neq g$ and, as in the real case~\cite{ekstrom193} where $g\coloneqq c_0$, we have $g=g^\Re+\mathbf{i}g^\Im\coloneqq c_0= c_0^\Re+\mathbf{i}c_0^\Im$.
We therefore formulate the following working hypothesis. 

\begin{wh}
Suppose that the eigenvalues of $T_n(f)$ are complex-valued for all $n$. 
Then, for every integer $\alpha\ge0$, every $n$, and every $j=1,\ldots,n$, the following asymptotic expansion holds:
\begin{align}\label{eq:introduction:hoapp}
\lambda_j(T_n(f))&=g^\Re(\theta_{j,n})+\sum_{k=1}^\alpha c_k^\Re(\theta_{j,n})h^k+E_{j,n,\alpha}^\Re+
\mathbf{i}\left\{g^\Im(\theta_{j,n})+\sum_{k=1}^\alpha c_k^\Im(\theta_{j,n})h^k+E_{j,n,\alpha}^\Im\right\},\nonumber\\
&=\sum_{k=0}^\alpha \left(c_k^\Re(\theta_{j,n})+\mathbf{i}c_k^\Im(\theta_{j,n})\right)h^k+E_{j,n,\alpha},
\end{align}
where:
\begin{itemize}
	\item the eigenvalues $\lambda_j(T_n(f))$ are arranged in a consistent order, as $n$ varies. For a monotone symbol $g^\Re$ or $g^\Im$ this is typically a non-decreasing order, that is, either 
  $\Re\{\lambda_1(T_n(f))\}\le\ldots\le\Re\{\lambda_n(T_n(f))\}$ or $\Im\{\lambda_1(T_n(f))\}\le\ldots\le\Im\{\lambda_n(T_n(f))\}$.
  However, we demonstrate ordering strategies for a case where neither $g^\Re$ nor $g^\Im$ are monotone. Also, subsets of the eigenvalues could be considered;
	\item $\{g^\Re\coloneqq c_0^\Re,c_1^\Re,c_2^\Re,c_3^\Re,\ldots\}$ and $\{g^\Im\coloneqq c_0^\Im,c_1^\Im,c_2^\Im,c_3^\Im,\ldots\}$ are sequences of functions from $(0,\pi)$ to $\mathbb R$ which depends only on $f$;
	\item $h=\frac{1}{n+1}$ and $\theta_{j,n}=\frac{j\pi}{n+1}=j\pi h$;
	\item $E_{j,n,\alpha}=E_{j,n,\alpha}^\Re+\mathbf{i}E_{j,n,\alpha}^\Im=O(h^{\alpha+1})$ is the remainder (the error), which satisfies the inequality $|E_{j,n,\alpha}|\le C_\alpha h^{\alpha+1}$ for some constant $C_\alpha$ depending only on $\alpha,f$.
\end{itemize}
\end{wh}
\begin{rmrk}
For notational purposes, since we typically have two different expansions for $g^\Re$ and $g^\Im$ we introduce $\xi_{j,n}^\Re$ and $\xi_{j,n}^\Im$ that denote two ``perfect'' sampling grids, typically not equispaced, such that $\lambda_j(T_n(f))=g^\Re(\xi_{j,n}^\Re)+\mathbf{i}g^\Im(\xi_{j,n}^\Im)$ for $j=1,\ldots,n$; asymptotic expansions of such grids are discussed for matrix sequences $\{T_n(f)\}_n$ with real eigenvalues in~\cite{ekstrom191}. Here, the grids $\xi_{j,n}^\Re$ and $\xi_{j,n}^\Im$ are only used for visualization.
\end{rmrk}

The paper is organized as follows. In Section~\ref{sec:motivation} we present four representative examples for testing the working hypothesis. In Section~\ref{sec:describing} we describe the numerical approach in Algorithm~\ref{algo:1} for approximating $c_k^\Re(\theta_{j,n_0})$ and $c_k^\Im(\theta_{j,n_0})$ in the working hypothesis. For completeness we include Algorithm~\ref{algo:2}, previously given in~\cite{ekstrom191}, for computing the Fourier coefficients. In Section~\ref{sec:numerical} we present numerical results for the four previously defined examples.  Finally, in the conclusion we discuss the presented results and possible future research avenues.

\section{Motivation and illustrative examples}
\label{sec:motivation}
In this section we present four examples in support of our working hypothesis.
We also briefly mention, as in~\cite{ekstrom193}, the fact that standard double precision eigenvalue solvers (such as LAPACK, \texttt{eig} in \textsc{Matlab}, and \texttt{eigvals} in \textsc{Julia}~\cite{bezanson171}) fail to give accurate eigenvalues of certain matrices $T_n(f)$; see, e.g., \cite{beam931,trefethen051}.
High-precision computations, by using packages such as \textsc{GenericLinearAlgebra.jl}~\cite{noack191} in \textsc{Julia}, can compute the true eigenvalues, but they are very expensive from the computational point of view. 
Therefore, approximating $g$ on the grid $\theta_{j,n}$ and using matrix-less methods~\cite{ekstrom183,ekstrom193} to compute the spectrum of $T_n(f)$ can be computationally very advantageous. Also, the presented approaches can be a valuable tool for the analysis of the spectra of non-normal Toeplitz matrices having complex eigenvalues; see~\cite{ekstrom193} for real eigenvalues. 

Following is a short summary of the four examples we consider in the current paper:
\begin{itemize}
\item \textbf{Example~\ref{exmp:1}:} $T_n(f)$ is non-symmetric complex-valued tridiagonal, $g\neq f$, and $g$ is known. Eigenvalues $\lambda_j(T_n(f))$ are known explicitly. Both symbols $g^\Re(\theta)$ and $g^\Im(\theta)$ are monotone;
\item \textbf{Example~\ref{exmp:2}:} $T_n(f)$ is complex symmetric pentadiagonal, and $g=f$. Eigenvalues $\lambda_j(T_n(f))$ are not known explicitly. The symbol $g^\Re(\theta)$ is non-monotone and $g^\Im(\theta)$ is monotone;
\item \textbf{Example~\ref{exmp:3}:} $T_n(f)$ is complex symmetric heptadiagonal, and $g=f$. Eigenvalues $\lambda_j(T_n(f))$ are not known explicitly. Both symbols $g^\Re(\theta)$ and $g^\Im(\theta)$ are non-monotone;
\item \textbf{Example~\ref{exmp:4}:} $T_n(f)$ is non-symmetric real-valued banded (shifted non-symmetric real-valued pentadiagonal), $g\neq f$, and $g$ is not known. Eigenvalues $\lambda_j(T_n(f))$ are not known explicitly. The symbol $g^\Re(\theta)$ is non-monotone and $g^\Im(\theta)$ is monotone.
\end{itemize}
In Section~\ref{sec:numerical} we perform numerical experiments in  Examples~\ref{exmp:5}--\ref{exmp:8} (corresponding to Examples~\ref{exmp:1}--\ref{exmp:4}), supporting the working hypothesis.

\begin{exmp}
\label{exmp:1}
In \cite[Examples 1 and 5]{ekstrom193} we studied the symbol 
$f(\theta)=-\E^{\mathbf{i}\theta}+2-2\E^{-\mathbf{i}\theta}$,
which generates a Toeplitz matrix $T_n(f)$ that has a real spectrum described by the symbol
$
g(\theta)=2-2\sqrt{2}\cos(\theta)
$.
The exact eigenvalues of $T_n(f)$ are given by $\lambda_j(T_n(f))=\lambda_j(T_n(g))=g\left(\theta_{j,n}\right)$,
where sampling grid is
\begin{align}
\theta_{j,n}=\frac{j\pi}{n+1},\qquad j=1,\ldots,n.\label{eq:exmp1:tau_grid}
\end{align}
Now, instead consider the symbol
\begin{align}
f(\theta)&=-\E^{\mathbf{i}\theta}+2+(-2+\mathbf{i})\E^{-\mathbf{i}\theta},\label{eq:exmp1:symbol}
\end{align}
and the corresponding symbol $g$ that describes the complex-valued spectrum of $T_n(f)$ (see (4) in \cite{ekstrom193}) is
\begin{align}
g(\theta)&=2+2\sqrt{-1}\sqrt{-2+\mathbf{i}}\cos(\theta)\nonumber\\
&=\underbrace{2}_{=\hat{g}_0}+2\underbrace{\sqrt[4]{5}\left(-\cos\left(\frac{\tan^{-1}(1/2)}{2}\right)+\mathbf{i}\sin\left(\frac{\tan^{-1}(1/2)}{2}\right)\right)}_{=\hat{g}_{\pm 1}}\cos(\theta)\nonumber\\
&=\underbrace{\underbrace{2}_{=\hat{g}_0^\Re}+2\underbrace{\left(-\sqrt[4]{5}\cos\left(\frac{\tan^{-1}(1/2)}{2}\right)\right)}_{=\hat{g}_{\pm 1}^\Re}\cos(\theta)}_{=g^\Re(\theta)}+\mathbf{i}\underbrace{2\underbrace{\sqrt[4]{5}\sin\left(\frac{\tan^{-1}(1/2)}{2}\right)}_{=\hat{g}_{\pm 1}^\Im}\cos(\theta)}_{=g^\Im(\theta)}\nonumber\\
&\approx 2-2.91069338\cos(\theta)+0.6871215\mathbf{i}\cos(\theta).\label{eq:exmp1:symbolg}
\end{align}
The eigenvalues of $T_n(f)$ (and $T_n(g)$) are given by \eqref{eq:exmp1:symbolg} with the grid \eqref{eq:exmp1:tau_grid}, where
\begin{align}
T_n(f)&=\left[
\begin{array}{ccccc}
2&-2+\mathbf{i}\\
-1&2&-2+\mathbf{i}\\
&\ddots&\ddots&\ddots\\
&&-1&2&-2+\mathbf{i}\\
&&&-1&2
\end{array}
\right]=
\left[
\begin{array}{ccccc}
2&-2\\
-1&2&-2\\
&\ddots&\ddots&\ddots\\
&&-1&2&-2\\
&&&-1&2
\end{array}
\right]
+\mathbf{i}
\left[
\begin{array}{ccccc}
0&1\\
&0&1\\
&\ddots&\ddots&\ddots\\
&&&0&1\\
&&&&0
\end{array}
\right]
,\nonumber
\end{align}
and
\begin{align}
T_n(g)&=\left[
\begin{array}{ccccc}
2&\hat{g}_1\\
\hat{g}_1&2&\hat{g}_1\\
&\ddots&\ddots&\ddots\\
&&\hat{g}_1&2&\hat{g}_1\\
&&&\hat{g}_1&2
\end{array}
\right]
=\left[
\begin{array}{ccccc}
2&\hat{g}_1^\Re\\
\hat{g}_1^\Re&2&\hat{g}_1^\Re\\
&\ddots&\ddots&\ddots\\
&&\hat{g}_1^\Re&2&\hat{g}_1^\Re\\
&&&\hat{g}_1^\Re&2
\end{array}
\right]+\mathbf{i}\left[
\begin{array}{ccccc}
0&\hat{g}_1^\Im\\
\hat{g}_1^\Im&0&\hat{g}_1^\Im\\
&\ddots&\ddots&\ddots\\
&&\hat{g}_1^\Im&0&\hat{g}_1^\Im\\
&&&\hat{g}_1^\Im&0
\end{array}
\right].\nonumber
\end{align}
\begin{figure}[!ht] 
\centering
\includegraphics[width=0.47\textwidth,valign=t]{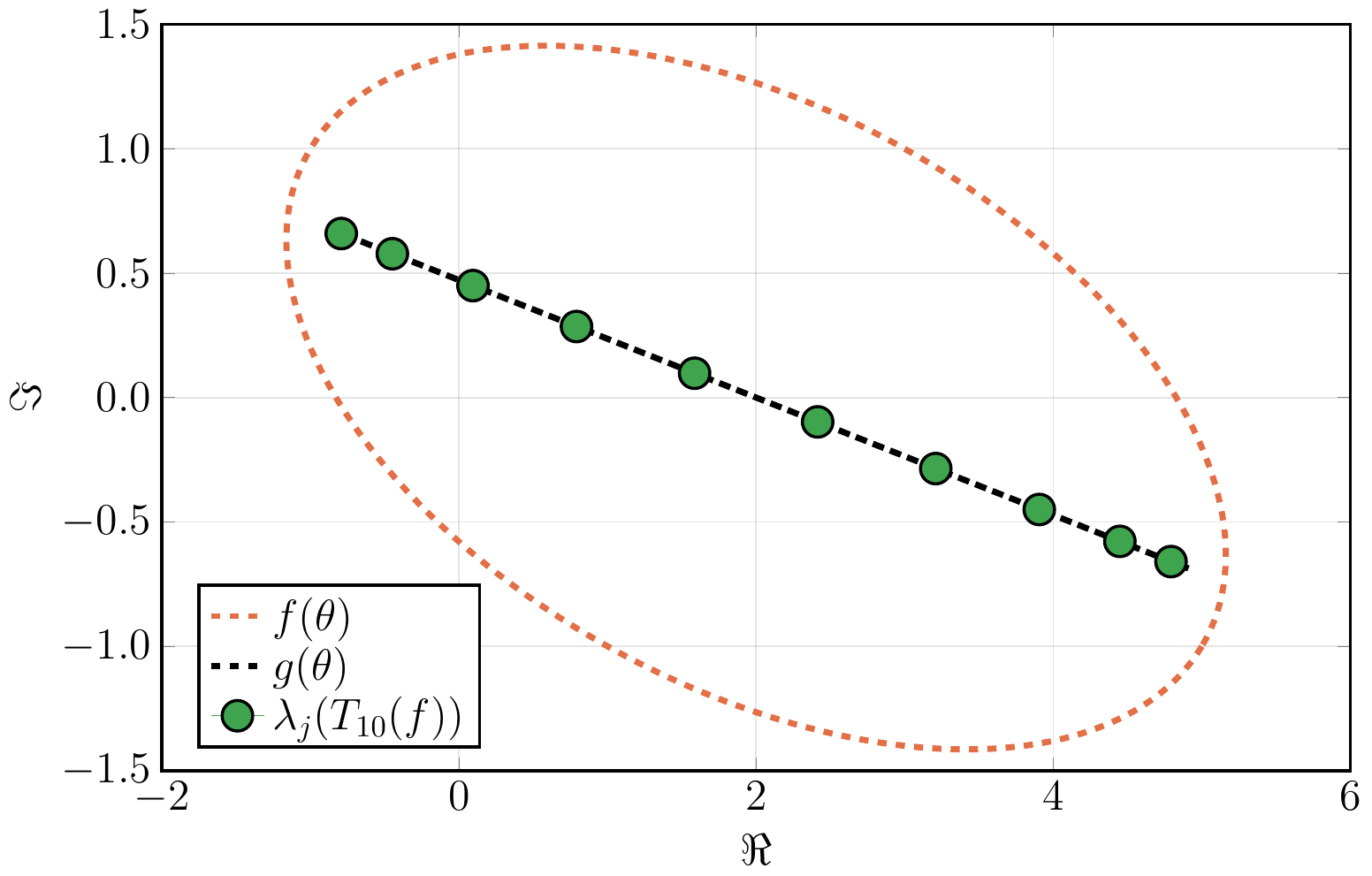}
\includegraphics[width=0.47\textwidth,valign=t]{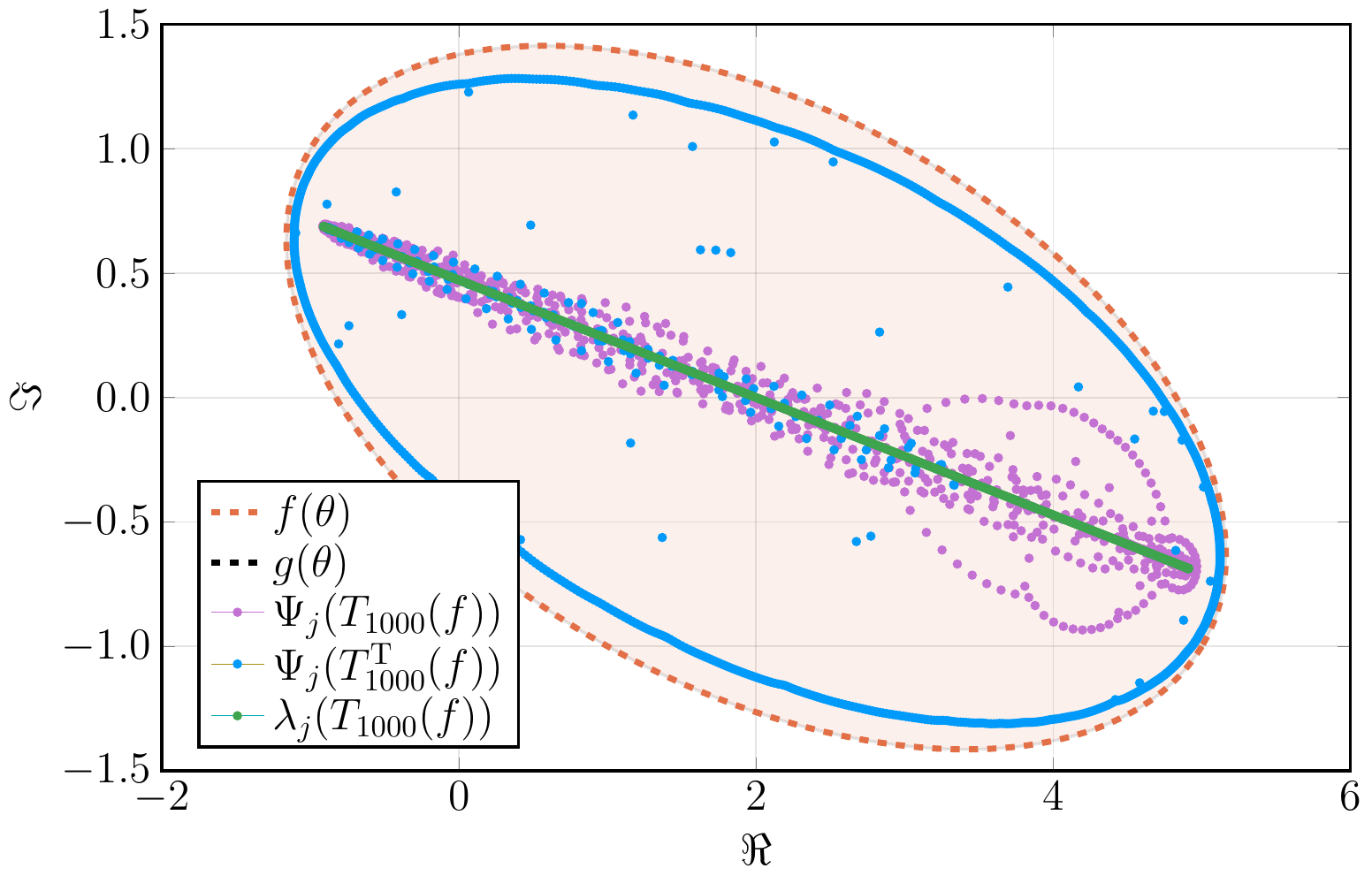}
\caption{[Example~\ref{exmp:1}: Symbol $f(\theta)=-\E^{\mathbf{i}\theta}+2+(-2+\mathbf{i})\E^{-\mathbf{i}\theta}$] Left: Symbols $f(\theta)$ (red dashed line), and $g(\theta)$ (black dashed line), and eigenvalues $\lambda_j(T_n(f))=\lambda_j(T_n(g))$ for $n=10$ (green circles). Right: Symbols $f$ and $g$, numerically computed eigenvalues, for $n=1000$,  $\Psi_j(T_{n}(f))$, $\Psi_j(T_{n}^{\mathrm{T}}(f))$, and $\lambda_j(T_{n}(f))=\lambda_j(T_{n}(g))=g(\theta_{j,n})$. The convex hull of $f$ is indicated in light red.}
\label{fig:exmp:1:symbols}
\end{figure}

In the left panel of Figure~\ref{fig:exmp:1:symbols} we show the functions $f$ (red dashed line) and $g$ (black dashed line) and the eigenvalues $\lambda_j(T_n(f))=\lambda_j(T_n(g))$ (green circles) for $n=10$. 
In the right panel of Figure~\ref{fig:exmp:1:symbols}  we show the functions $f$ (red dashed line), $g$ (black dashed line), although not visible since it is covered, and for $n=1000$ the numerically computed eigenvalues $\Psi_j(T_{n}(f))$ and $\Psi_j(T_{n}^{\mathrm{T}}(f))$. These numerically computed eigenvalues $\Psi_j(A_n)$ are related to the pseudospectrum, discussed for example in \cite{trefethen051,beam931,reichel921} and similarly in~\cite{ekstrom193}.  Furthermore, the true eigenvalues $\lambda_j(T_{n}(f))=\lambda_j(T_{n}(g))=g(\theta_{j,n})$ are shown. The computation of the eigenvalues of $\lambda_j(T_{n}(f))$ require high precision computation, whereas, for $\lambda_j(T_{n}(g))$ standard double precision is sufficient (and the exact expression is given by sampling \eqref{eq:exmp1:symbolg} with the grid \eqref{eq:exmp1:tau_grid}).

In Figure~\ref{fig:exmp:1:spectrum} we present the real (left panel) and imaginary (right panel) part of the spectrum of $T_n(f)$. For $n=10$ we see that the eigenvalues are equispaced samplings of $g^\Re$ and $g^\Im$. We present the eigenvalues on the grid \eqref{eq:exmp1:tau_grid} since both $g^\Re$ and $g^\Im$ are even functions (whereas neither the real nor the imaginary part of $f$ are even).
\begin{figure}[!ht]
\centering
\includegraphics[width=0.47\textwidth,valign=t]{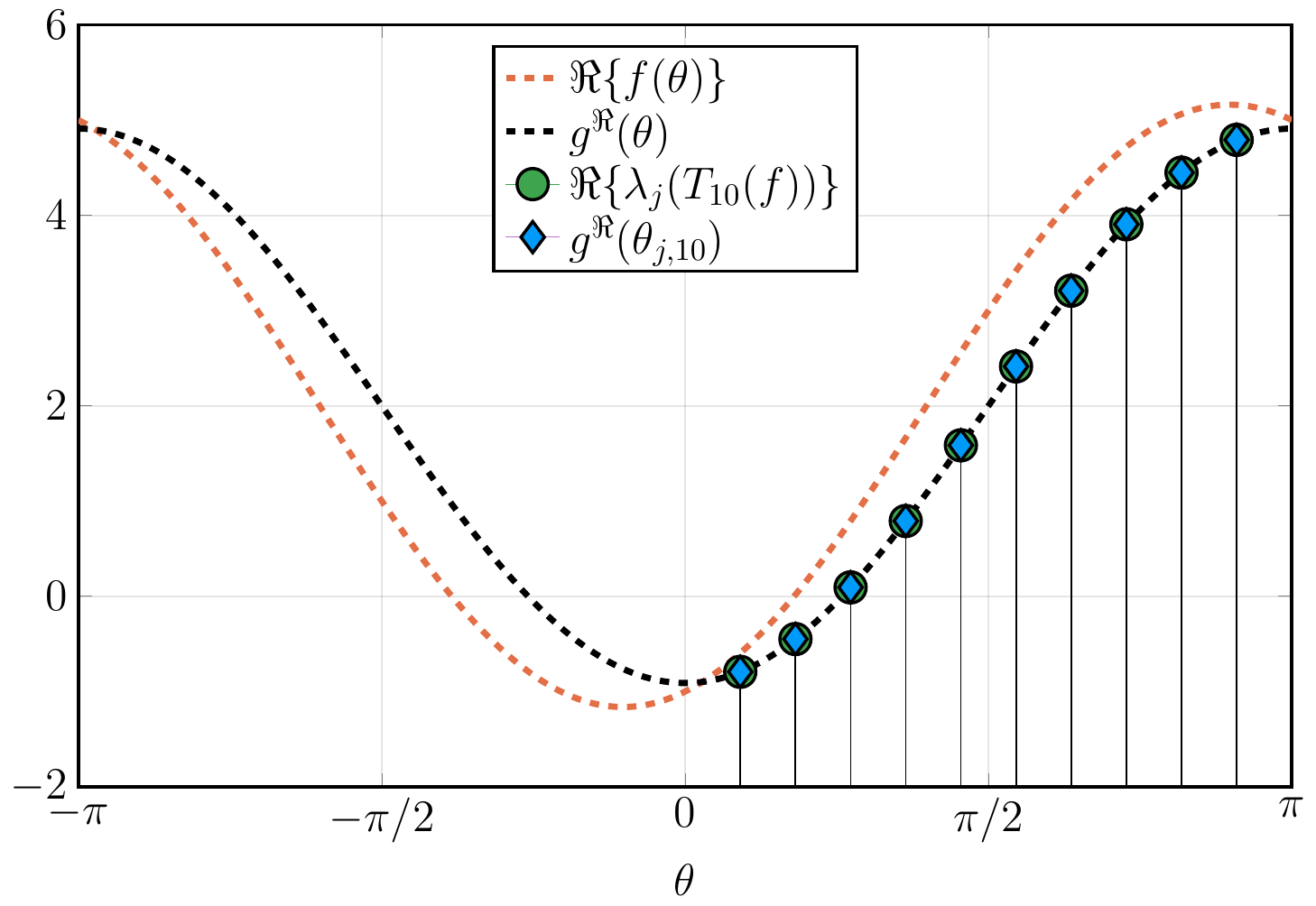}
\includegraphics[width=0.482\textwidth,valign=t]{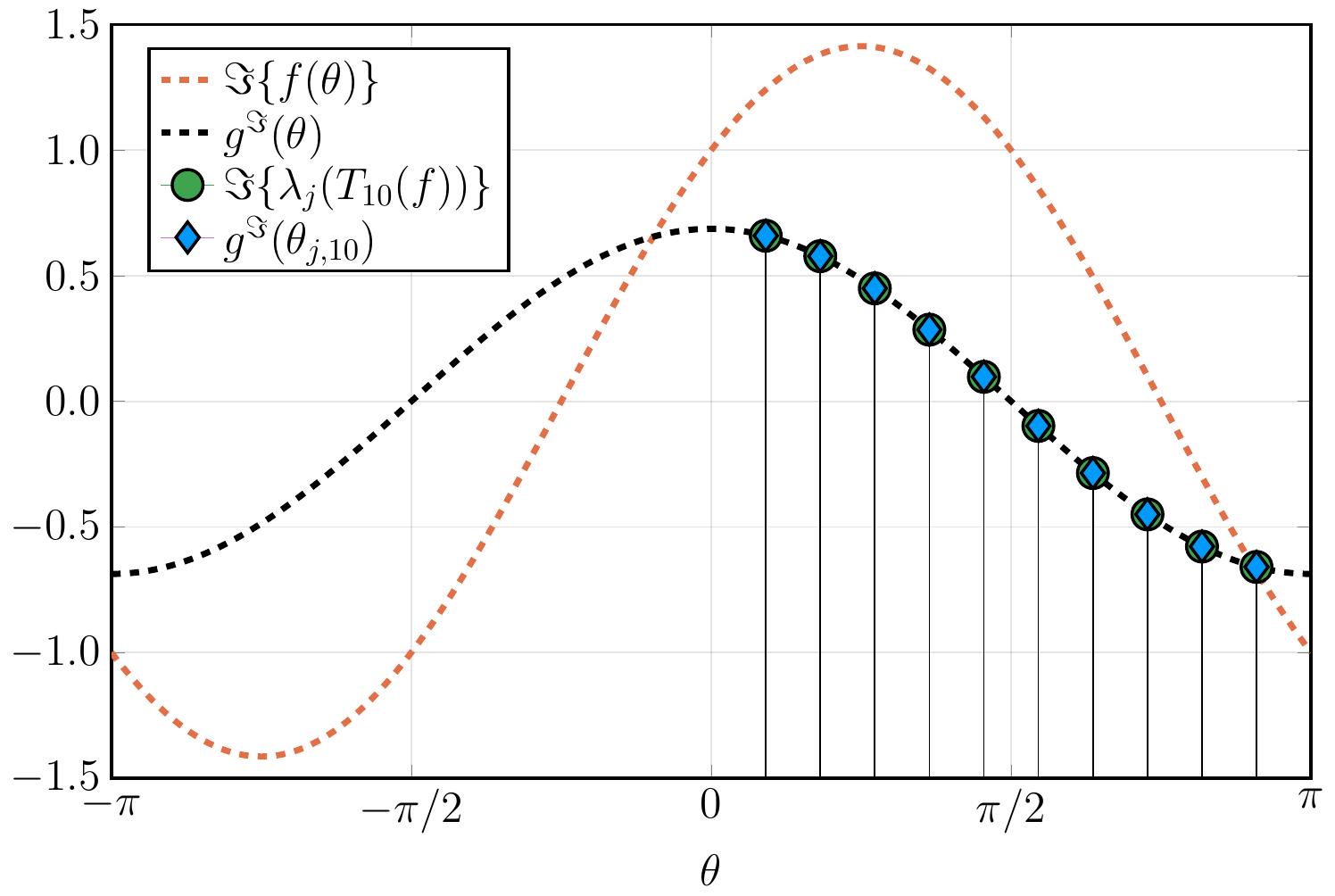}
\caption{[Example~\ref{exmp:1}:  Symbol $f(\theta)=-\E^{\mathbf{i}\theta}+2+(-2+\mathbf{i})\E^{-\mathbf{i}\theta}$] 
Left:
The real part of the symbol $f(\theta)$ (red dashed line), the symbol $g^\Re$ (black dashed line), the real part of the eigenvalues of $T_{10}(f)$, and the sampling of $g^\Re$ with the grid \eqref{eq:exmp1:tau_grid}.
Right: The corresponding imaginary counterparts of the left panel.}
\label{fig:exmp:1:spectrum}
\end{figure}
\end{exmp}

\begin{exmp}
\label{exmp:2}
In this example we construct a function that generates complex-valued matrices with a pentadiagonal real part (generated by a non-monotone function), and a pentadiagonal imaginary part (generated by a monotone function).  The function is chosen to be
\begin{align}
f(\theta)=\underbrace{2\cos(\theta)-2\cos(2\theta)}_{=g^\Re(\theta)}+\mathbf{i}\underbrace{\left(6-8\cos(\theta)+2\cos(2\theta)\right)}_{=g^\Im(\theta)=(2-2\cos(\theta))^2}=g(\theta).\label{eq:exmp2:symbol}
\end{align}
Thus, we have
\begin{align}
T_n(f)&=T_n(g^\Re)+\mathbf{i}T_n(g^\Im)
=\left[
\begin{array}{rrrrrrrrrrr}
0&1&-1\\
1&0&1&-1\\
-1&1&0&1&-1\\
&\ddots&\ddots&\ddots&\ddots&\ddots\\
&&-1&1&0&1&-1\\
&&&-1&1&0&1\\
&&&&-1&1&0
\end{array}
\right]+\mathbf{i}\left[
\begin{array}{rrrrrrrrr}
6&-4&1\\
-4&6&-4&1\\
1&-4&6&-4&1\\
&\ddots&\ddots&\ddots&\ddots&\ddots\\
&&1&-4&6&-4&1\\
&&&1&-4&6&-4\\
&&&&1&-4&6
\end{array}
\right].\nonumber
\end{align}
The spectra of the generated Toeplitz matrices $T_n(f)$ are complex-valued.
\begin{rmrk}
We note that in this example we have $f=g$. And as seen in Figure~\ref{fig:exmp:2:symbols} there are no numerical issues in computing the spectrum of $T_n(f)$ using double precision, and as expected the spectrum converges towards the function $f$ as $n$ increases.
\end{rmrk}
\begin{figure}[!ht] 
\centering
\includegraphics[width=0.47\textwidth,valign=t]{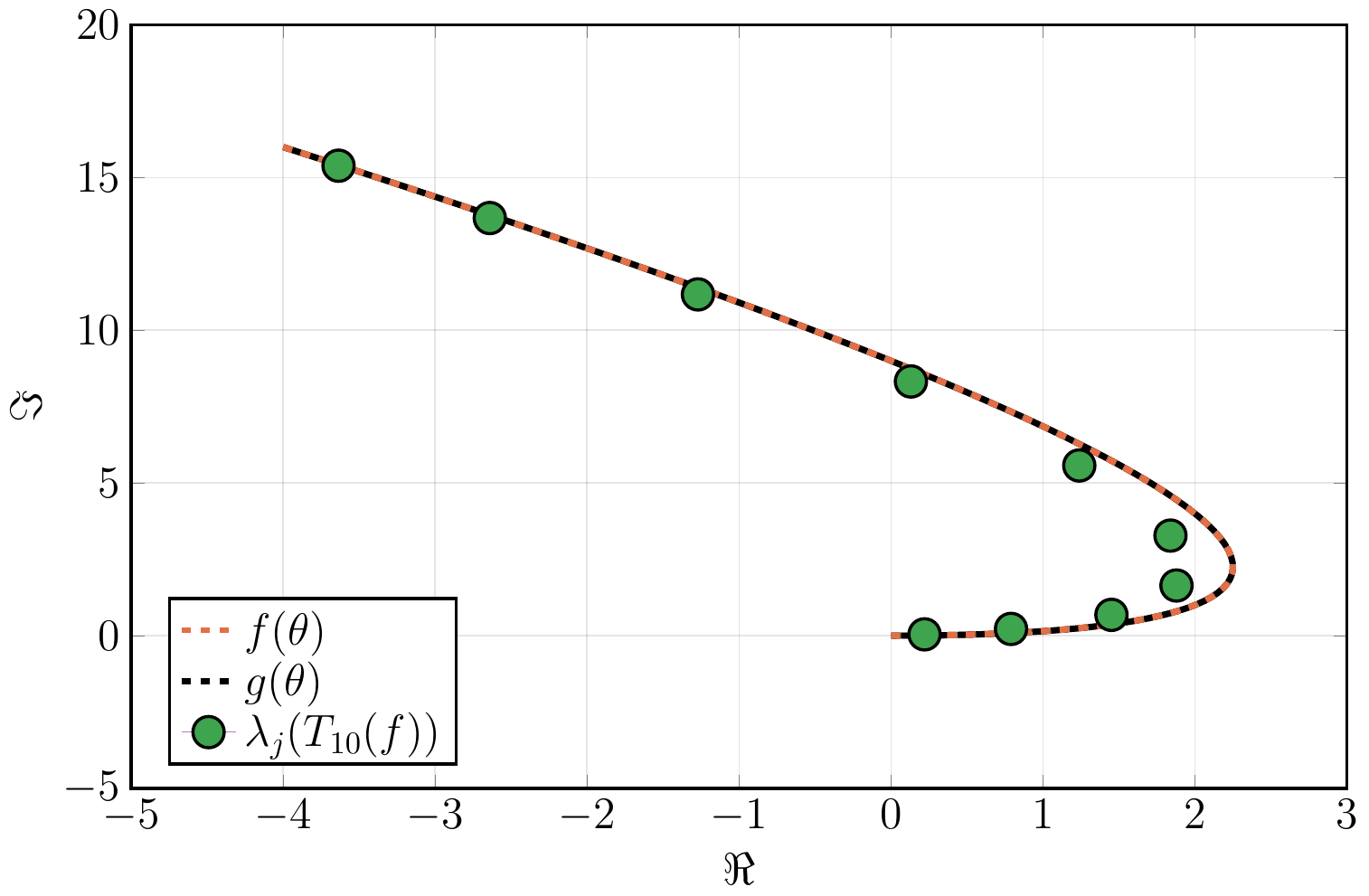}
\includegraphics[width=0.47\textwidth,valign=t]{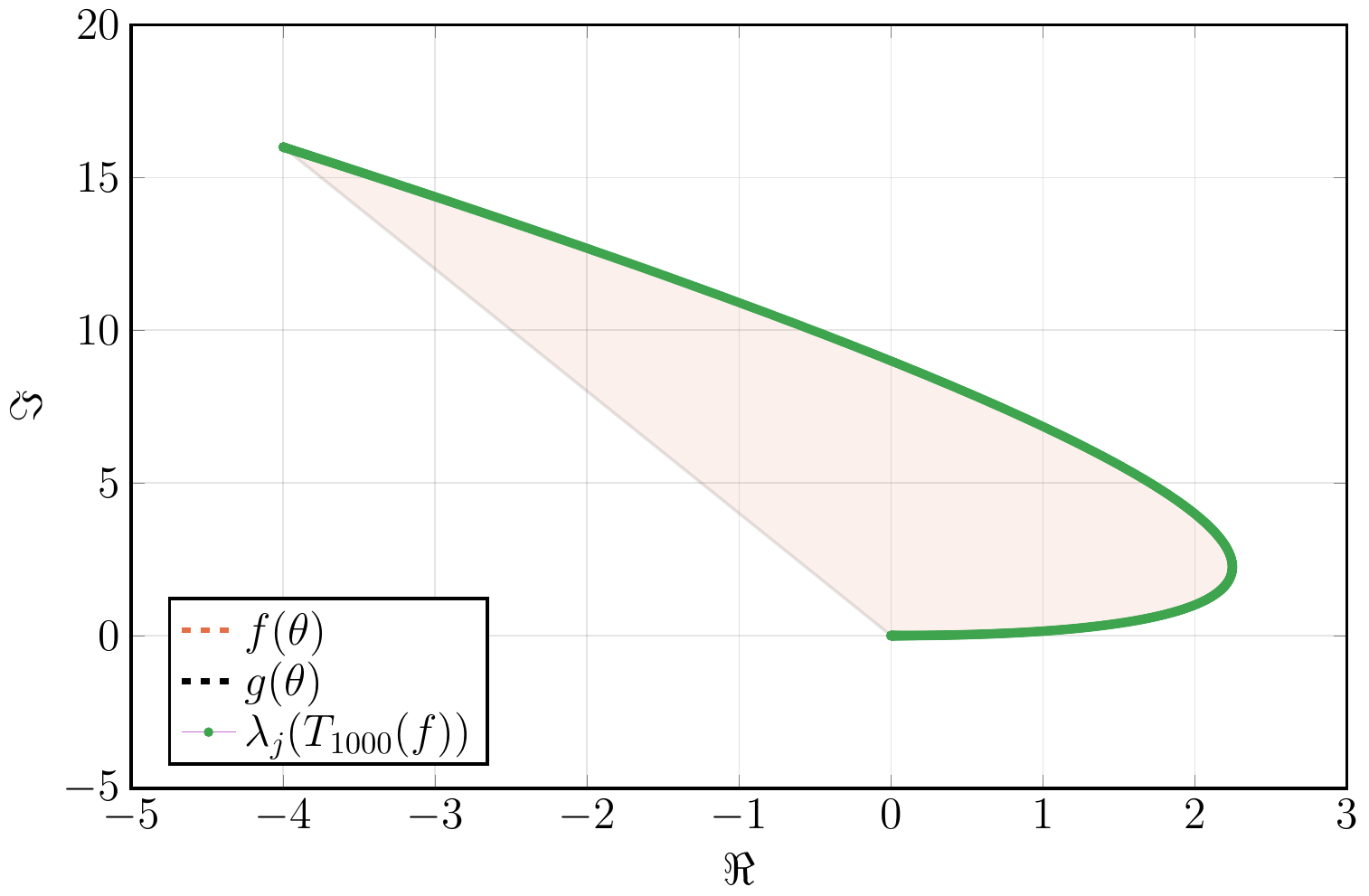} 
\caption{[Example~\ref{exmp:2}: Symbol $f(\theta)=2\cos(\theta)-2\cos(2\theta)+\mathbf{i}(2-2\cos(\theta))^2$] Left: Symbols $f=g$ (red and black dashed lines) and $\lambda_j(T_{10}(f))$ (green circles). 
Right: Convex hull of $f$ (light red) and eigenvalues $\lambda_j(T_{1000}(f))$.}
\label{fig:exmp:2:symbols}
\end{figure}
\begin{figure}[!ht]
\centering
\includegraphics[width=0.46\textwidth,valign=t]{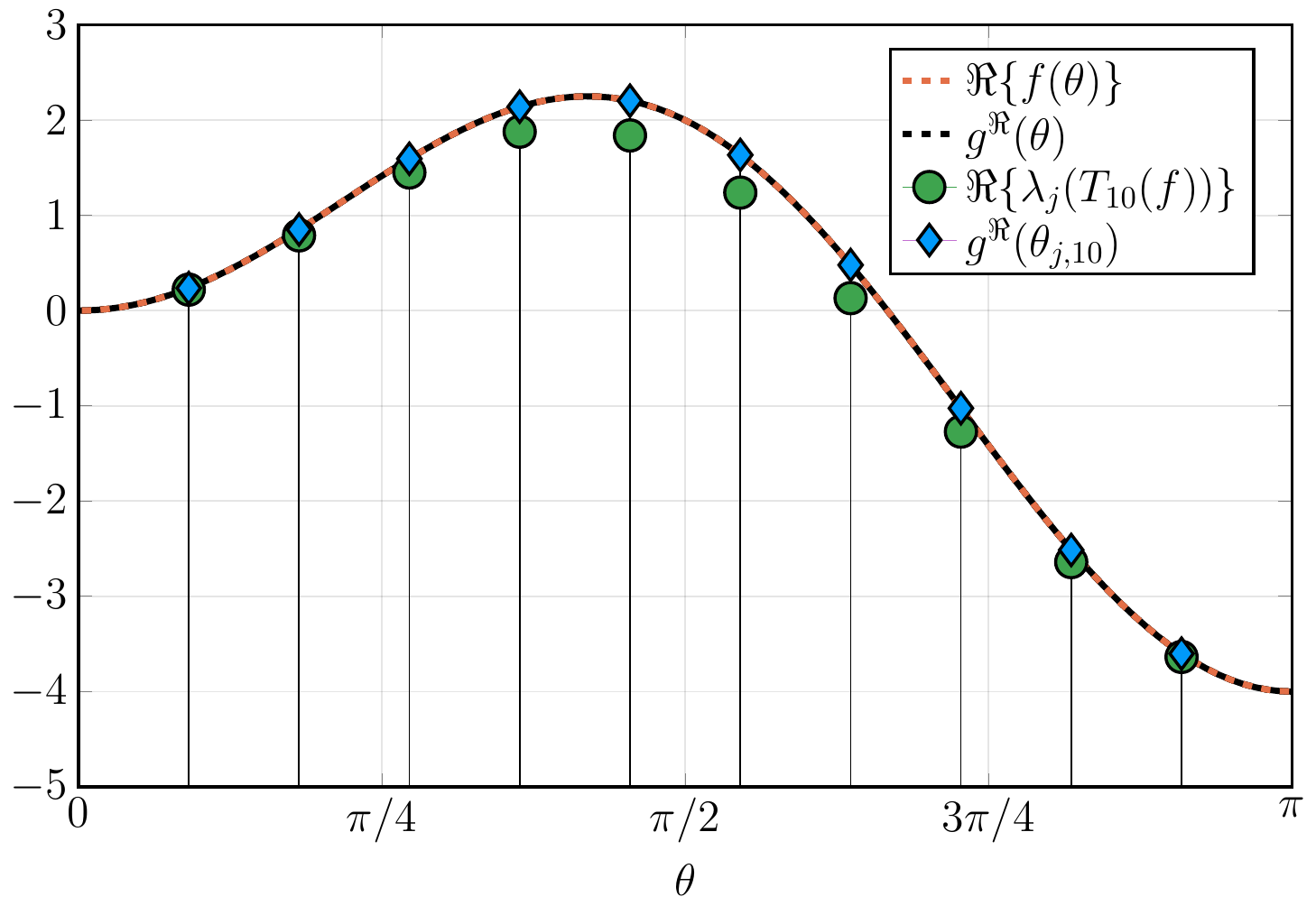}
\includegraphics[width=0.46\textwidth,valign=t]{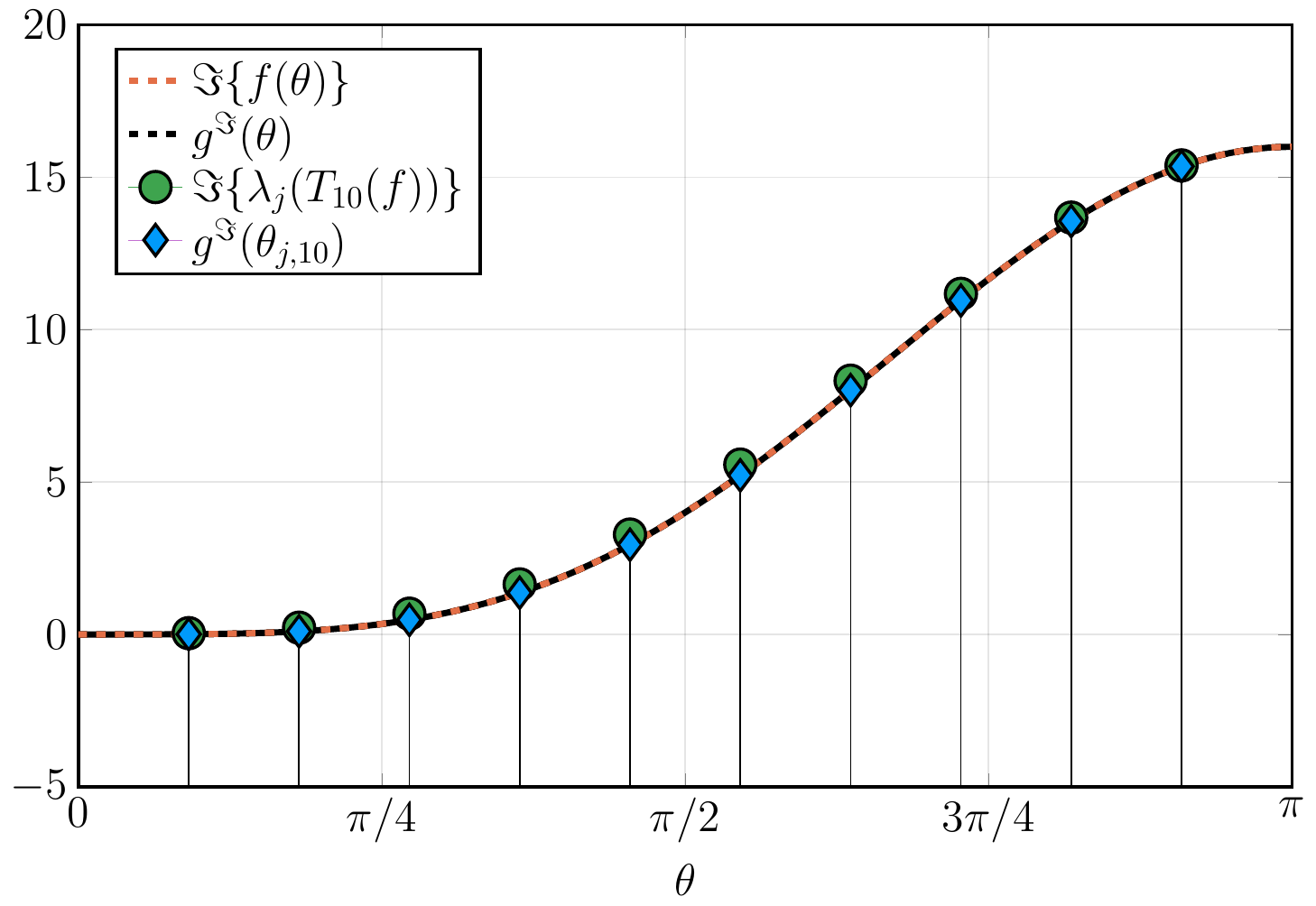}

\includegraphics[width=0.46\textwidth,valign=t]{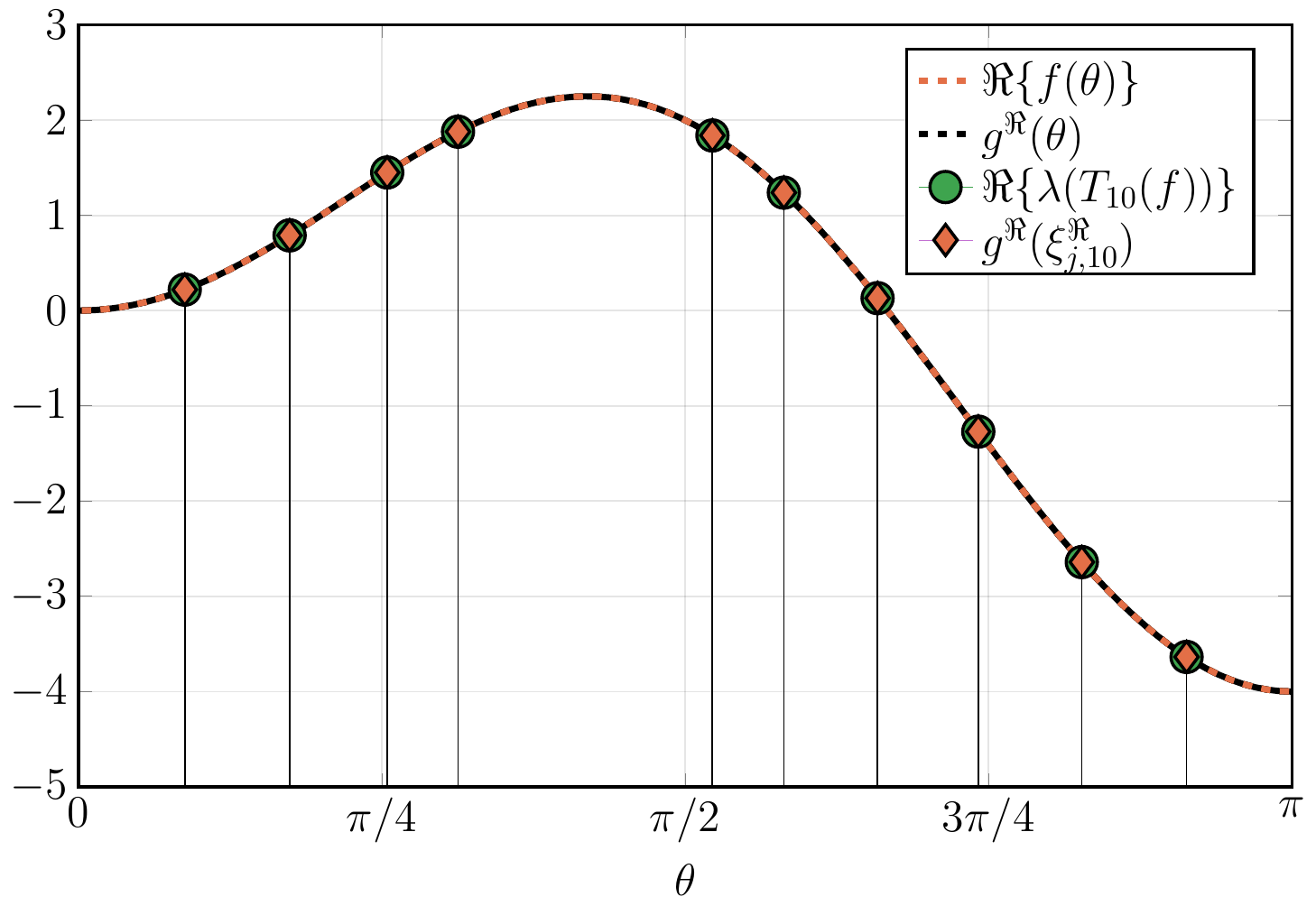}
\includegraphics[width=0.46\textwidth,valign=t]{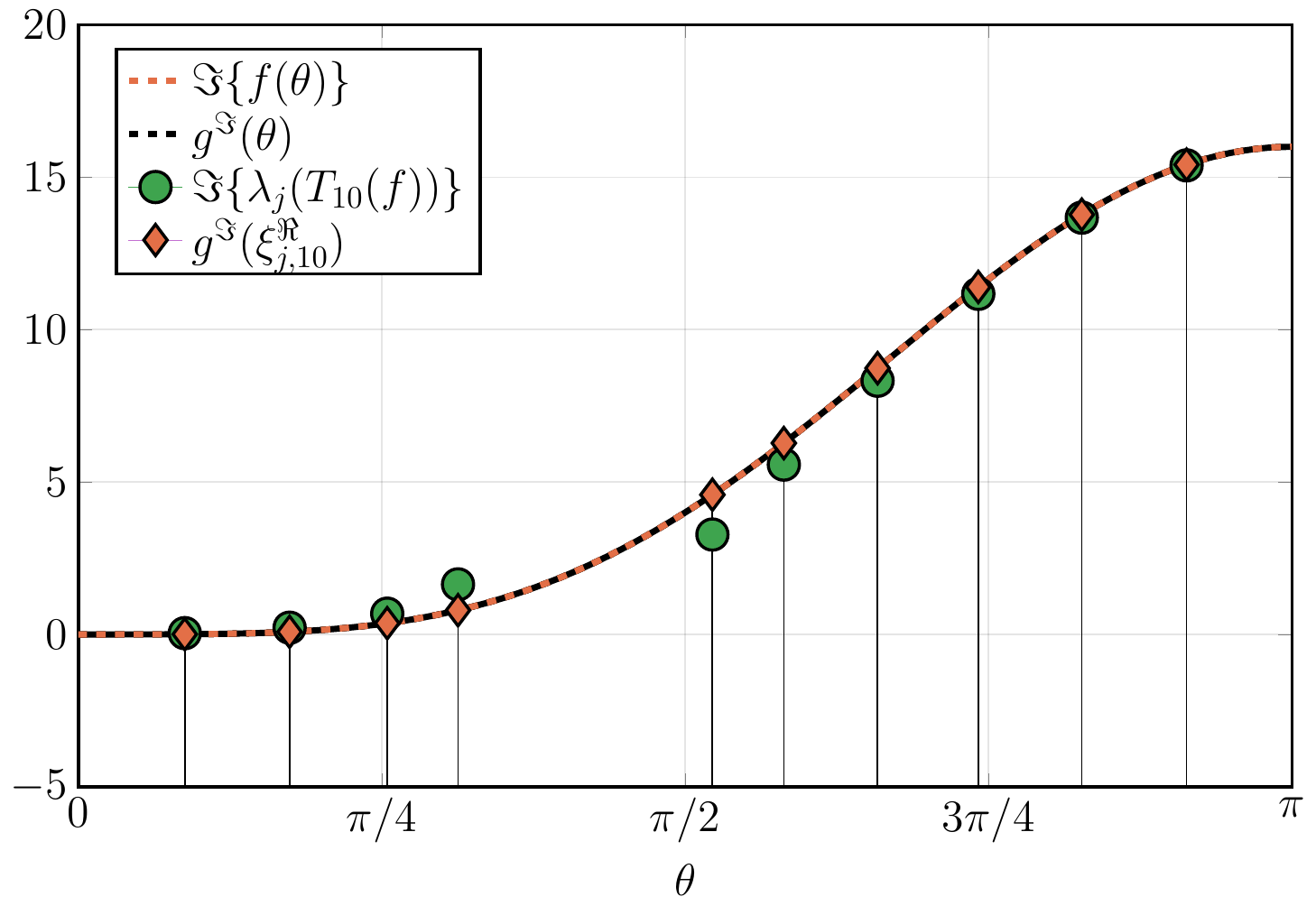}

\includegraphics[width=0.46\textwidth,valign=t]{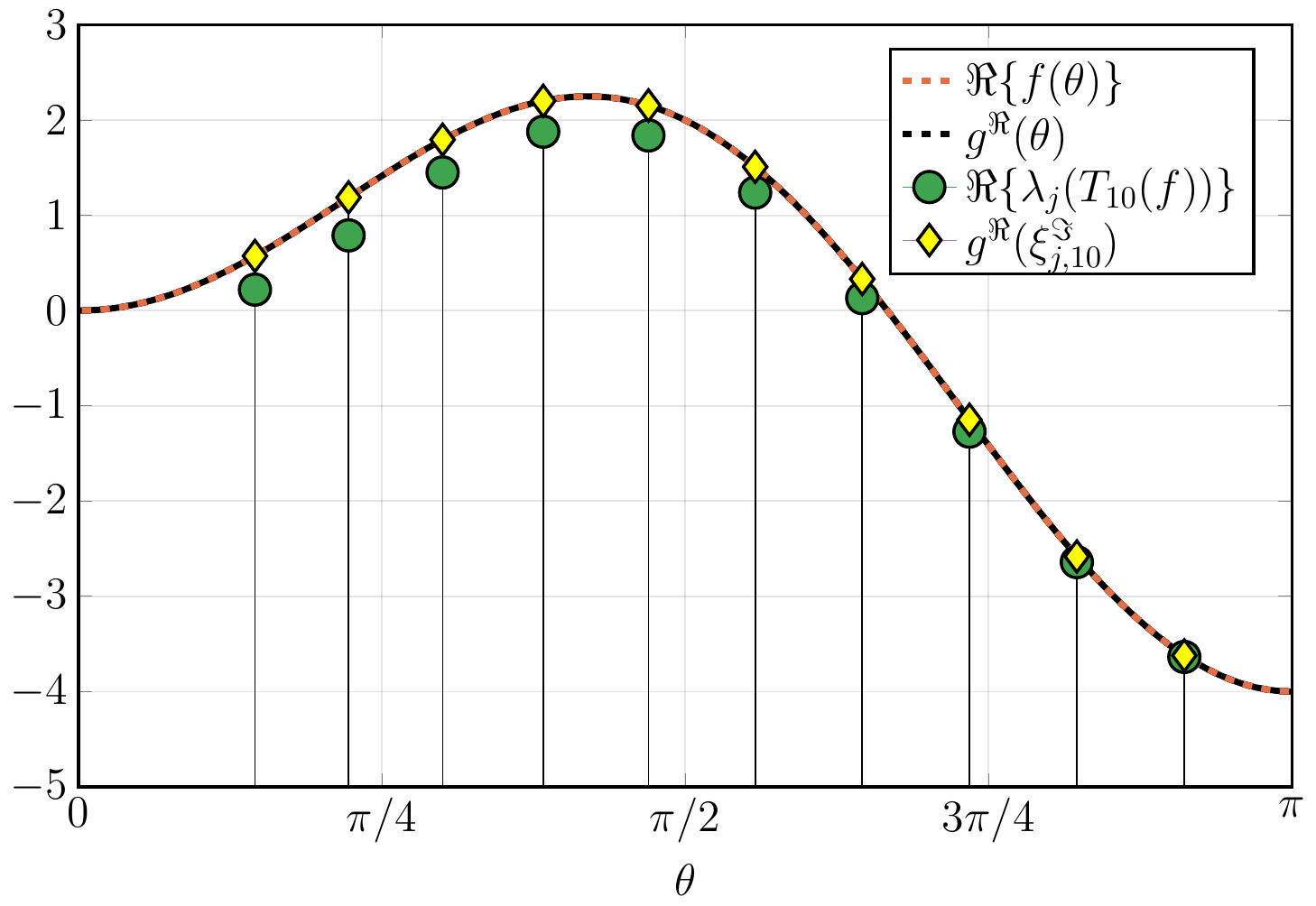}
\includegraphics[width=0.46\textwidth,valign=t]{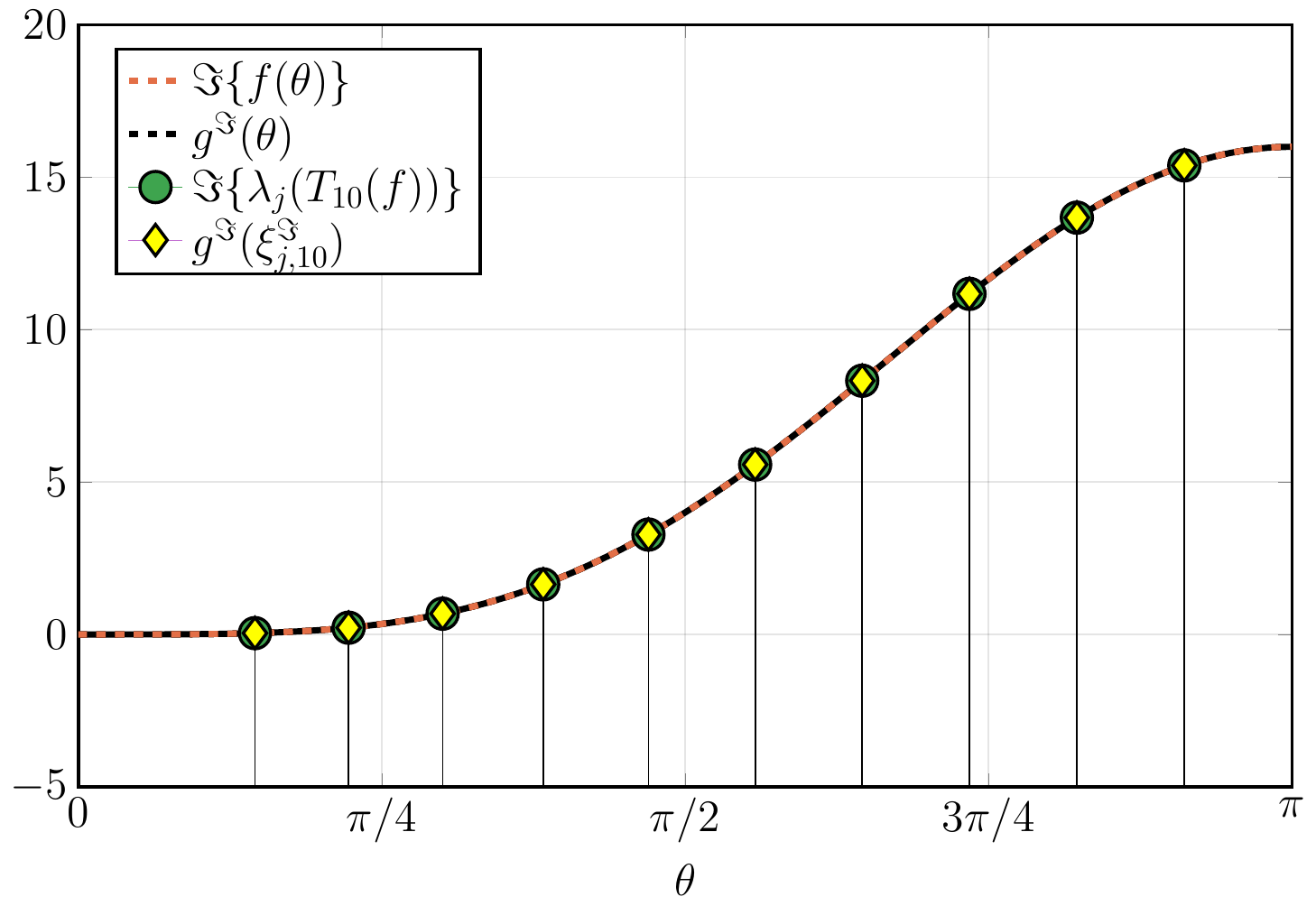}
\caption{[Example~\ref{exmp:2}:  Symbol $f(\theta)=2\cos(\theta)-2\cos(2\theta)+\mathbf{i}(2-2\cos(\theta))^2$]  Left: Real part of the function $f$ (red dashed line) (and $g^\Re=\Re\{f\}$ (black dashed line)). Eigenvalues for $n=10$ $\lambda_j(T_{n}(f))$, and sampling grids $\theta_{j,n}$ (top), $\xi_{j,n}^\Re$ (middle), and $\xi_{j,n}^\Im$ (bottom). Right: Imaginary counterparts of the left panels.}
\label{fig:exmp:2:spectrum}
\end{figure}

\noindent In Figure~\ref{fig:exmp:2:spectrum} we present the spectrum of $T_n(f)$ for $n=10$, presented with different sampling grids $\theta_{j,n}$ (top), $\xi_{j,n}^\Re$ (middle), and $\xi_{j,n}^\Im$ (bottom). The left panels concern the real part of the spectrum of $T_n(f)$ and the right panels the imaginary part.

The top panels of Figure~\ref{fig:exmp:2:spectrum} show that we get an error in the eigenvalue approximations when sampling the function $g$ using the grid $\theta_{j,10}$, defined in  \eqref{eq:exmp1:tau_grid}, both in the real (left panel) and the imaginary (right panel) parts of the spectrum.
In the middle panels of Figure~\ref{fig:exmp:2:spectrum} we use the perfect grid $\xi_{j,10}^\Re$ for sampling $g$, that is, we get a perfect of the real part of the spectrum. However, the imaginary part of the spectrum is not exact.
The bottom panels of Figure~\ref{fig:exmp:2:spectrum} show the results when sampling $g$ with the perfect grid $\xi_{j,10}^\Im$. The left panel shows the erroneous approximated real part, and the right panel shows the exact imaginary part, of the spectrum of $T_{10}(f)$.
\end{exmp}

\begin{exmp}
\label{exmp:3}
In this example we construct the following complex-valued function,
\begin{align}
f(\theta)=\underbrace{2\cos(\theta)-2\cos(2\theta)}_{=g^\Re(\theta)}+\mathbf{i}\underbrace{\left(2\cos(2\theta)-2\cos(3\theta)\right)}_{=g^\Im(\theta)}=g(\theta),\label{eq:exmp3:symbol}
\end{align}
which generates matrices with complex spectra.
We have
\begin{align}
T_n(f)&=T_n(g^\Re)+\mathbf{i}T_n(g^\Im)\nonumber\\
&=\left[
\begin{array}{rrrrrrrrr}
0&1&-1\\
1&0&1&-1\\
-1&1&0&1&-1\\
&\ddots&\ddots&\ddots&\ddots&\ddots\\
&&-1&1&0&1&-1\\
&&&-1&1&0&1\\
&&&&-1&1&0
\end{array}
\right]+\mathbf{i}\left[
\begin{array}{rrrrrrrrrr}
0&0&1&-1\\
0&0&0&1&-1\\
1&0&0&0&1&-1\\
-1&1&0&0&0&1&-1\\
&\ddots&\ddots&\ddots&\ddots&\ddots&\ddots&\ddots\\
&&-1&1&0&0&0&1&-1\\
&&&-1&1&0&0&0&1\\
&&&&-1&1&0&0&0\\
&&&&&-1&1&0&0\\
\end{array}
\right].\nonumber
\end{align}
We note that both $g^\Re$ and $g^\Im$ are non-monotone symbols. Thus, a consistent ordering of the eigenvalues $\lambda_j(T_n(f))$, for varying $n$, can not be done the standard fashion. Instead, we employ the following ordering strategy
\begin{enumerate}
\item Compute the eigenvalues $\tilde{\lambda}_j$, for $j=1,\ldots,n$, of $T_n(f)$  using an appropriate eigenvalue solver. Store the eigenvalues (in any order) in a vector $\tilde{\Lambda}=[\tilde{\lambda}_1,\ldots, \tilde{\lambda}_n]$.
\item Choose the eigenvalue in $\tilde{\Lambda}$  with the smallest absolute value, $\min_j|\tilde{\lambda}_j|$, as the first eigenvalue $\lambda_1$ in a vector $\Lambda$ which will contain the ordered eigenvalues. Remove the corresponding eigenvalue from the vector $\tilde{\Lambda}$.
\item For each eigenvalue $\lambda_j$ with $j>1$ choose $\lambda_{j}$ from $\tilde{\Lambda}$ as the eigenvalue with smallest distance from $\lambda_{j-1}$. After each choice $\lambda_j$ is added to $\Lambda$, remove the corresponding eigenvalue from $\tilde{\Lambda}$.
\item The resulting vector $\Lambda$ with the ordered eigenvalues $\lambda_1,\lambda_2,\ldots,\lambda_n$ will be consistent for any $n$.
\end{enumerate}
Alternatively we could in this case define the ordering by minimizing $\left|\lambda_{j}(T_n(f))-f(\theta_{j,n})\right|$, for $j=1,\ldots,n$, which would give the same result as the scheme defined above.
\begin{rmrk}
This type of ``customized'' ordering would be necessary for many complex-valued spectra, to get consistent ordering as $n$ varies. Also a splitting of the spectrum into multiple parts, described by several ``sub symbols'' is a viable path. For spectra with complex conjugate pairs, one can obviously focus on one branch and reconstruct the rest of the eigenvalues. Further research in this direction is warranted, to devise new standardized and heuristic approaches to consistent eigenvalue ordering.
\end{rmrk}

\noindent In Figures~\ref{fig:exmp:3:symbols} and \ref{fig:exmp:3:spectrum} we present the same information as in Figures~\ref{fig:exmp:2:symbols} and \ref{fig:exmp:2:spectrum} in Example~\ref{exmp:2}, but for \eqref{eq:exmp3:symbol}.
\begin{figure}[!ht] 
\centering
\includegraphics[width=0.47\textwidth,valign=t]{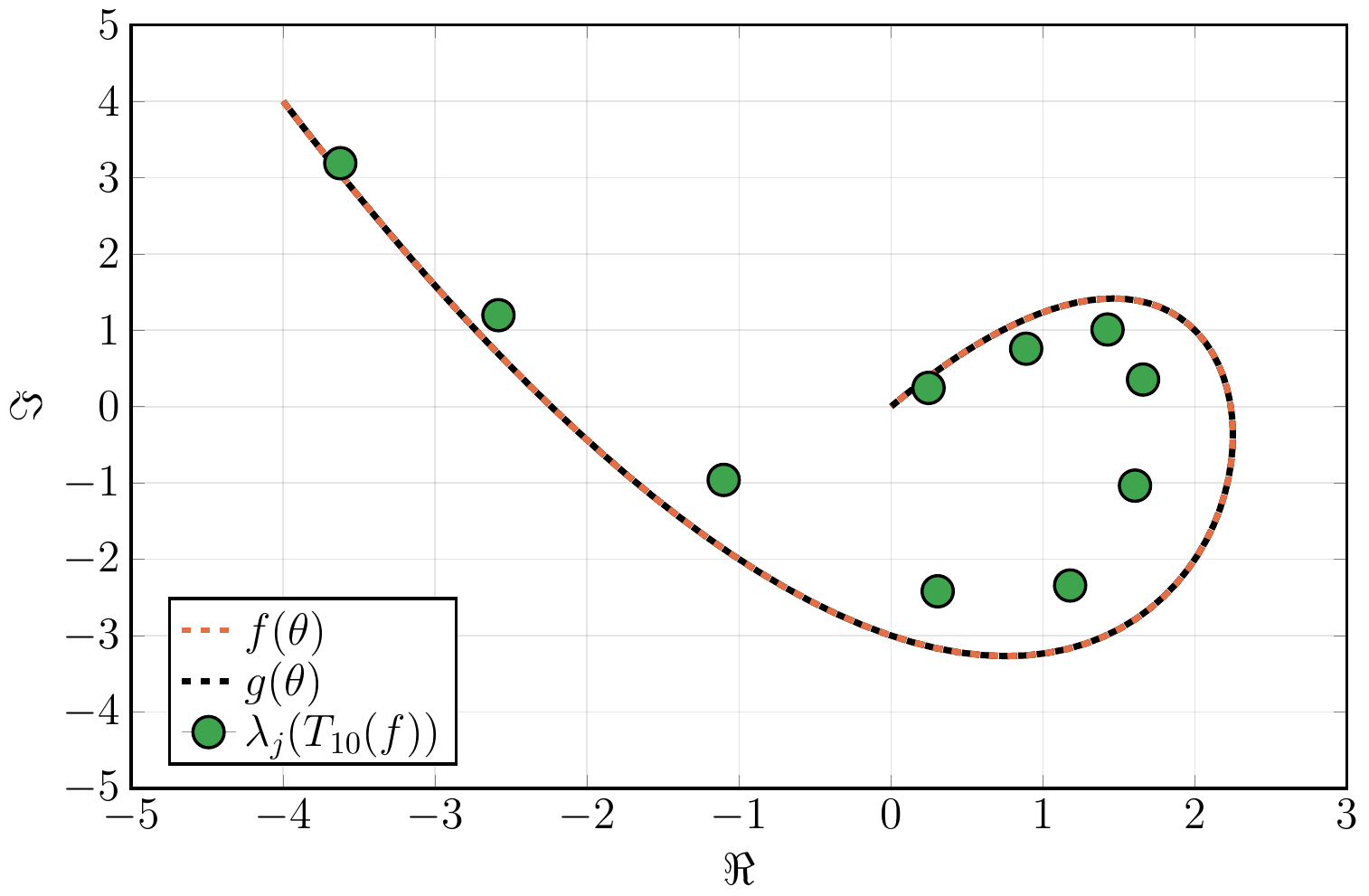}
\includegraphics[width=0.47\textwidth,valign=t]{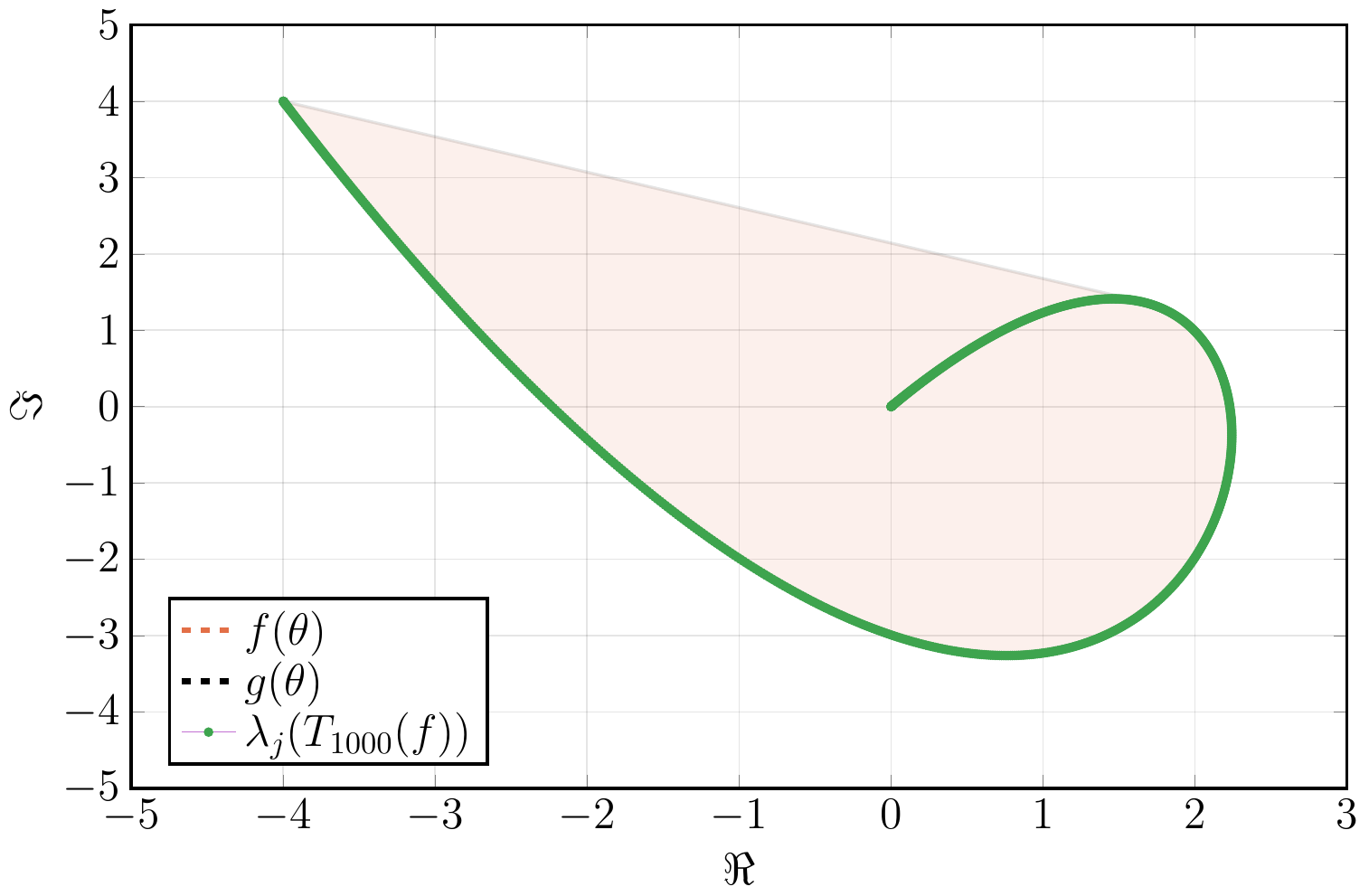}
\caption{[Example~\ref{exmp:3}: Symbol $f(\theta)=2\cos(\theta)-2\cos(2\theta)+\mathbf{i}\left(2\cos(2\theta)-2\cos(3\theta)\right)$] Left: Symbols $f(\theta)$ (red dashed line), and $g(\theta)=f(\theta)$ (black dashed line), and $\lambda_j(T_{10}(f))$ (green circles). 
Right: Symbols $f$ and $g$, and $\lambda_j(T_{1000}(f))$. The convex hull of $f$ is shown in light red.}
\label{fig:exmp:3:symbols}
\end{figure}

\begin{rmrk}
Inspecting the distribution of the eigenvalues in Figure~\ref{fig:exmp:3:spectrum}, for the ``perfect grids'' (middle left and bottom right panels), we note that there are ``distinct grids'' between the points where the derivative of the symbols are zero. See the similar behavior in middle left panel of Figure~\ref{fig:exmp:2:spectrum} of Example~\ref{exmp:2}. 

This observation leads to the conjecture that for real-valued non-monotone symbols an avenue of research is to use matrix-less methods locally between the all the points in $[-\pi,\pi]$ where the derivative of the symbols are zero.
Furthermore, the real part of the eigenvalues of Example~\ref{exmp:4}, shown in the middle left panel of Figure~\ref{fig:exmp4:spectrum}, shows the similar behaviour, but $n=10$ is not large enough such so that any grid points are present in the ``middle part'' centered at $\theta=\pi/2$.  Because the eigenvalues $\lambda_j(T_n(f))$ in Example~\ref{exmp:4} consist of complex-conjugate pairs (for $n$ even), the real parts are of multiplicity two. Assuming the conjecture stated above is correct, one could in theory separate four different grids (or sub symbols with their separate asymptotic expansions), for this problem and then reconstruct the full spectrum.

A possibility to customize the ordering of eigenvalues of Toeplitz matrices generated by real-valued non-monotone symbols, is to add a small monotone imaginary symbol to perturb the symbol into the complex plane, which could yield a possibility to order the eigenvalues consistently by associating them to the original real-valued eigenvalues.
\end{rmrk}

\begin{figure}[!ht]
\centering
\includegraphics[width=0.46\textwidth,valign=t]{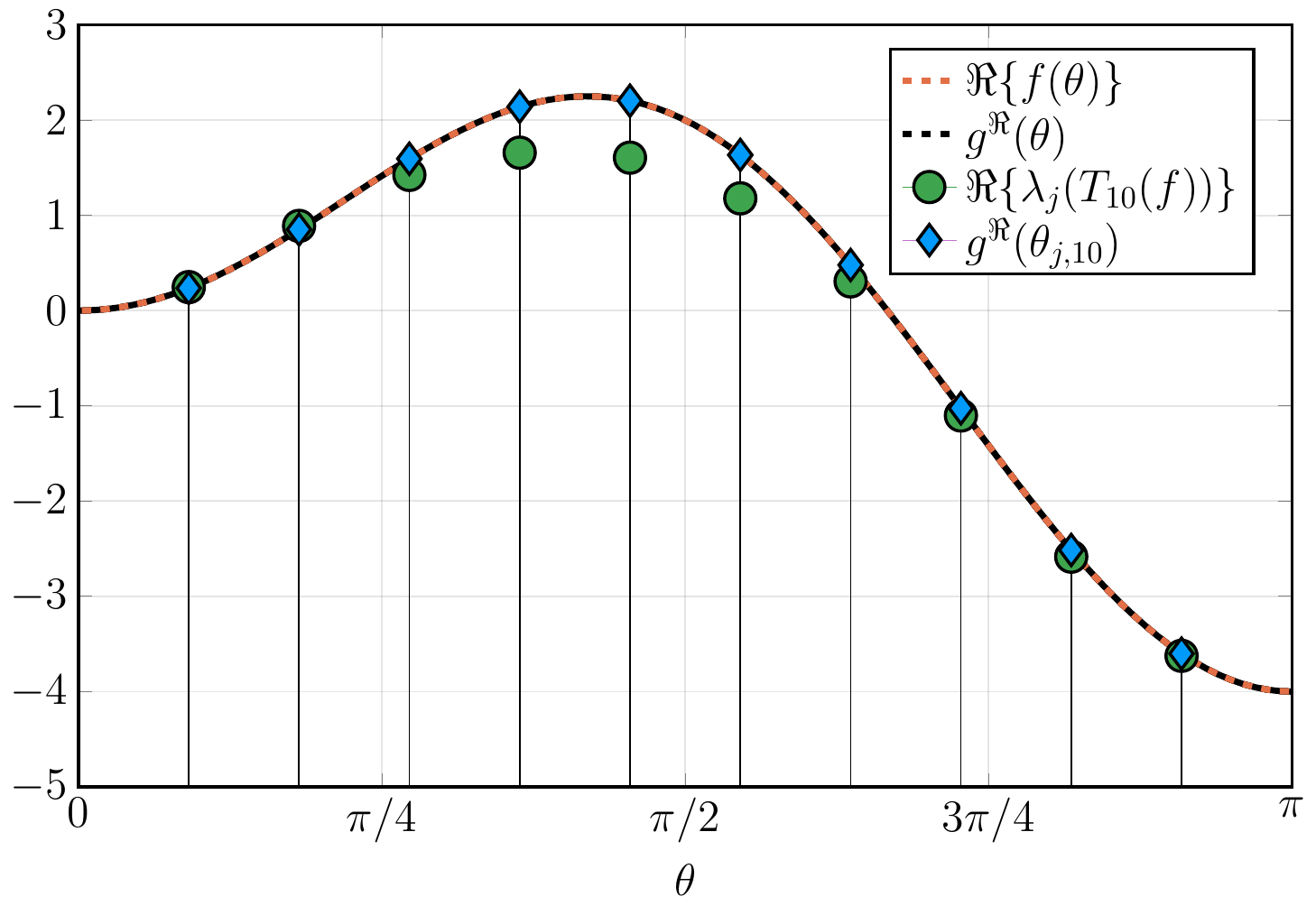}
\includegraphics[width=0.46\textwidth,valign=t]{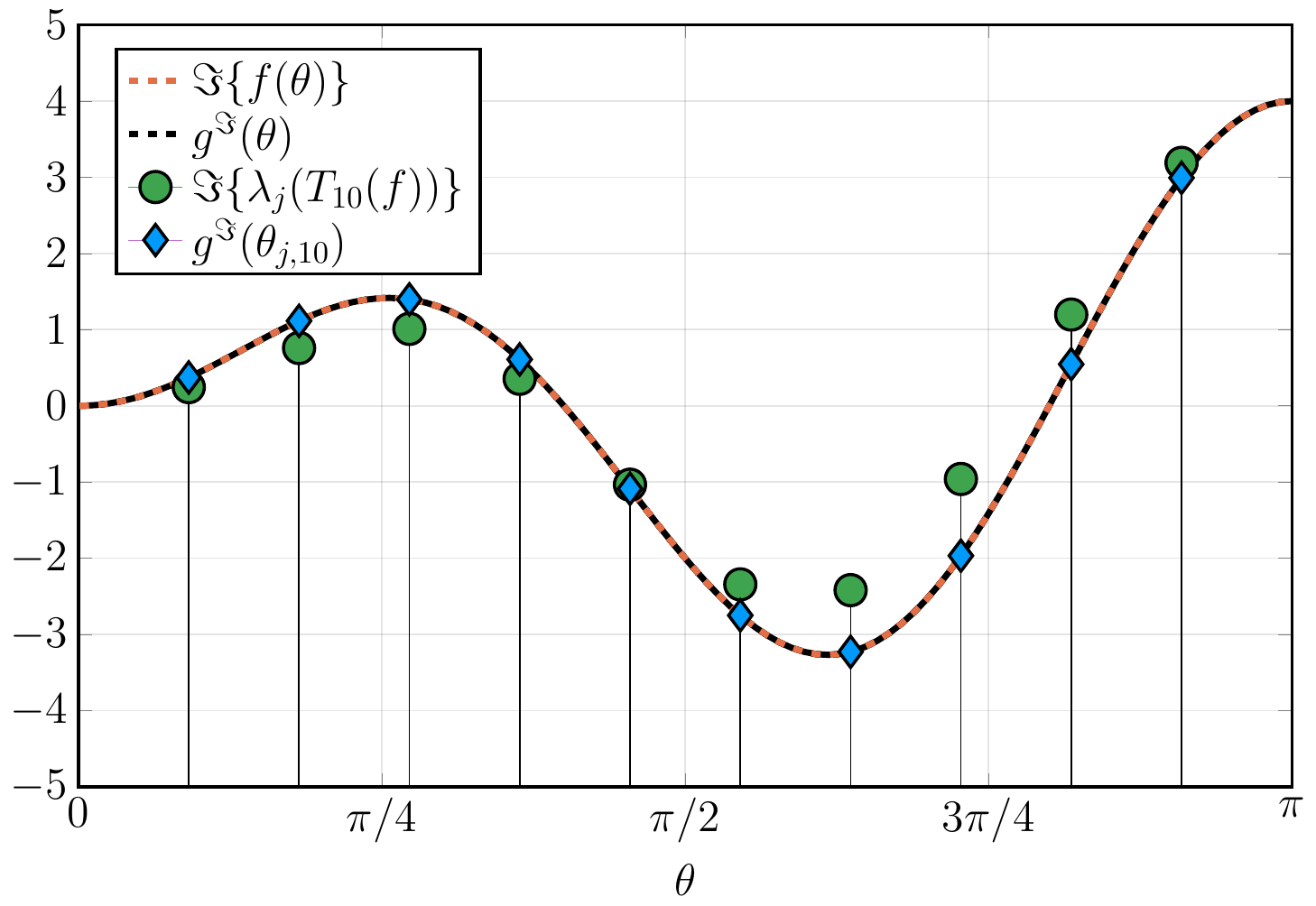}

\includegraphics[width=0.46\textwidth,valign=t]{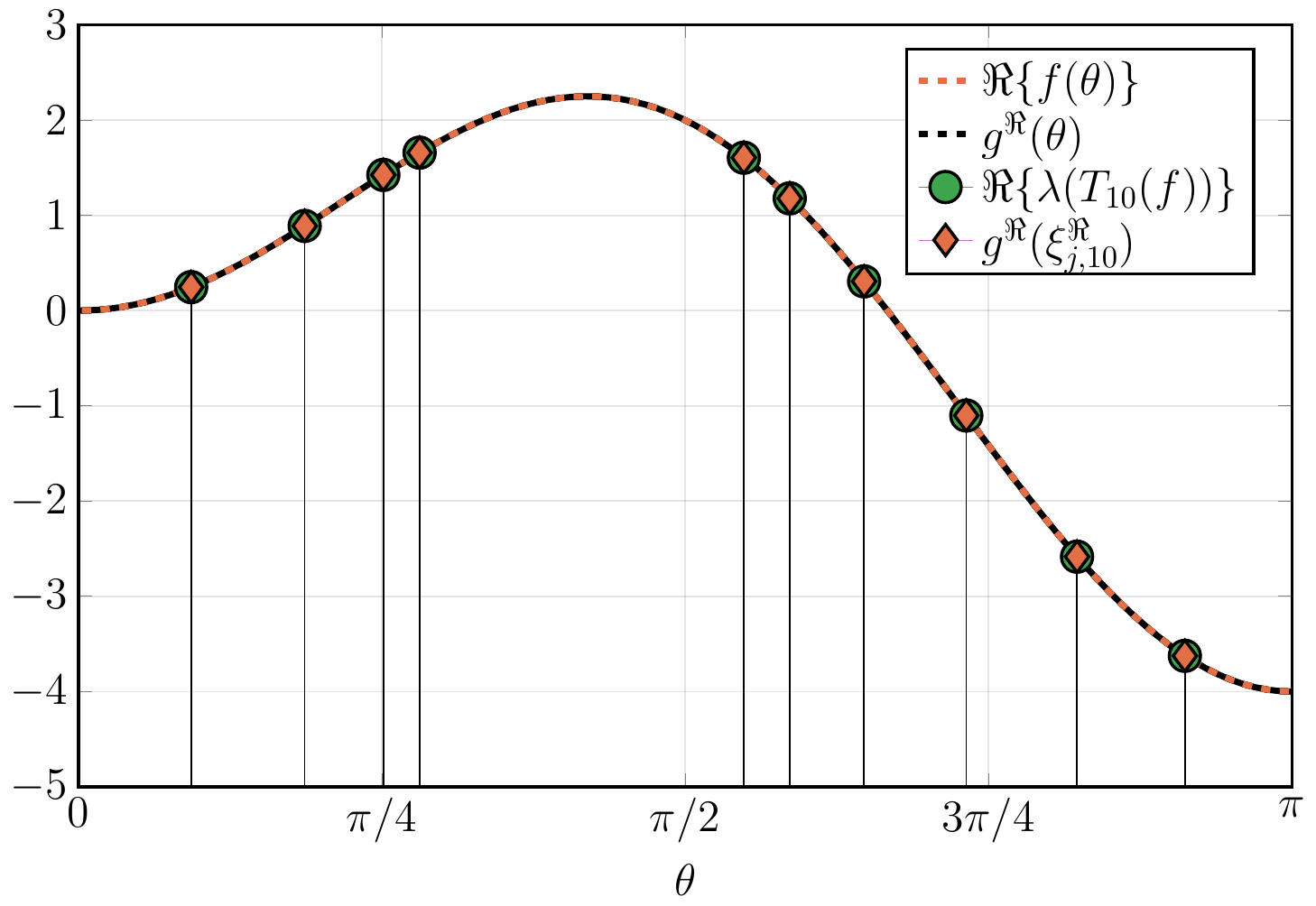}
\includegraphics[width=0.46\textwidth,valign=t]{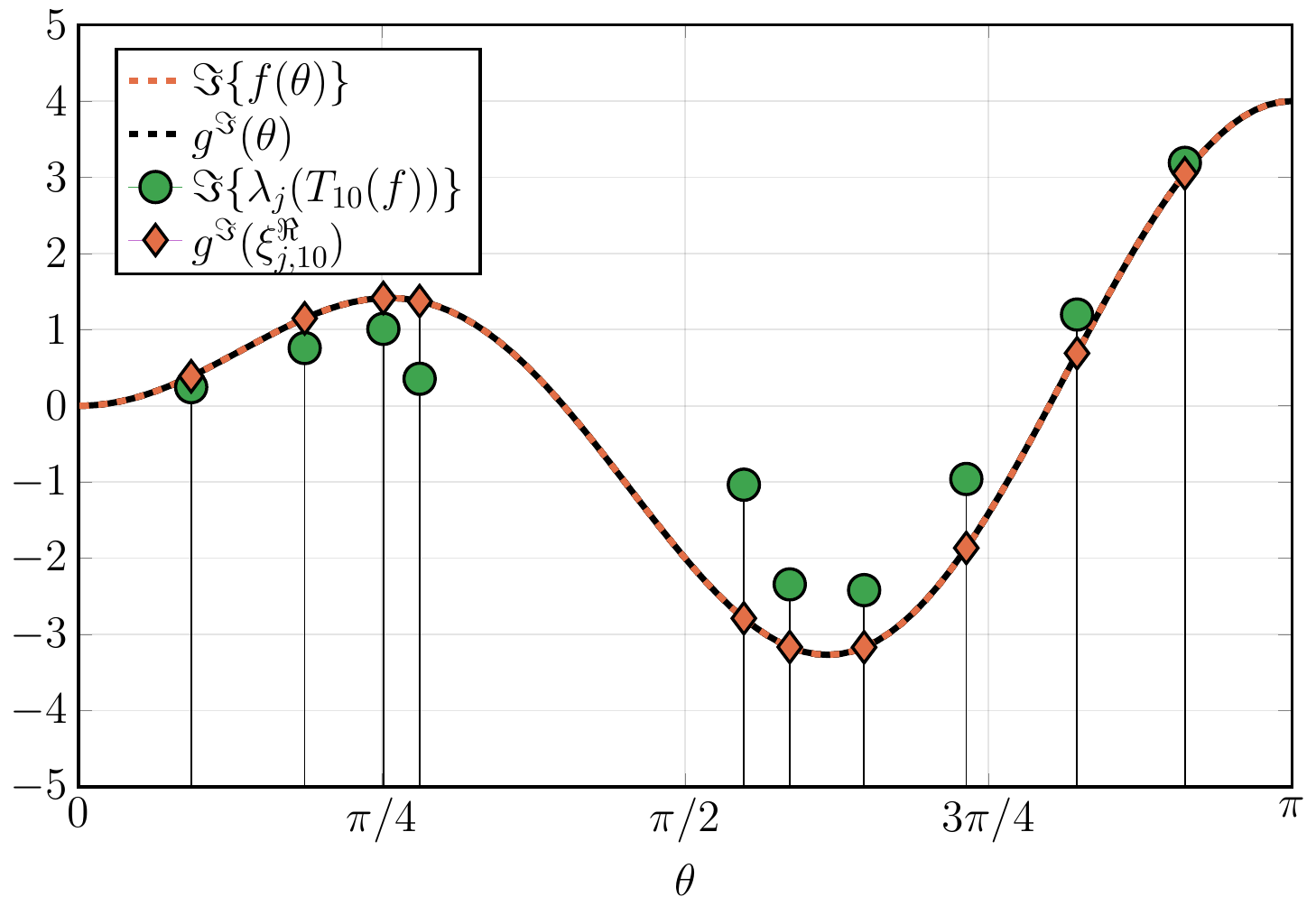}

\includegraphics[width=0.46\textwidth,valign=t]{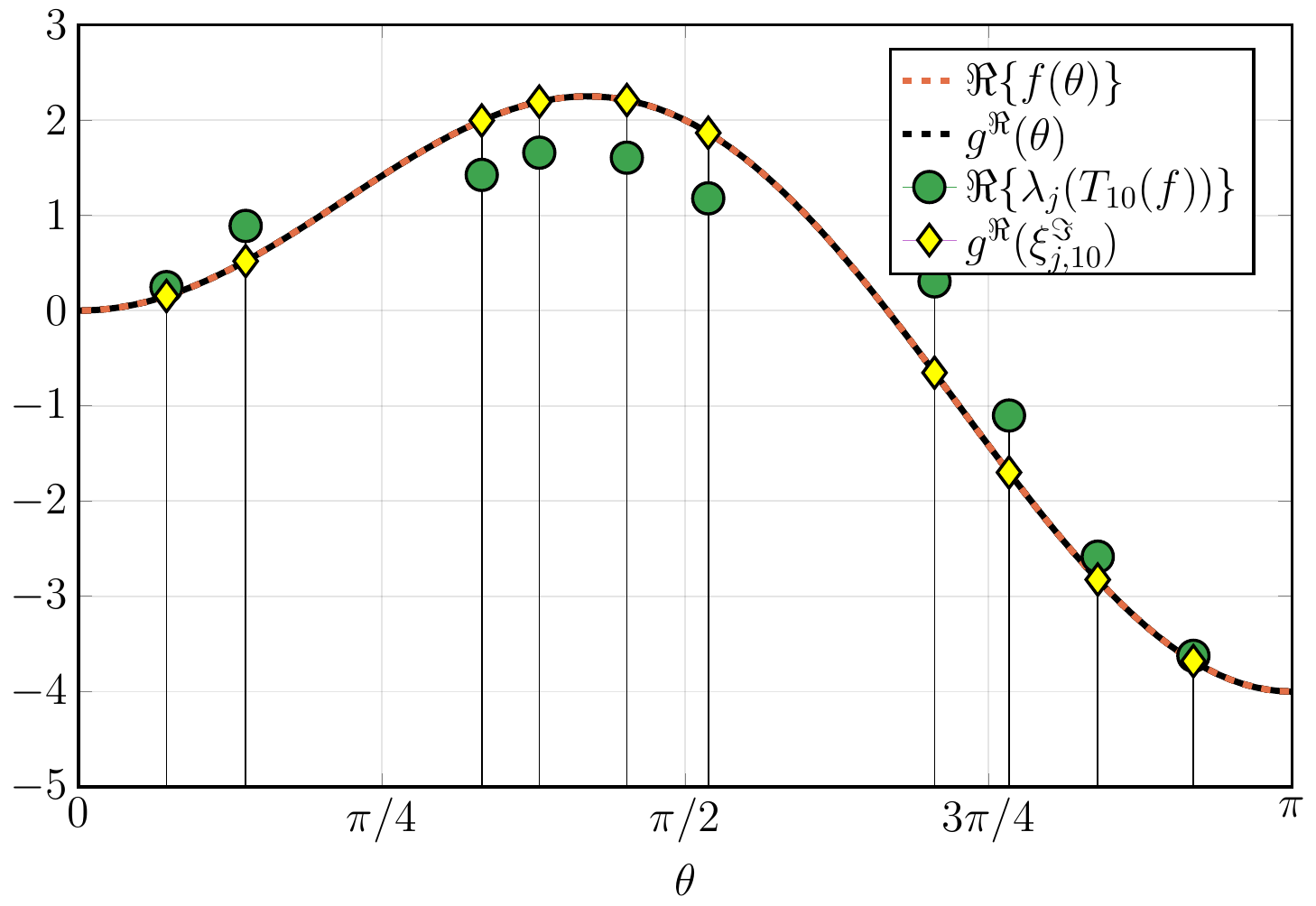}
\includegraphics[width=0.46\textwidth,valign=t]{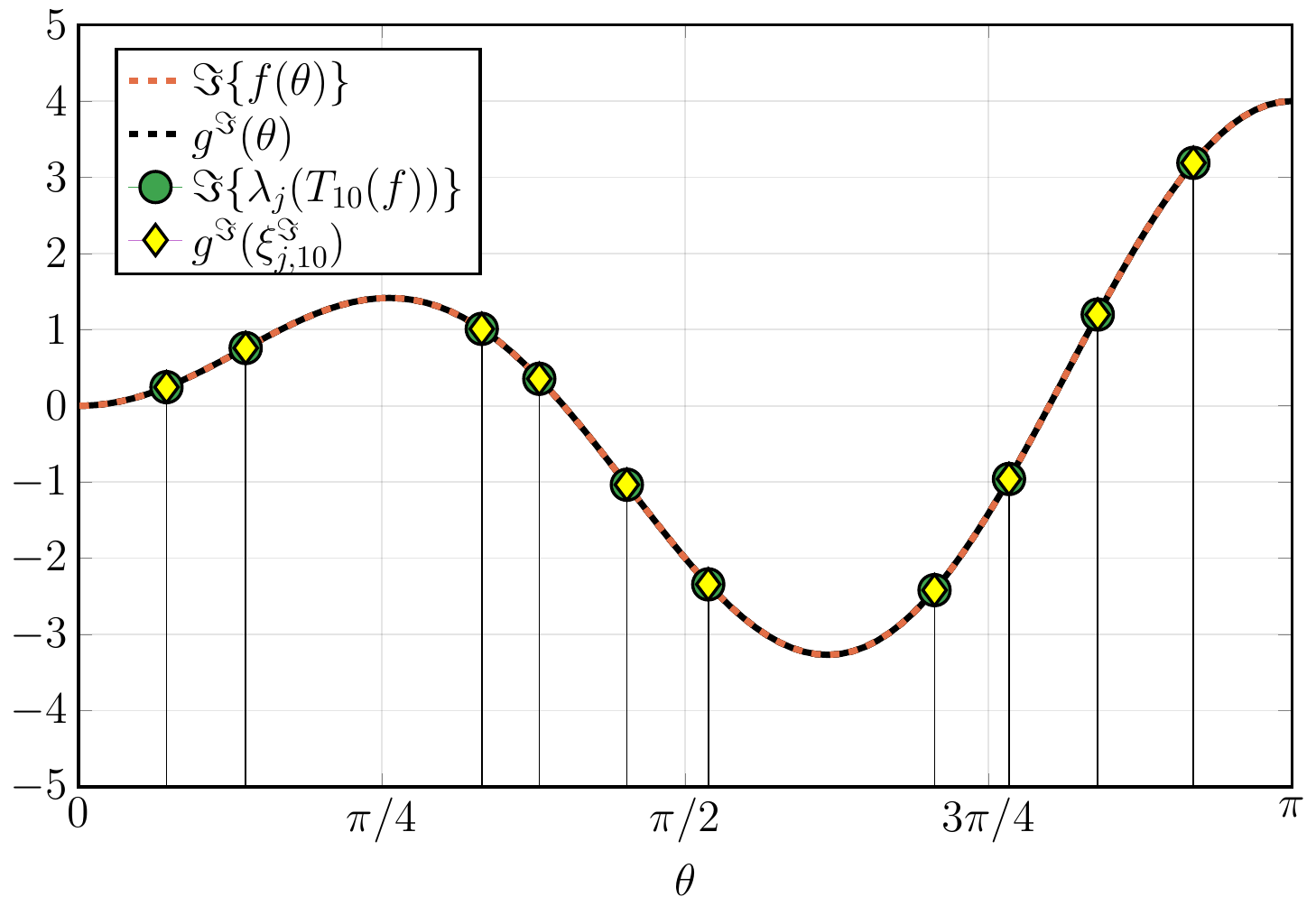}
\caption{[Example~\ref{exmp:3}:  Symbol $f(\theta)=2\cos(\theta)-2\cos(2\theta)+\mathbf{i}\left(2\cos(2\theta)-2\cos(3\theta)\right)$]  Left: Real part of the symbol $f$ (red dashed line) (and $g^\Re=\Re\{f\}$ (black dashed line)). Eigenvalues $\lambda_j(T_{n}(f))$ for $n=10$, and sampling grids $\theta_{j,n}$ (top), $\xi_{j,n}^\Re$ (middle), and $\xi_{j,n}^\Im$ (bottom). Right: Imaginary counterparts of the left panels.}
\label{fig:exmp:3:spectrum}
\end{figure}
\end{exmp}

\clearpage

\begin{exmp}
\label{exmp:4}
The so called Grcar matrix, a real-valued non-normal Toeplitz matrix, was first proposed in~\cite{grcar891} and studied in for example \cite{trefethen91}. The symbol  is
\begin{align}
f(\theta)=-\E^{\mathbf{i}\theta}+1+\E^{-\mathbf{i}\theta}+\E^{-2\mathbf{i}\theta}+\E^{-3\mathbf{i}\theta}\label{eq:exmp4:symbol}
\end{align}
and generates the Toeplitz matrix
\begin{align}
T_n(f)=\left[
\begin{array}{rrrrrrrrrrrrr} 
1&1&1&1\\
-1&1&1&1&1\\
&\ddots&\ddots&\ddots&\ddots&\ddots\\
&&-1&1&1&1&1\\
&&&-1&1&1&1\\
&&&&-1&1&1\\
&&&&&-1&1\\
\end{array}
\right].\nonumber
\end{align}
In the left panel of Figure~\ref{fig:exmp4:symbols} we present the symbol $f$ (red dashed line) and the spectrum for $n=10$ (green circles). The symbol $g$ that describes the eigenvalues of $T_n(f)$ is not known in closed form, and instead we here show a numerical approximation of it, $\tilde{g}$ (black line), computed in Example~\ref{exmp:8} of this article.

In the right panel of Figure~\ref{fig:exmp4:symbols} we present the symbol $f$, the numerical approximation $\tilde{g}$, and the convex hull of $f$ in light red. Again, as in Example~\ref{exmp:1}, the numerical computation of eigenvalues of $T_n(f)$, for $n=1000$, yields numerically unstable results closely related to the pseudospectrum, denoted here by $\Psi_j(T_{n}(f))$ (pink circles) and $\Psi_j(T_{n}^\mathrm{T}(f))$ (blue circles).  By $\lambda_j(T_{n}(f))$ (green circles) we denote the true eigenvalues, here computed using high precision computations with \textsc{GenericLinearAlgebra.jl} in \textsc{Julia}.
\begin{figure}[!ht]
\centering
\includegraphics[width=0.47\textwidth]{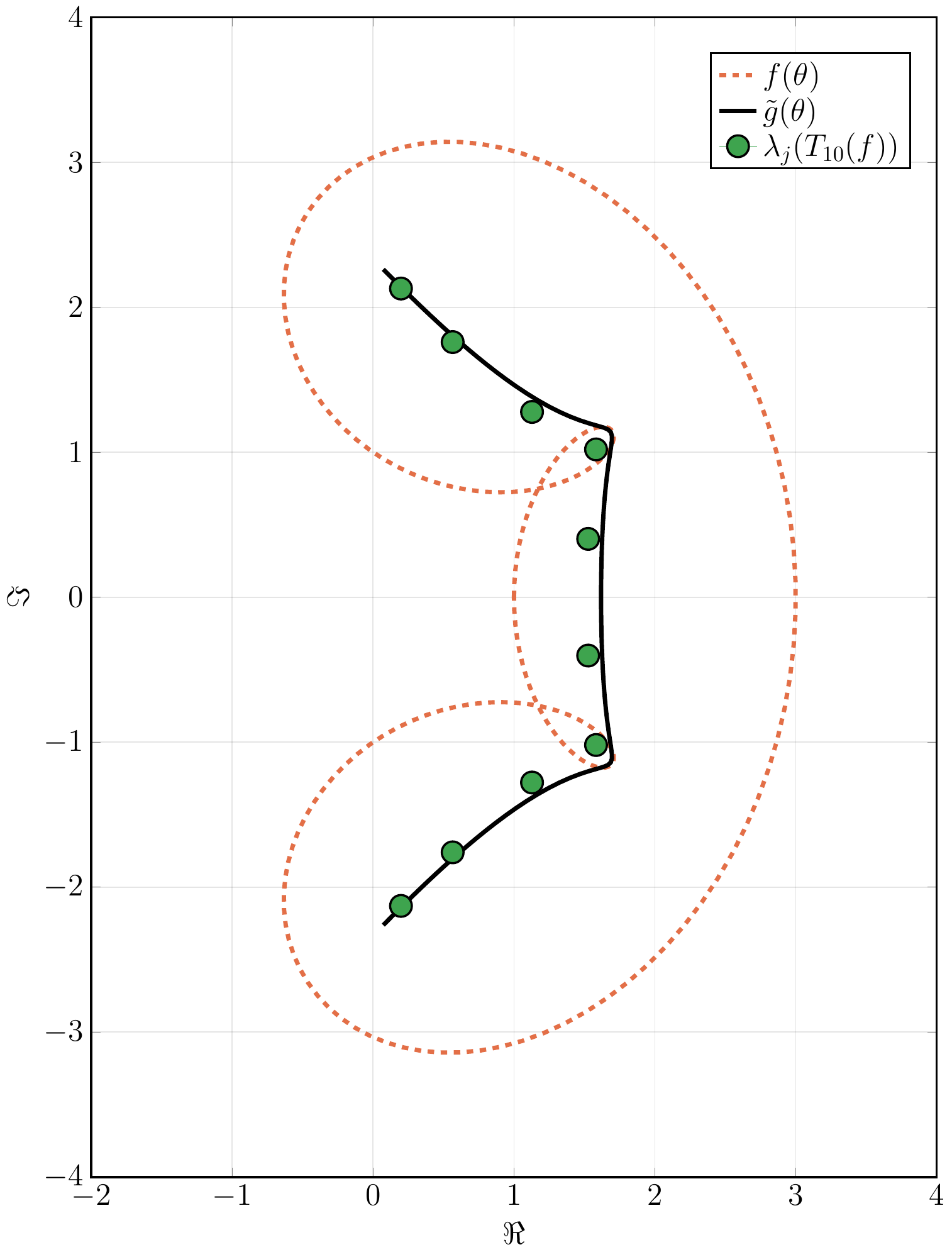}
\includegraphics[width=0.47\textwidth]{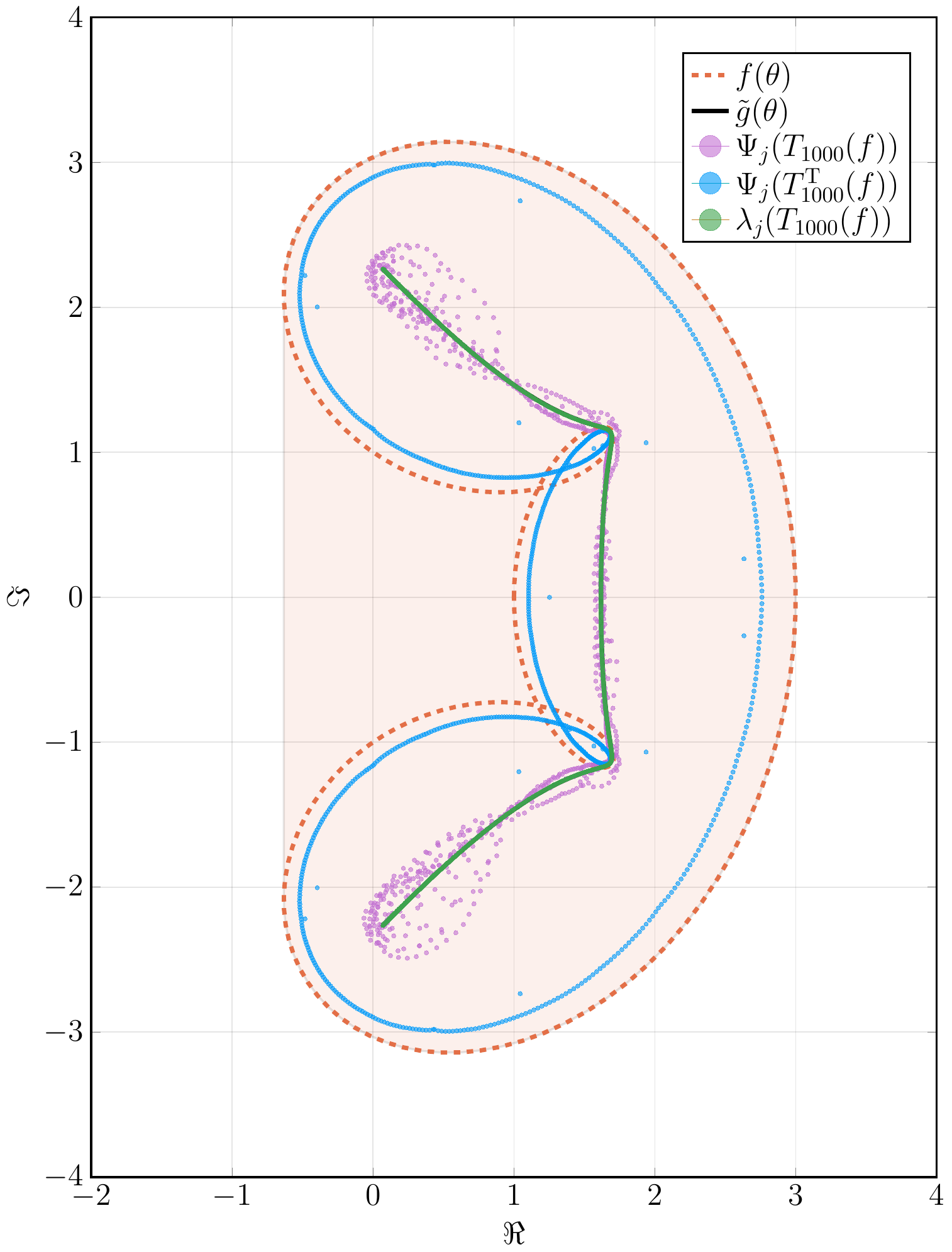}
\caption{[Example \ref{exmp:4}: Symbol $f(\theta)=-\E^{\mathbf{i}\theta}+1+\E^{-\mathbf{i}\theta}+\E^{-2\mathbf{i}\theta}+\E^{-3\mathbf{i}\theta}$] Left: Symbols $f(\theta)$ (red dashed line), and the numerical approximation $\tilde{g}(\theta)$ (black line), and $\lambda_j(T_{10}(f))$ (green circles). 
Right: Symbols $f$ and $\tilde{g}$, numerically computed eigenvalues, for $n=1000$,  $\Psi_j(T_{n}(f))$, $\Psi_j(T_{n}^{\mathrm{T}}(f))$, and $\lambda_j(T_{n}(f))$. The convex hull of $f$ is indicated in light red.}
\label{fig:exmp4:symbols}
\end{figure}

\noindent In Figure~\ref{fig:exmp4:spectrum} we again present the same information as in Figures \ref{fig:exmp:2:spectrum} and \ref{fig:exmp:3:spectrum}.  Since $g$ is not known, we use, as explained in Example~\ref{exmp:8}, Algorithms~\ref{algo:1} and \ref{algo:2}, and the \textsc{Julia} package \textsc{ApproxFun.jl}~\cite{olver141} to construct the approximation $\tilde{g}(\theta)$, which is used to approximate $\xi_{j,10}^\Re$ and $\xi_{j,10}^\Im$ used in Figure~\ref{fig:exmp4:spectrum}.

\begin{figure}[!ht]
\centering
\includegraphics[width=0.46\textwidth]{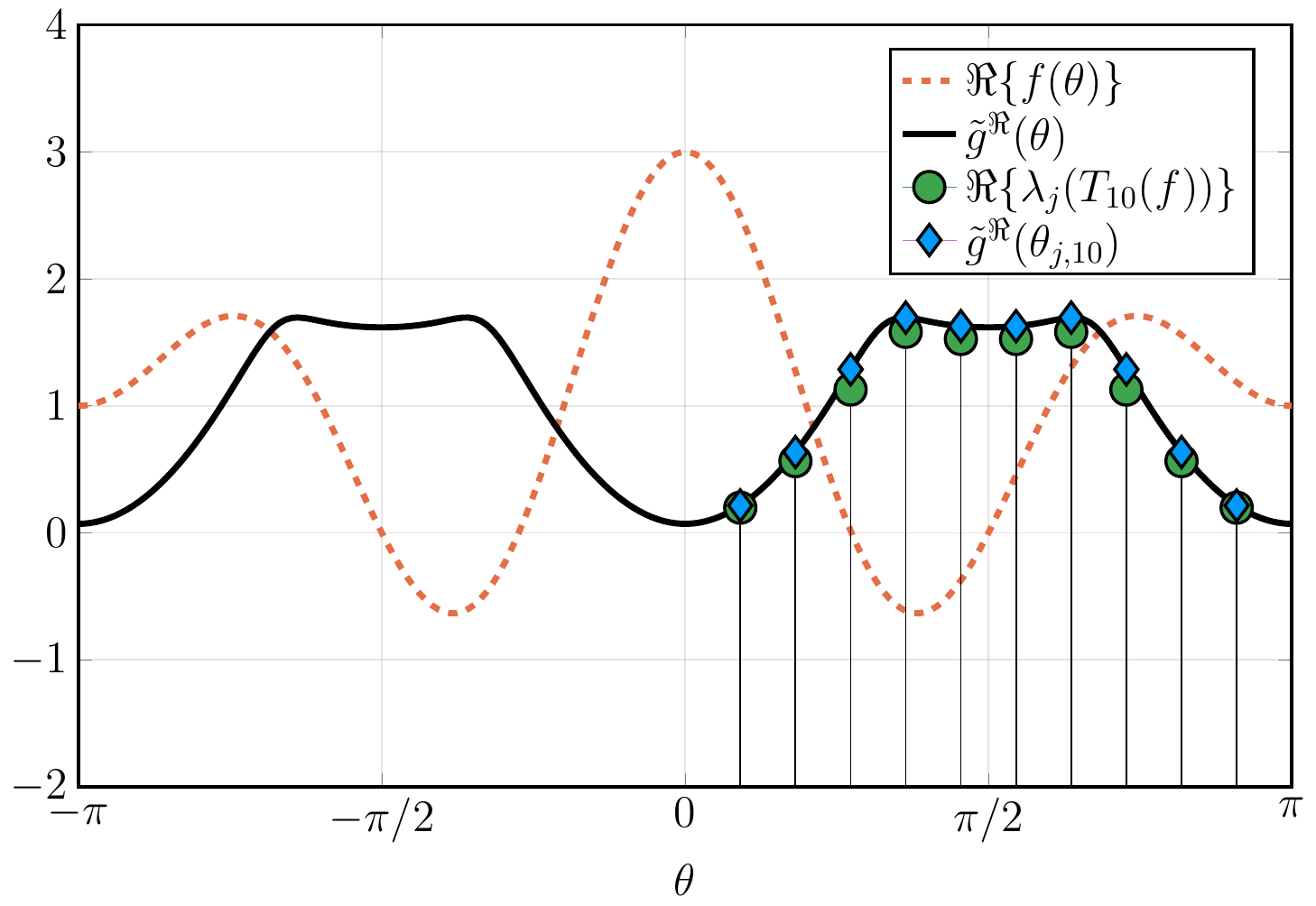}
\includegraphics[width=0.46\textwidth]{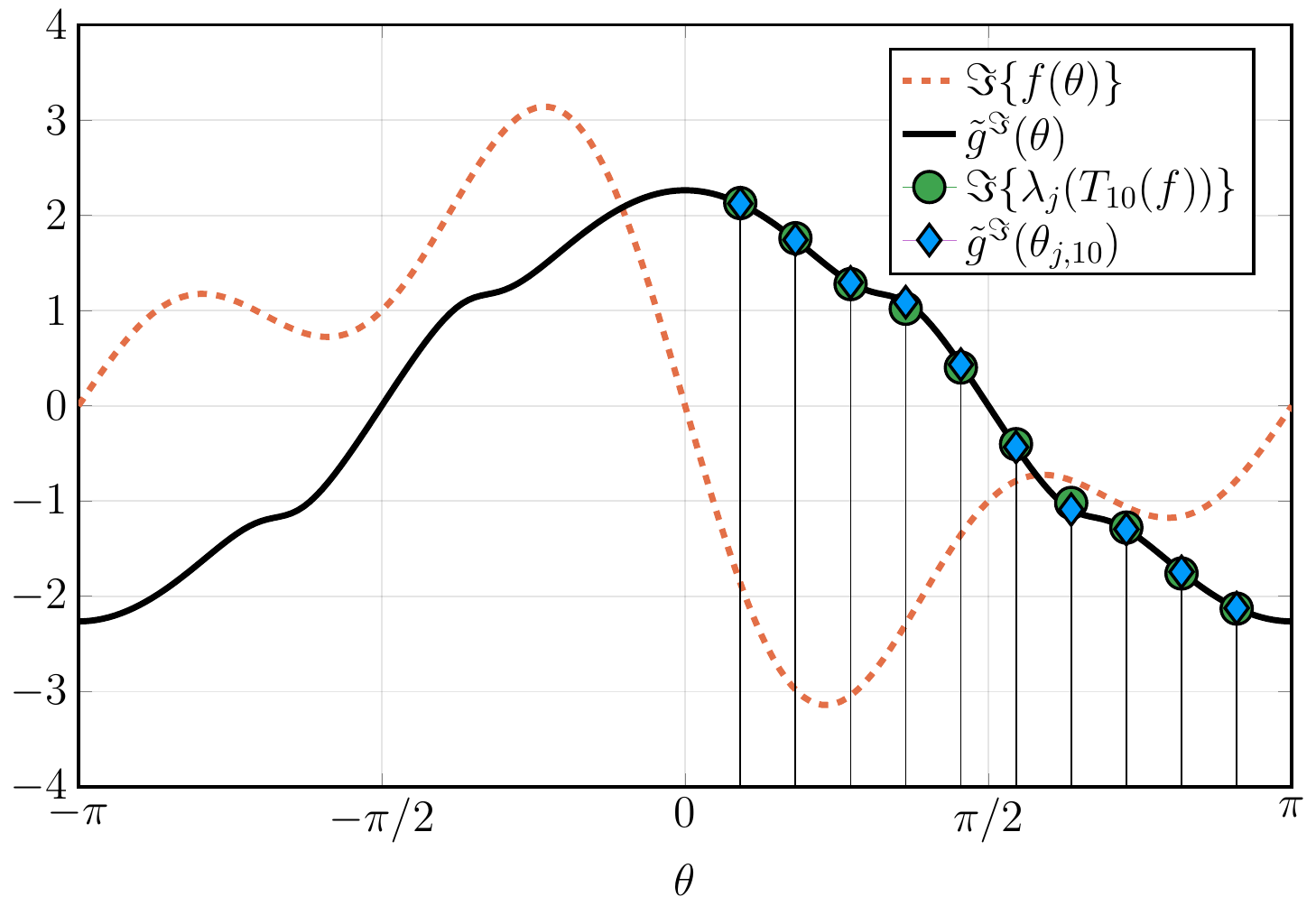}

\includegraphics[width=0.46\textwidth]{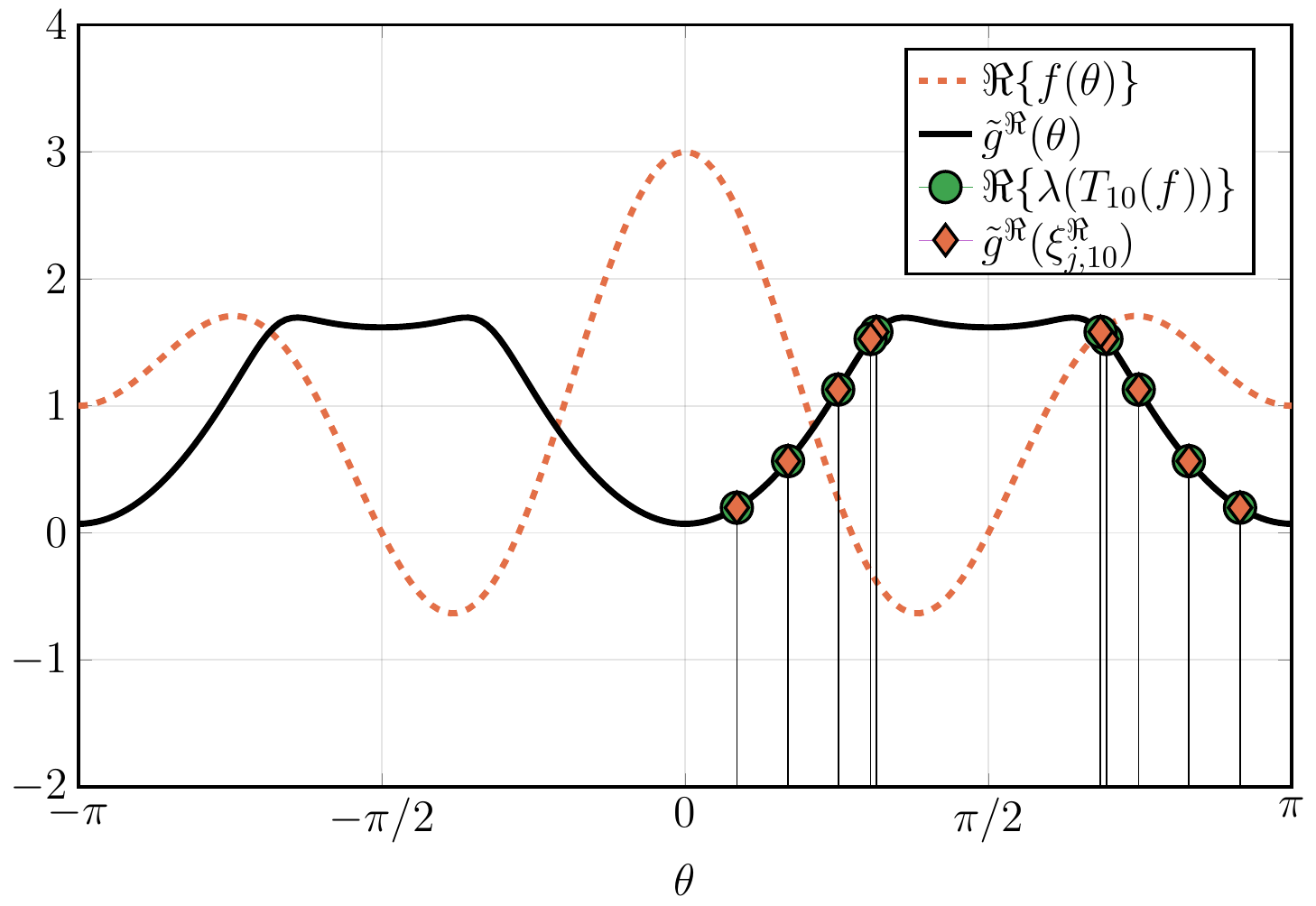}
\includegraphics[width=0.46\textwidth]{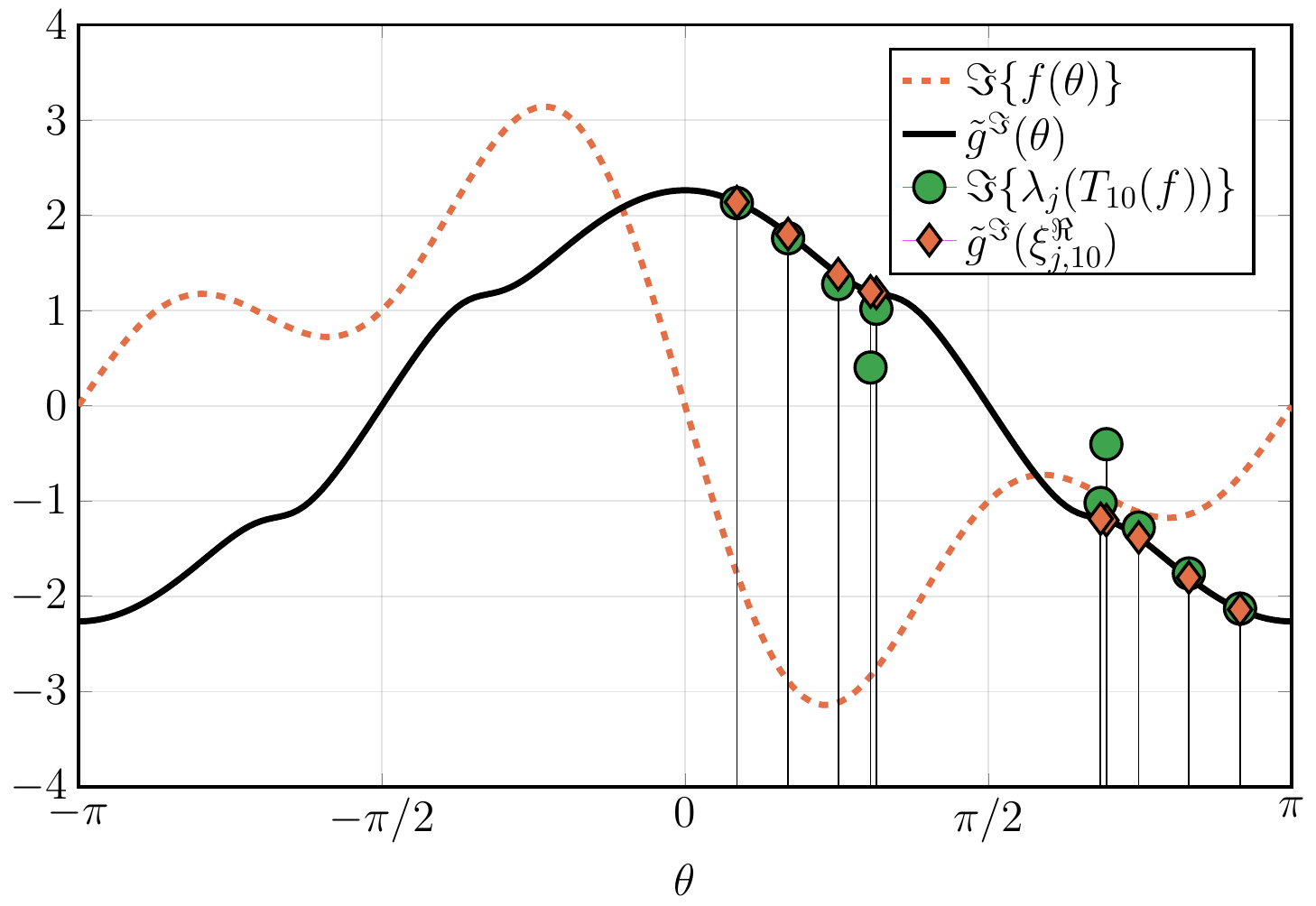}

\includegraphics[width=0.46\textwidth]{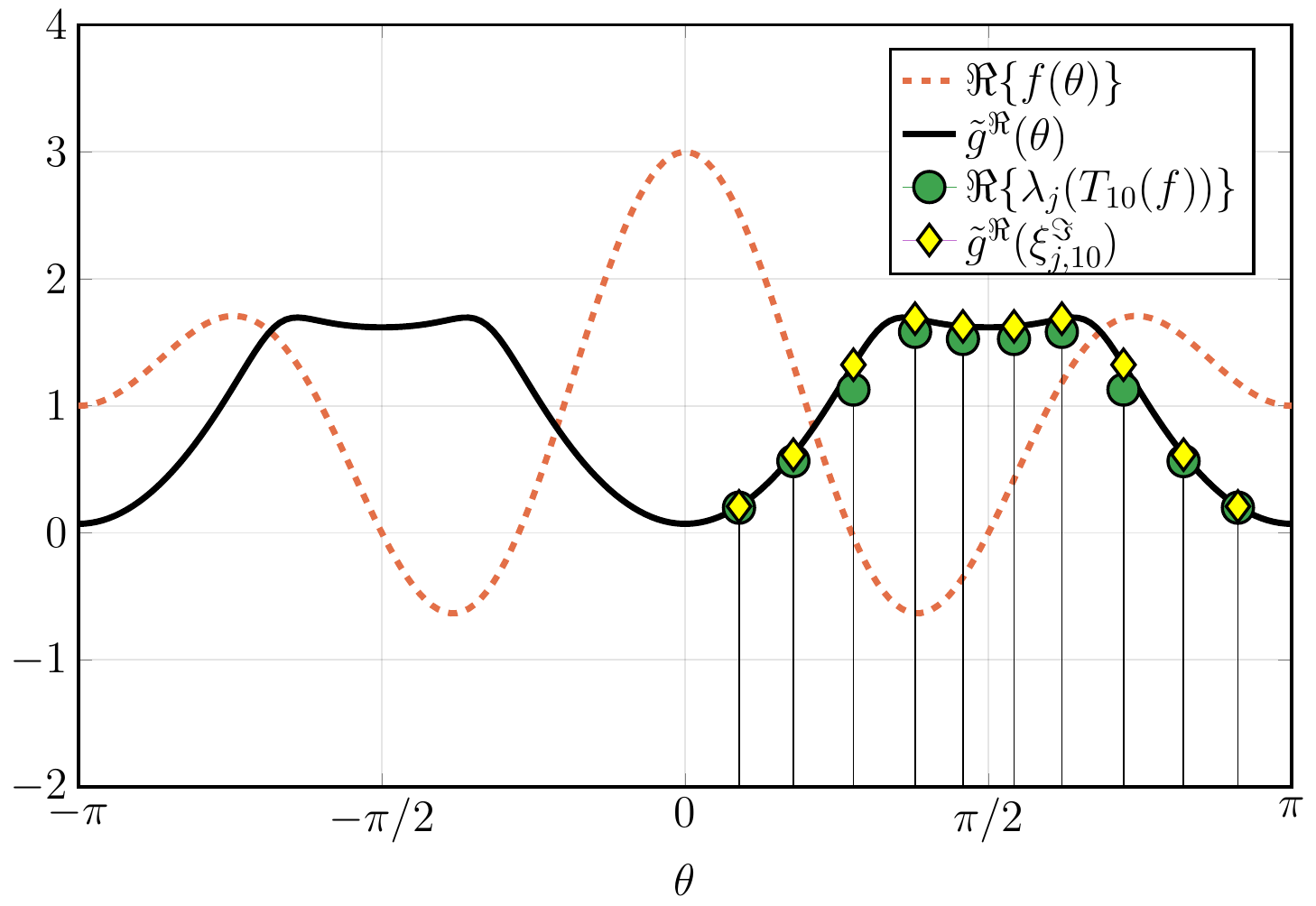}
\includegraphics[width=0.46\textwidth]{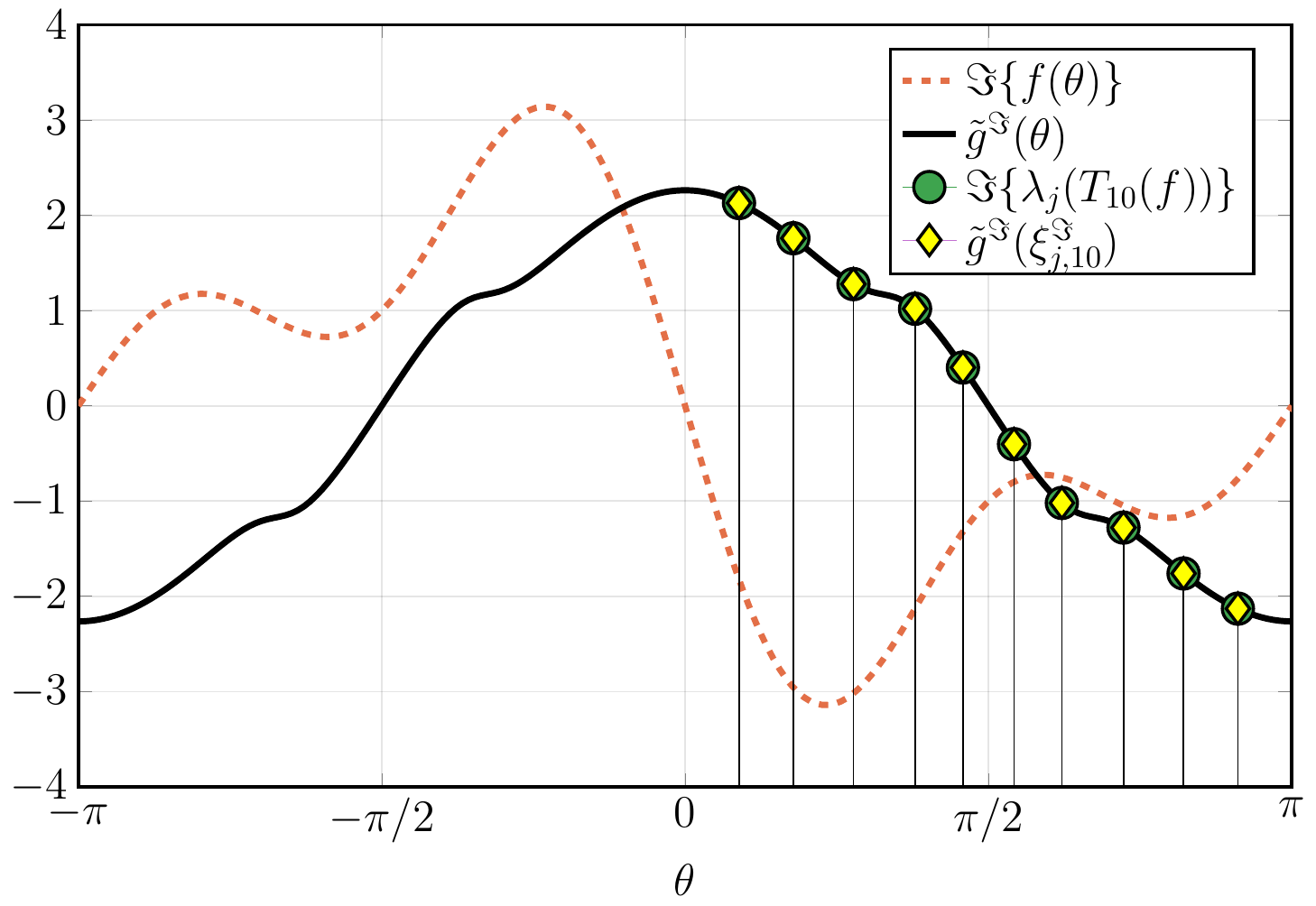}
\caption{[Example \ref{exmp:4}: Symbol $f(\theta)=-\E^{\mathbf{i}\theta}+1+\E^{-\mathbf{i}\theta}+\E^{-2\mathbf{i}\theta}+\E^{-3\mathbf{i}\theta}$] Left: Real part of the symbol $f$ (red dashed line), and approximation $\tilde{g}^\Re$ (black line) from Example~\ref{exmp:8}. Eigenvalues $\lambda_j(T_{n}(f))$ for $n=10$, and sampling grids $\theta_{j,n}$ (top), $\xi_{j,n}^\Re$ (middle), and $\xi_{j,n}^\Im$ (bottom). Right: Imaginary counterparts of the left panels.}
\label{fig:exmp4:spectrum}
\end{figure}
\end{exmp}

\section{Describing the complex-valued eigenvalue distribution}
\label{sec:describing}

Assuming that $g^\Re$ and $g^\Im$ are real cosine trigonometric (RCTP) symbols associated with a symbol $f$ as in the working hypothesis, we introduce in Section~\ref{sec:describing:approximating} a new \textit{matrix-less} method to accurately approximate the expansion functions $c_k^\Re, c_k^\Im$, for $k=0,\ldots,\alpha$, where we recall that $c_0^\Re\eqqcolon g^\Re$ and $c_0^\Im\eqqcolon g^\Im$. We note that Algorithm~\ref{algo:1} is a further modified and extended version of Algorithm 1 in~\cite{ekstrom193}, than just treating the real and imaginary part separately. We here introduce an ``eigenvalue function'', \texttt{eigfun}, as a ``black box'' argument for the algorithm, to accommodate customized ordering for more complicated spectra. This will also allow for future more complicated symbols and matrices of study, for example, preconditioning and block matrices generated by matrix-valued symbols.
Subsequently, in Section~\ref{sec:describing:constructing} we present a procedure in Algorithm~\ref{algo:2} to obtain an approximation of the symbol $g$, by approximating its Fourier series.

\subsection{Approximating the expansion functions $c_k^\Re$ and $c_k^\Im$ in grid points $\theta_{j,n_0}$}
\label{sec:describing:approximating}
We refer the reader to~\cite{ekstrom171,ahmad171,ekstrom181,ekstrom183,ekstrom193,ekstrom184,ekstrom185,batalshchikov192,barrera181,bogoya151,bogoya171,bottcher101} and the references therein, for the details on the \textit{matrix-less} methods, and the asymptotic expansion of eigenvalues using the spectral symbol.
In~\cite{ekstrom193} we extended these methods and no longer required the spectral symbol as an input argument in the algorithms.

Assuming that the complex eigenvalues of the matrices $T_n(f)$ in the sequence $\{T_n(f)\}_n$ admit an asymptotic expansion in terms of a unknown (or known) function $g(\theta)=g^\Re(\theta)+\mathbf{i}g^\Im(\theta)$ instead of $f$ (or $f=g$), as in our working hypothesis, we can use Algorithm~\ref{algo:1} in order to find approximations of both $g$ and the functions $c_k^\Re$ and $c_k^\Im$ the following formula,
\begin{align}
\lambda_{j}(T_{n_0}(f)) &\approx\sum_{k=0}^\alpha \left\{\tilde{c}_k^\Re(\theta_{j,n_0})h_0^k+\mathbf{i}\tilde{c}_k^\Im(\theta_{j,n_0})h_0^k\right\}\nonumber\\
&=\underbrace{\tilde{g}^\Re(\theta_{j,n_0})+\mathbf{i}\tilde{g}^\Im(\theta_{j,n_0})}_{=\tilde{g}(\theta_{j,n_0})}+\sum_{k=1}^\alpha \left\{\tilde{c}_k^\Re(\theta_{j,n_0})h_0^k+\mathbf{i}\tilde{c}_k^\Im(\theta_{j,n_0})h_0^k\right\},
\label{eq:describing:approximating:expansion}
\end{align}  
where the approximations $\tilde{c}_k^\Re(\theta_{j,n_0})$ and $\tilde{c}_k^\Im(\theta_{j,n_0})$
(where $\tilde g^\Re(\theta_{j,n_0})\coloneqq\tilde c_0^\Re(\theta_{j,n_0})$ and $\tilde g^\Im(\theta_{j,n_0})\coloneqq\tilde c_0^\Im(\theta_{j,n_0})$) are obtained from $\alpha+1$ small matrices $T_{n_0}(f),\ldots,T_{n_\alpha}(f)$. 
The approximation of the eigenvalues of $T_n(f)$, for arbitraty $n\gg n_0$, can be derived by using \eqref{eq:describing:approximating:expansion} and the interpolation--extrapolation technique described in~\cite{ekstrom183}, once for real and once for the imaginary part of the eigenvalues.
In Algorithm~\ref{algo:1} is shown an implementation in \textsc{Julia} of the algorithm that computes the approximations $\tilde c_k^\Re(\theta_{j,n_0})$ and $\tilde c_k^\Im(\theta_{j,n_0})$ for $k=0,\ldots,\alpha$, where the input arguments are $n_0$, $\alpha$, \texttt{eigfun}, and the data type for computation \texttt{T}; the algorithm is written for clarity and not performance. All computations in this article are made with \textsc{Julia} 1.2.0~\cite{bezanson171}, using \texttt{Float64} or \texttt{BigFloat} data types, and the \textsc{GenericLinearAlgebra.jl} package~\cite{noack191}.
\begin{algo}
\label{algo:1}
Approximate expansion functions $c_k^\Re(\theta)$ and $c_k^\Im(\theta)$ for $k=0,\ldots,\alpha$ on the grid $\theta_{j,n_0}$.

\normalfont
{\footnotesize
\begin{lstlisting}[backgroundcolor = \color{jlbackground},
                   language = Julia,
                   xleftmargin = 0.1em,
                   framexleftmargin = 0.1em]
using LinearAlgebra, GenericLinearAlgebra 
setprecision(BigFloat,128)

# Example: C = compute_c(100, 3, eigfun_example_1_and_5, Complex{BigFloat})               
function compute_c(n0 :: Integer, alpha :: Integer, eigfun, T :: DataType) 
  j0 = 1:n0
  E  = zeros(T,alpha+1,n0)
  hs = zeros(real(T),alpha+1)
  for kk = 0:alpha
    nk        = (2^kk)*(n0+1)-1
    jk        = (2^kk)*j0
    hs[kk+1]  = convert(T,1)/(nk+1)
    eTnk      = eigfun(nk,T)
    E[kk+1,:] = eTnk[jk]
  end
  V = zeros(T,alpha+1,alpha+1)
  for ii = 1:alpha+1, jj = 1:alpha+1
    V[ii,jj] = hs[ii]^(jj-1)
  end
  return C=V\E
end
\end{lstlisting}
}
\end{algo}

\subsection{Constructing a function $\tilde{g}\approx g$ from approximations $\tilde{c}_0^\Re(\theta_{j,n_0})$ and $\tilde{c}_0^\Im(\theta_{j,n_0})$}
\label{sec:describing:constructing}
We here assume, for the sake of simplicity, that the sought functions $g^\Re$ and $g^\Im$ 
are real av even, so that we have cosine Fourier series of the form 
\begin{align}
g^\Re(\theta)=\hat{g}_0^\Re+2\sum_{k=1}^\infty \hat{g}_k^\Re\cos(k\theta),\quad \hat{g}_k^\Re\in\mathbb{R},\qquad\qquad
g^\Im(\theta)=\hat{g}_0^\Im+2\sum_{k=1}^\infty \hat{g}_k^\Im\cos(k\theta),\quad \hat{g}_k^\Im\in\mathbb{R}.\nonumber
\end{align}
As we shall see in Examples~\ref{exmp:5}--\ref{exmp:7}, if $g^\Re$ and $g^\Im$ are RCTP,
then we are able to recover the exact expression of $g$ up to machine precision; otherwise, as in Example~\ref{exmp:8}, we will get a truncated representation of the Fourier series of $g$.
More specifically, what we do is the following: we consider the approximations $\tilde{c}_0^\Re(\theta_{j,n_0})$ and $\tilde{c}_0^\Im(\theta_{j,n_0})$ provided by Algorithm~1 and we approximate the first $n_0$ Fourier coefficients $\hat g_0^\Re,\ldots,\hat g_{n_0}^\Re$ and $\hat g_0^\Im,\ldots,\hat g_{n_0}^\Im$ with the numbers $\tilde{\hat{g}}_0^\Re,\ldots,\tilde{\hat{g}}_{n_0}^\Re$ and $\tilde{\hat{g}}_0^\Im,\ldots,\tilde{\hat{g}}_{n_0}^\Im$ obtained by solving the following two linear systems
\begin{align}
\tilde{\hat{g}}_0^\Re+2\sum_{k=1}^{n_0} \tilde{\hat{g}}_k^\Re\cos(k\theta_{j,n_0})=\tilde c_0^\Re(\theta_{j,n_0}), \qquad\qquad \tilde{\hat{g}}_0^\Im+2\sum_{k=1}^{n_0} \tilde{\hat{g}}_k^\Im\cos(k\theta_{j,n_0})=\tilde c_0^\Im(\theta_{j,n_0}),\qquad j=1,\ldots,n_0.\nonumber
\end{align}
The approximated Fourier coefficients and Fourier series of $g=g^\Re+\mathbf{i}g^\Im$ can then, for example, be used to approximate samplings of the function $g$; as for computing the perfect grids $\xi_{j,n}^\Re$ and $\xi_{j,n}^\Im$ used for visualization in Figure~\ref{fig:exmp4:spectrum}.
\clearpage

\begin{algo}
\label{algo:2}
Compute approximations $\tilde{\hat{g}}_k$ of the Fourier coefficients $\hat{g}_k$ of $g(\theta)$.
\normalfont
{\footnotesize
\begin{lstlisting}[backgroundcolor = \color{jlbackground},
                   language = Julia,
                   xleftmargin = 0.1em,
                   framexleftmargin = 0.1em]
# Example: ghattildeRe = compute_ghattilde(real.(C[1,:]))
function compute_ghattilde(c0 :: Array{T,1}) where T
  n0     = length(c0)
  t      = LinRange(convert(T,pi)/(n0+1), n0*convert(T,pi)/(n0+1), n0)
  G      = zeros(T,n0,n0)
  G[:,1] = ones(T,n0)
  for jj = 2:n0
    G[:,jj] = 2*cos.((jj-1)*t)
  end
  return ghattilde = G\c0
end
\end{lstlisting}
}
\end{algo}

\noindent For completeness, we also define the following support function to construct ($s\times s$ block) Toeplitz matrices given a size $n$ and the first (block) column and (block) row vector of the matrix. For all examples in the current article $s=1$.
{\normalfont
{\footnotesize
\begin{lstlisting}[backgroundcolor = \color{jlbackground},
                   language = Julia,
                   xleftmargin = 0.1em,
                   framexleftmargin = 0.1em]
# Example: Tn = toeplitz(100, Float64[2, -1], Float64[2, -1])
function toeplitz(n :: Integer, vc :: Array{T,1}, vr :: Array{T,1}) where T
  s  = size(vc[1],2)
  Tn = zeros(eltype(T),s*n,s*n)
  for ii = 1:length(vc)
    Tn = Tn + kron(diagm(-ii+1=>ones(eltype(T),n-ii+1)),vc[ii])
  end
  for jj = 2:length(vr)
    Tn = Tn + kron(diagm( jj-1=>ones(eltype(T),n-jj+1)),vr[jj])
  end
  return Tn
end
\end{lstlisting}
}}

\section{Numerical examples}
\label{sec:numerical}
We now employ the proposed Algorithms~\ref{algo:1} and \ref{algo:2} on the matrix sequqnces $\{T_n(f)\}_n$ generated by the symbols $f$ discussed in Examples~\ref{exmp:1}--\ref{exmp:4} to highlight the applicability of the approach and validity of the working hypothesis, in the respective Examples~\ref{exmp:5}--\ref{exmp:8}.

\begin{itemize}
\item Example~\ref{exmp:5}: Symbol $f\neq g$. Only $\tilde{g}=\tilde{c}_0$ is non-zero, since $g(\theta_{j,n})$ gives exact the eigenvalues of $T_n(f)$. The function $g$ is constructed to machine precision with the approximation of its Fourier series;
\item Example~\ref{exmp:6}: Symbol $f=g=c_0^\Re+\mathbf{i}c_0^\Im$. The approximations $\tilde{c}_k^\Re$ and $\tilde{c}_k^\Im$, for $k=0,\ldots,3$, are approximated accurately. The function $g$ is constructed to machine precision with the approximation of its Fourier series;
\item Example~\ref{exmp:7}: Symbol $f=g=c_0^\Re+\mathbf{i}c_0^\Im$. The approximations $\tilde{c}_k^\Re$ and $\tilde{c}_k^\Im$, for $k=0,\ldots,3$, are approximated accurately. The function $g$ is constructed to machine precision with the approximation of its Fourier series;
\item  Example~\ref{exmp:8}: Symbol $f\neq g$. 
The approximations $\tilde{c}_k^\Re$ and $\tilde{c}_k^\Im$, for $k=0,\ldots,3$, are approximated, but have discontinuities for $k>0$.
A truncated RCTP representation of of $g$ is constructed with the approximation of its Fourier series.
\end{itemize}

\begin{exmp}
\label{exmp:5}
We return to the symbol \eqref{eq:exmp1:symbol} in Example~\ref{exmp:1}, and use the proposed Algorithm~\ref{algo:1}. 
First we define the eigenvalue function, \texttt{eigfun}, for Examples~\ref{exmp:1} and \ref{exmp:5} named \texttt{eigfun\_example\_1\_and\_5}, which is used as the third argument in Algorithm~\ref{algo:1}. 
The two arguments are the size $n$ and the data type of the generated matrix, of which the eigenvalues are computed. We here assume that we do not know the symbol $g$. The required data type is then \texttt{Complex\{BigFloat\}} with floating point precision depending on $(n_0,\alpha)$ to get correct eigenvalue approximations. By uncommenting lines 5 and 6 of the function \texttt{eigfun\_example\_1\_and\_5} we can instead use standard double precision \texttt{Complex\{Float64\}}.
{\normalfont
{\footnotesize
\begin{lstlisting}[backgroundcolor = \color{jlbackground},
                   language = Julia,
                   xleftmargin = 0.1em,
                   framexleftmargin = 0.1em]
function eigfun_example_1_and_5(n :: Integer, T :: DataType)
  # f(theta)
  vc = convert.(T,[2+0im, -1+0im])
  vr = convert.(T,[2+0im, -2+1im])
  # g(theta)
  # vc = convert.(T,[2+0im,sqrt(-1+0im)*sqrt(-2+1im)])
  # vr = vc
  Tn  = toeplitz(n,vc,vr)
  eTn = eigvals(Tn) 
  p   = sortperm(real.(eTn))
  return eTn[p]
end
\end{lstlisting}
}}

\noindent In Figure~\ref{fig:exmp:5:expansion} we present the approximations of $\tilde{c}_k^\Re(\theta_{j,n_0})$ (left panel) and $\tilde{c}_k^\Im(\theta_{j,n_0})$ (right panel), $k=0,\ldots,\alpha$, $(n_0,\alpha)=(100,3)$.
As is seen, the only non-zero $\tilde{c}_k^\Re$ and $\tilde{c}_k^\Im$ are the first functions  $\tilde{c}_0^\Re$ and $\tilde{c}_0^\Im$, which is expected since the exact eigenvalues of $T_n(f)$ are given by $g(\theta_{j,n})=c_0^\Re(\theta_{j,n})+\mathbf{i}c_0^\Im(\theta_{j,n})$.

\begin{figure}[!ht]
\centering
\includegraphics[width=0.47\textwidth,valign=t]{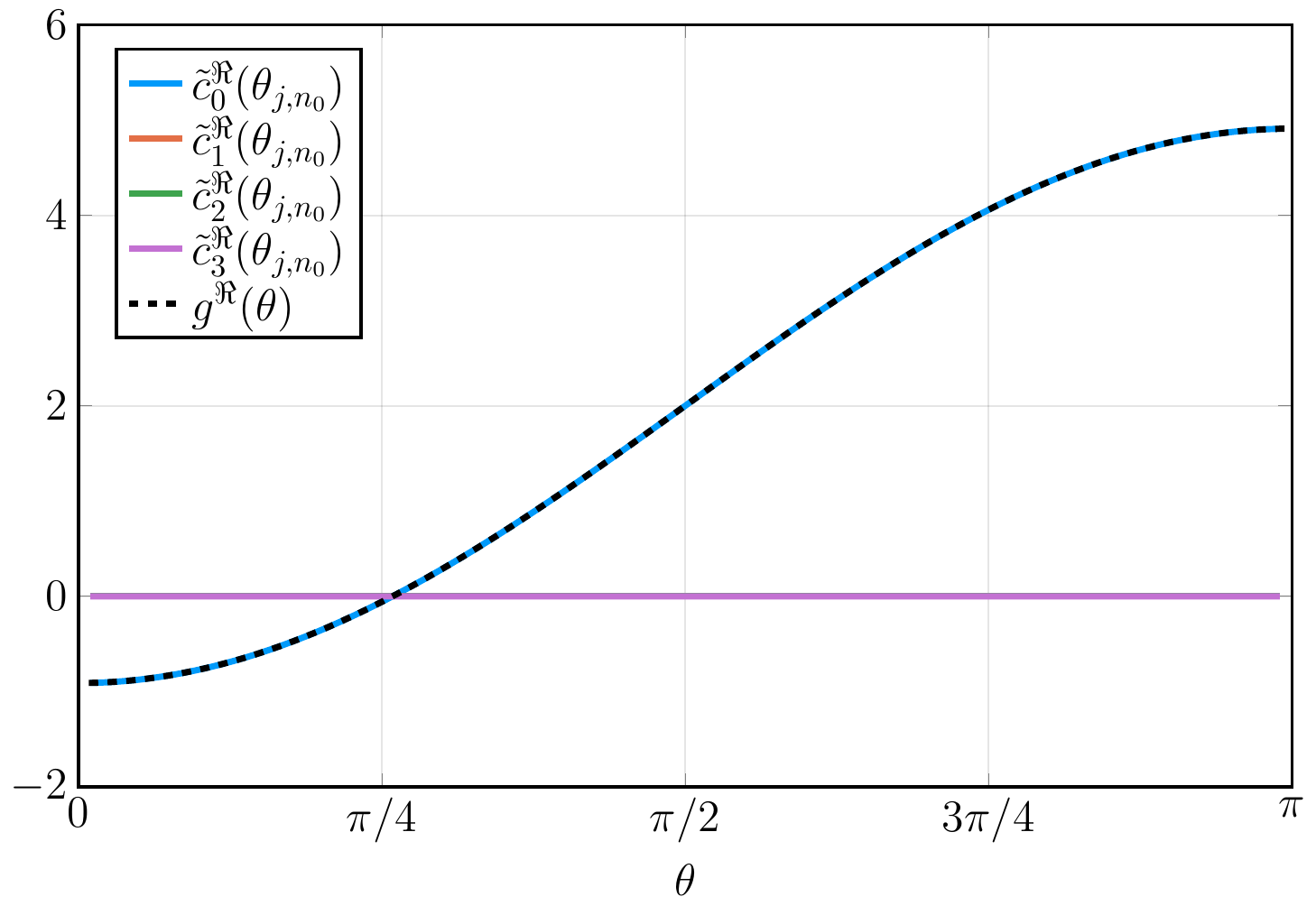}
\includegraphics[width=0.482\textwidth,valign=t]{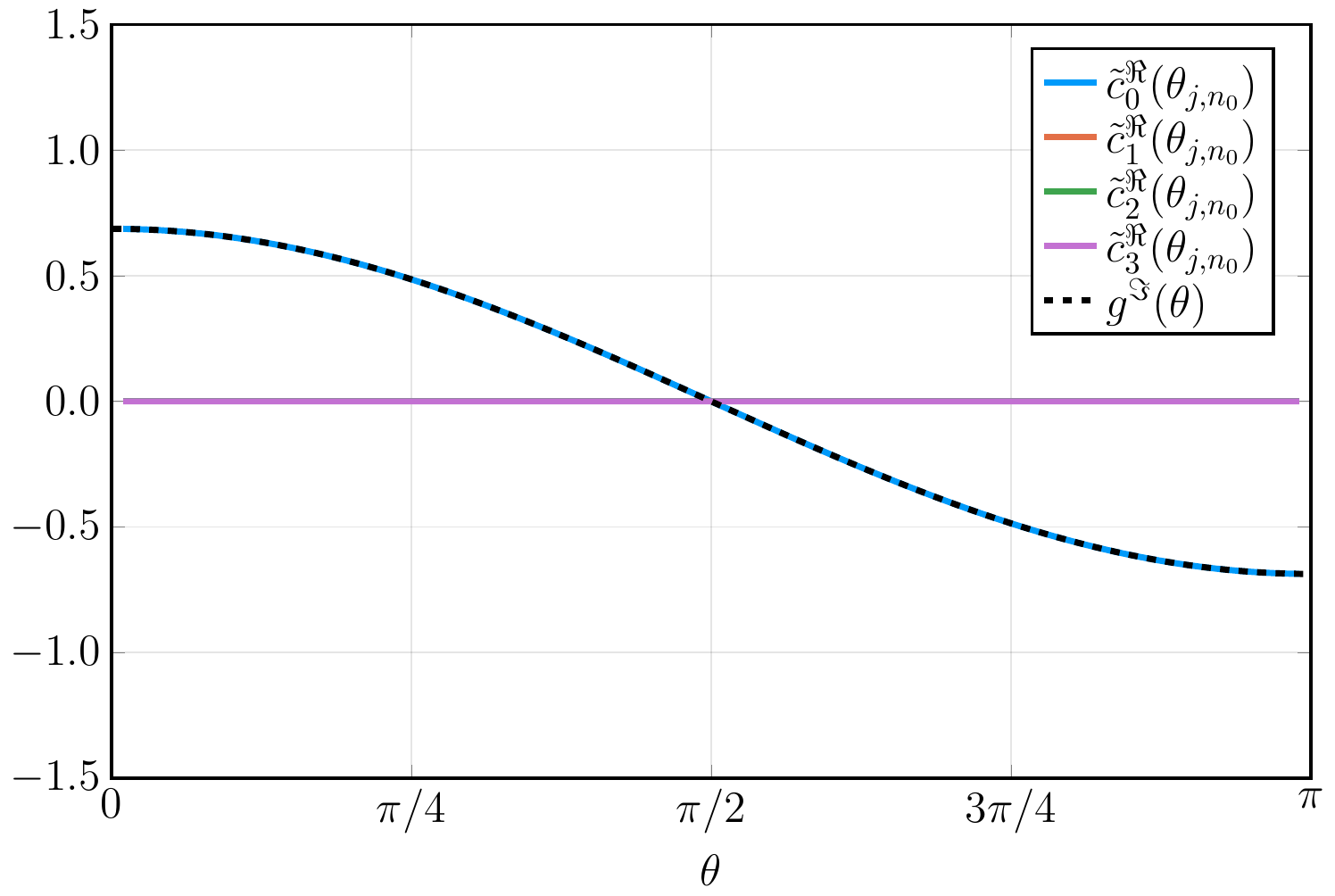}
\caption{[Example~\ref{exmp:5}: Symbol $f(\theta)=-\E^{\mathbf{i}\theta}+2+(-2+\mathbf{i})\E^{-\mathbf{i}\theta}$]
Left:
The computed $\tilde{c}_k^\Re(\theta_{j,n_0})$, $k=0,\ldots\alpha$, with $(n_0,\alpha)=(100,3)$ using Algorithm~\ref{algo:1}.  Only $\tilde{c}_0^\Re(\theta_{j,n_0})$ are non-zero, and matches $g^\Re(\theta)$ defined in \eqref{eq:exmp1:symbolg}.
Right: The corresponding $\tilde{c}_k^\Im(\theta_{j,n_0})$ as in the left panel.
}
\label{fig:exmp:5:expansion}
\end{figure}

If we use Algorithm~\ref{algo:2} to compute the Fourier coefficients of the symbol $g$ we indeed get the numerical approximation of \eqref{eq:exmp1:symbolg} to machine precision.
\end{exmp}

\begin{exmp}
\label{exmp:6}
We here return to the symbol defined in \eqref{eq:exmp2:symbol} in Example~\ref{exmp:2}.
The \texttt{eigfun} function used as an argument in Algorithm~\ref{algo:1} is:
{\normalfont
{\footnotesize
\begin{lstlisting}[backgroundcolor = \color{jlbackground},
                   language = Julia,
                   xleftmargin = 0.1em,
                   framexleftmargin = 0.1em]
function eigfun_example_2_and_6(n :: Integer, T :: DataType)
  vc1 = convert.(T,[0+0im, 1+0im, -1+0im])
  vr1 = vc1
  vc2 = convert.(T,[0+6im, 0-4im,  0+1im])
  vr2 = vc2
  Tn  = toeplitz(n,vc1,vr1)+toeplitz(n,vc2,vr2)
  eTn = eigvals(Tn)
  p   = sortperm(imag.(eTn))
  return eTn[p]
end
\end{lstlisting}
}}

\noindent In Figure~\ref{fig:exmp:6:expansion} we show in the left panel the approximated expansion functions $\tilde{c}_k^\Re(\theta_{j,n_0})$ for $k=0,\ldots,\alpha$, computed using $(n_0,\alpha)=(100,3)$. We see that $\tilde{c}_0^\Re(\theta_{j,n_0})$ (blue line) and $g^\Re(\theta)=\Re\{f(\theta)\}$ (black dashed line) overlap.
In the right panel of Figure~\ref{fig:exmp:6:expansion} present the corresponding functions $\tilde{c}_k^\Im(\theta_{j,n_0})$. Note that $\tilde{c}_k^\Im$ for $k>0$ do not match the expansion functions $\tilde{c}_k$ if only computing the expansion for the matrix sequence $\{T_n(6-8\cos(\theta)+2\cos(2\theta))\}_n$; e.g., shown in \cite[Figure 2.1.3]{ekstrom184} and \cite[Figure 9]{ekstrom193}.

Using Algorithm~\ref{algo:2} we recover, to machine precision, the non-zero Fourier coefficients  $\hat{g}_{\pm 1}^\Re=1$, $\hat{g}_{\pm 2}^\Re=-1$, $\hat{g}_0^\Im=6$, $\hat{g}_{\pm 1}^\Im=-4$, and $\hat{g}_{\pm 2}^\Im=1$.

\begin{figure}[!ht]
\centering
\includegraphics[width=0.4824\textwidth,valign=t]{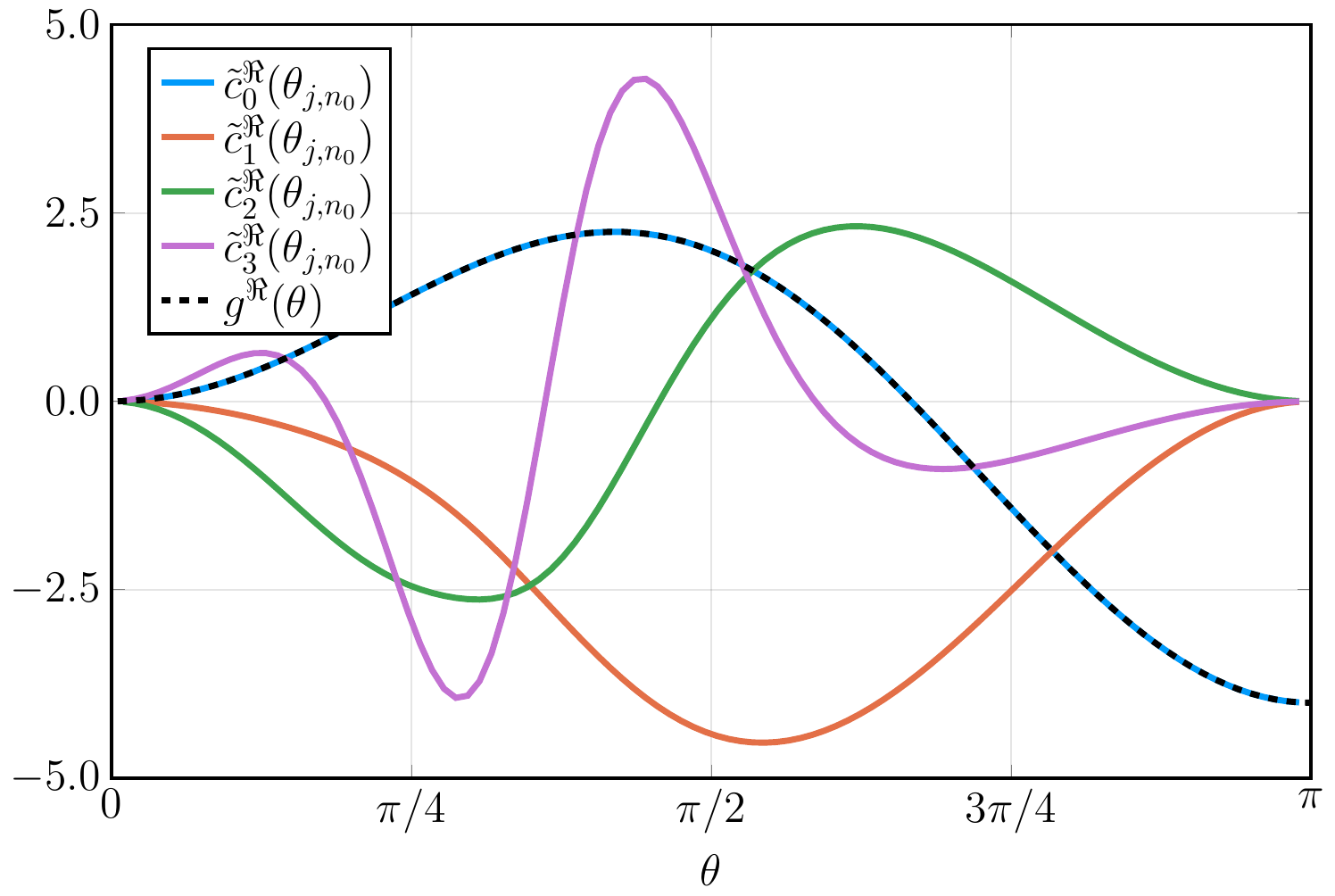}
\includegraphics[width=0.47\textwidth,valign=t]{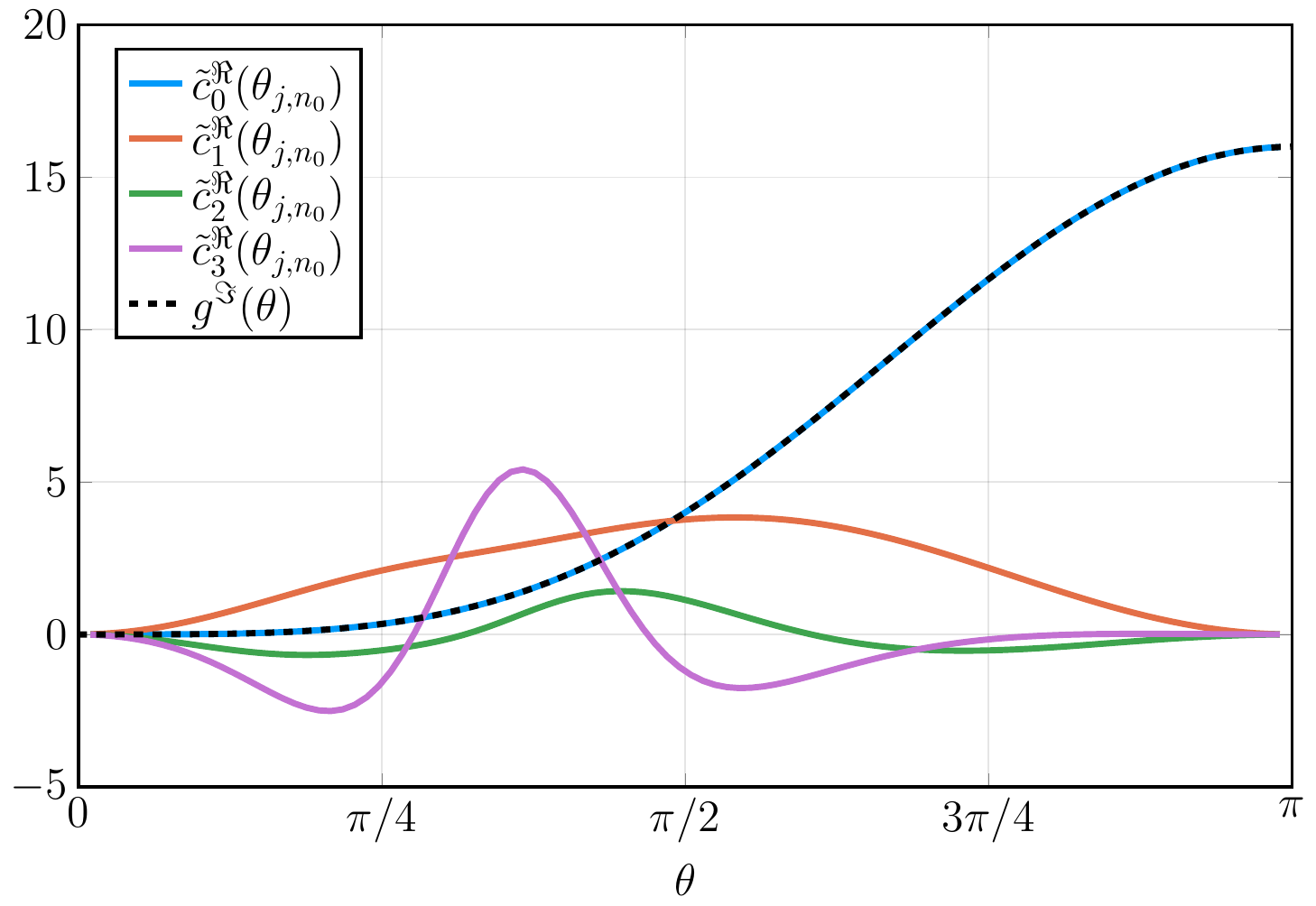}

\caption{[Example~\ref{exmp:6}: Symbol $f(\theta)=2\cos(\theta)-2\cos(2\theta)+\mathbf{i}\left(6-8\cos(\theta)+2\cos(2\theta)\right)$]
Left:
The computed $\tilde{c}_k^\Re(\theta_{j,n_0})$, $k=0,\ldots\alpha$, with $(n_0,\alpha)=(100,3)$ using Algorithm~\ref{algo:1}.
Right: The corresponding $\tilde{c}_k^\Im(\theta_{j,n_0})$ as in the left panel.
}
\label{fig:exmp:6:expansion}
\end{figure}

\end{exmp}

\begin{exmp}
\label{exmp:7}
In this example we continue the investigation of the symbol \eqref{eq:exmp3:symbol} from Example~\ref{exmp:3}.
For this example, the \texttt{eigfun} function used as an argument in Algorithm~\ref{algo:1} is:
{\normalfont
{\footnotesize
\begin{lstlisting}[backgroundcolor = \color{jlbackground},
                   language = Julia,
                   xleftmargin = 0.1em,
                   framexleftmargin = 0.1em]
function eigfun_example_3_and_7(n :: Integer, T :: DataType)
  vc   = convert.(T,[0+0im, 1+0im, -1+1im, 0-1im])
  vr   = vc
  Tn   = toeplitz(n,vc,vr)
  eTn  = eigvals(Tn)
  seTn = zeros(T,n)
  x0   = convert(T, 0+0im)
  for jj = 1:n
    idx      = argmin(abs.(eTn .- x0))
    seTn[jj] = copy(eTn[idx])
    x0       = seTn[jj]
    deleteat!(eTn,idx)
  end
  return seTn
end
\end{lstlisting}
}}

\noindent In Figure~\ref{fig:exmp:7:expansion} we show the approximated expansion functions $\tilde{c}_k^\Re(\theta_{j,n_0})$ (left panels) and $\tilde{c}_k^\Im(\theta_{j,n_0})$ (right panels) for $(n_0,\alpha)=(100,3)$. The two bottom panels show a close-up of the expansion function, and again as expected, Algorithm~\ref{algo:1} approximates the known $g^\Re(\theta_{j,n_0})$ (bottom left panel) and $g^\Im(\theta_{j,n_0})$ (bottom right panel) well with $\tilde{c}_0^\Re(\theta_{j,n_0})$ and $\tilde{c}_0^\Im(\theta_{j,n_0})$.
Again, by using Algorithm~\ref{algo:2} we recover, to machine precision, the non-zero Fourier coefficients of the symbol $g$, namely, $\hat{g}_{\pm 1}^\Re=1$, $\hat{g}_{\pm 2}^\Re=-1$, $\hat{g}_{\pm 2}^\Im=1$, and $\hat{g}_{\pm 3}^\Im=-1$.
\begin{figure}[!ht]
\centering
\includegraphics[width=0.47\textwidth,valign=t]{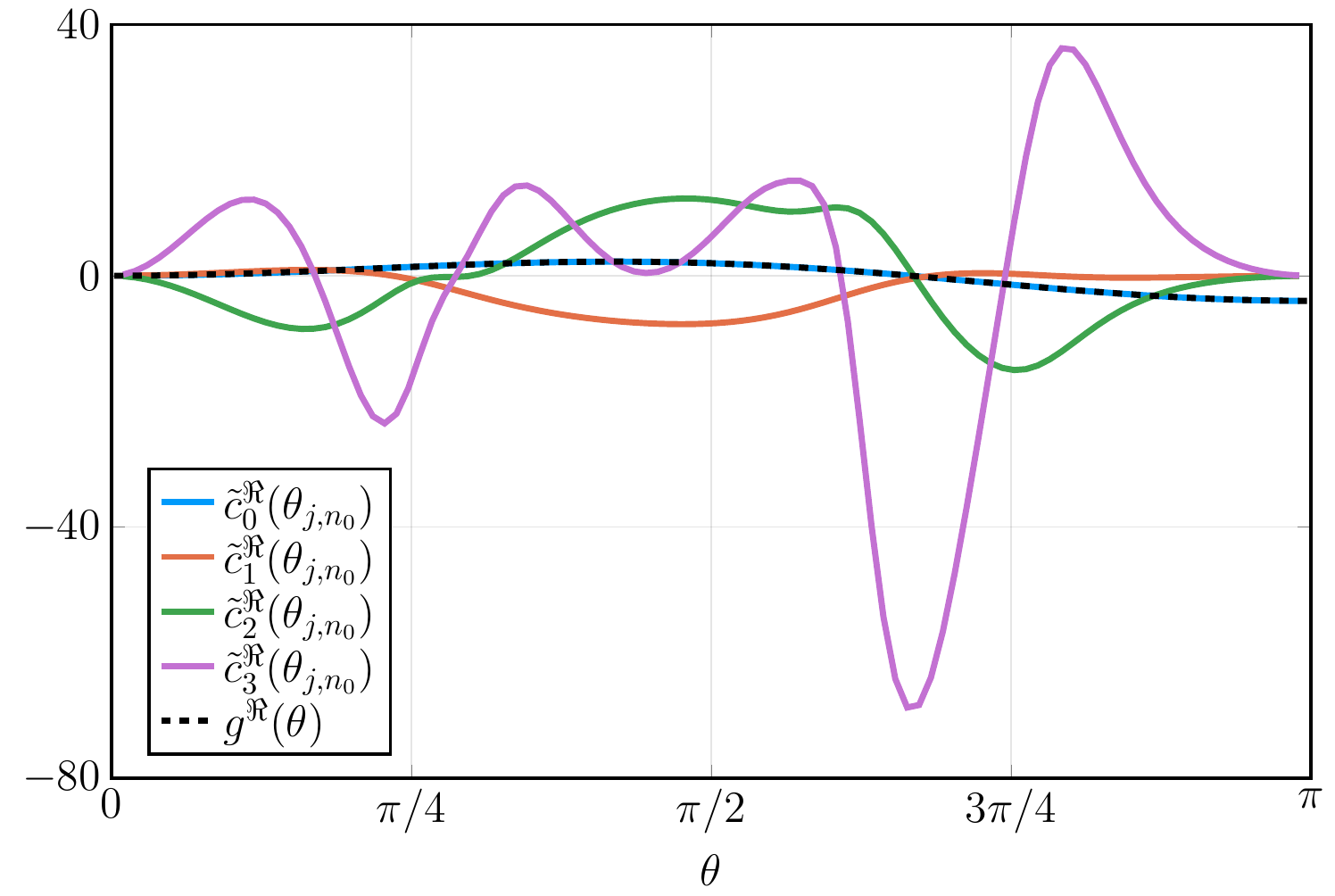}
\includegraphics[width=0.47\textwidth,valign=t]{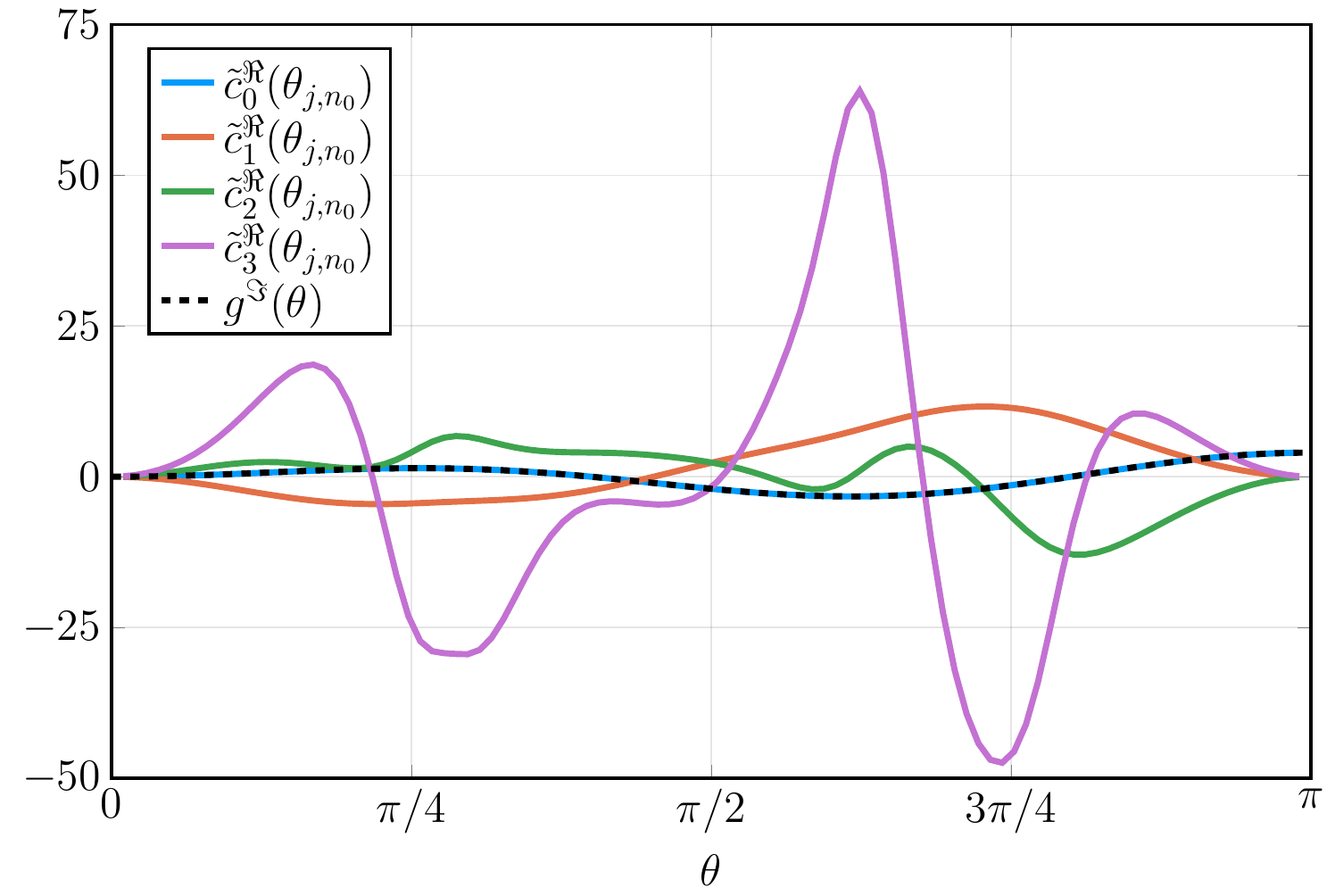}

\includegraphics[width=0.47\textwidth,valign=t]{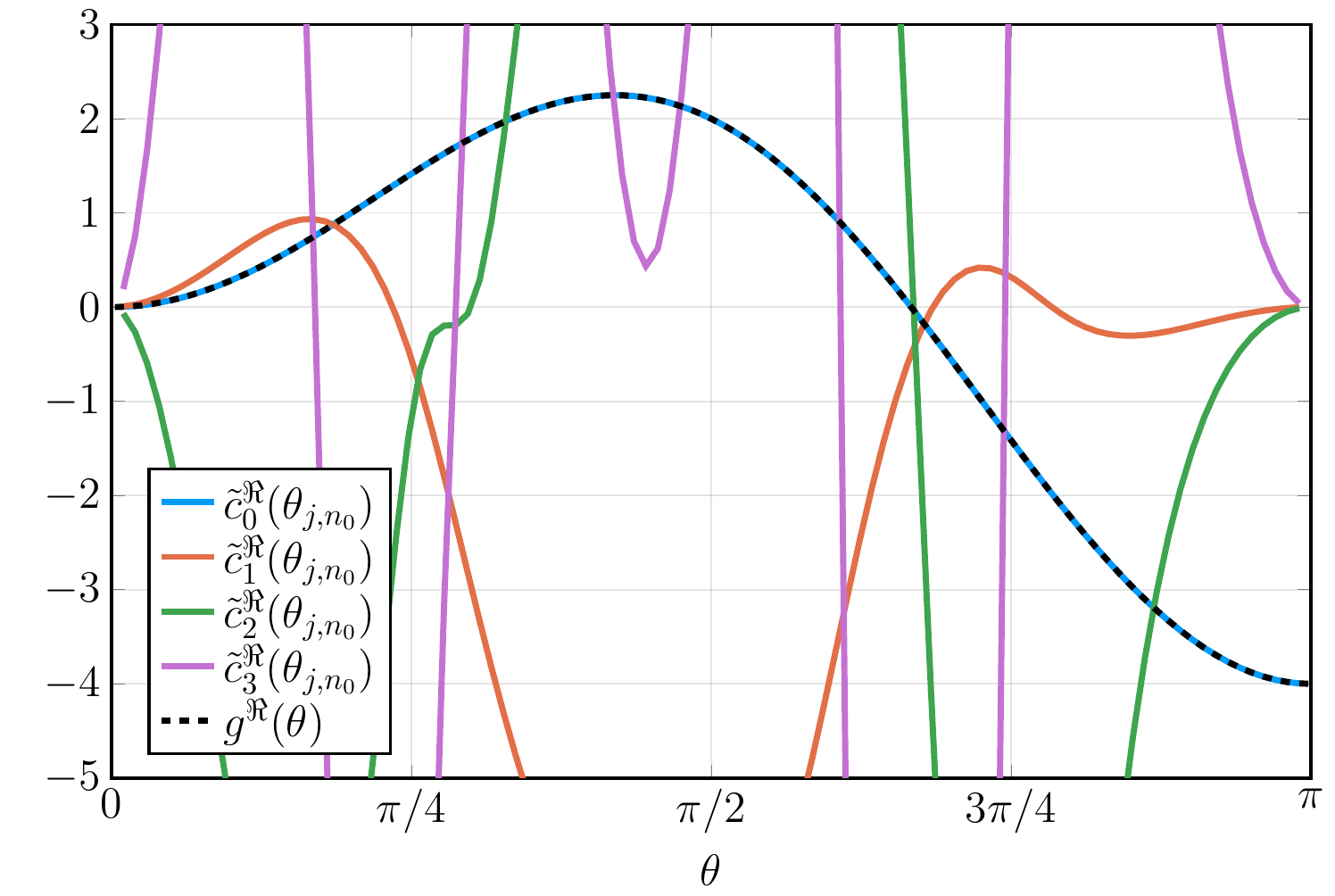}
\includegraphics[width=0.47\textwidth,valign=t]{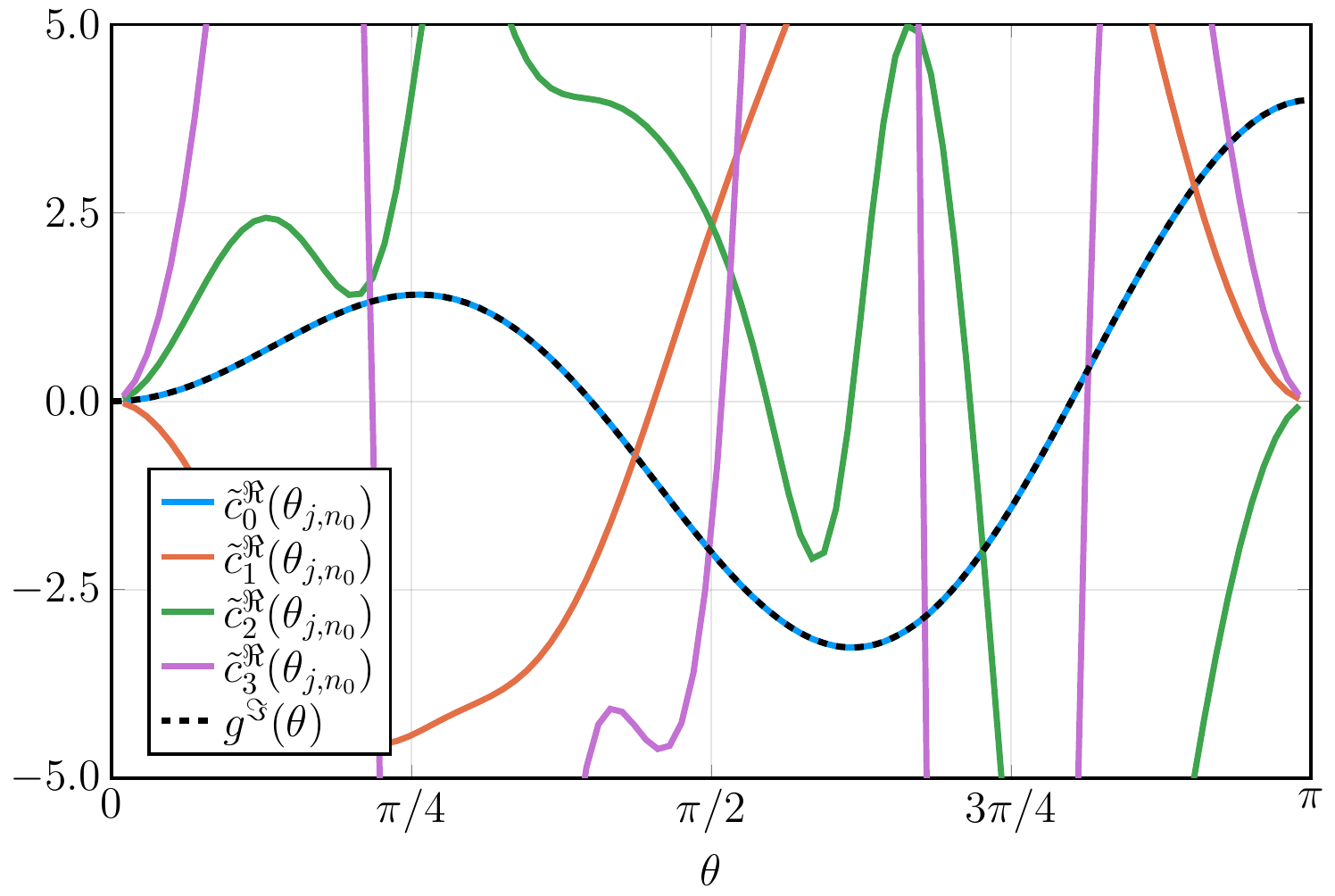}
\caption{[Example~\ref{exmp:7}: Symbol $f(\theta)=2\cos(\theta)-2\cos(2\theta)+\mathbf{i}\left(2\cos(2\theta)-2\cos(3\theta)\right)$]
Top Left:
The computed $\tilde{c}_k^\Re(\theta_{j,n_0})$, $k=0,\ldots\alpha$, with $(n_0,\alpha)=(100,3)$ using Algorithm~\ref{algo:1}.
Top Right: The corresponding $\tilde{c}_k^\Im(\theta_{j,n_0})$ as in the top left panel.
Bottom: Detail of the top panels, clearly showing the overlap of $\tilde{c}_0^\Re(\theta_{j,n_0})$ and $g^\Re(\theta)$ (bottom left) and $\tilde{c}_0^\Im(\theta_{j,n_0})$ and $g^\Im(\theta)$ (bottom right).}
\label{fig:exmp:7:expansion}
\end{figure}
\end{exmp}

\begin{exmp}
\label{exmp:8}
Finally, we return to the Grcar matrix discussed in Example~\ref{exmp:4}, generated by the symbol \eqref{eq:exmp4:symbol}. The \texttt{eigfun} function used as an argument in Algorithm~\ref{algo:1} is:
{\normalfont
{\footnotesize
\begin{lstlisting}[backgroundcolor = \color{jlbackground},
                   language = Julia,
                   xleftmargin = 0.1em,
                   framexleftmargin = 0.1em]
function eigfun_example_4_and_8(n :: Integer, T :: DataType)
  vc  = convert.(T,[1+0im, -1+0im])
  vr  = convert.(T,[1+0im,  1+0im, 1+0im, 1+0im])
  Tn  = toeplitz(n,vc,vr)
  eTn = eigvals(Tn)
  p   = sortperm(imag.(eTn), rev=true)
  return eTn[p]
end
\end{lstlisting}
}}

\noindent In Figure~\ref{fig:exmp:8:expansion} we present the approximated expansion functions in the working hypothesis, for $(n_0,\alpha)=(200,3)$ (512 bits). In the left panel of Figure~\ref{fig:exmp:8:expansion} we show the approximations $\tilde{c}_k^\Re$, and in the right panel $\tilde{c}_k^\Im$. We note that the resolution is rather poor in both $\tilde{c}_k^\Re$ and $\tilde{c}_k^\Im$ for $k=2,3$, at the discontinuities, and $n_0$ should be increased.

\begin{figure}[!ht]
\centering
\includegraphics[width=0.47\textwidth,valign=t]{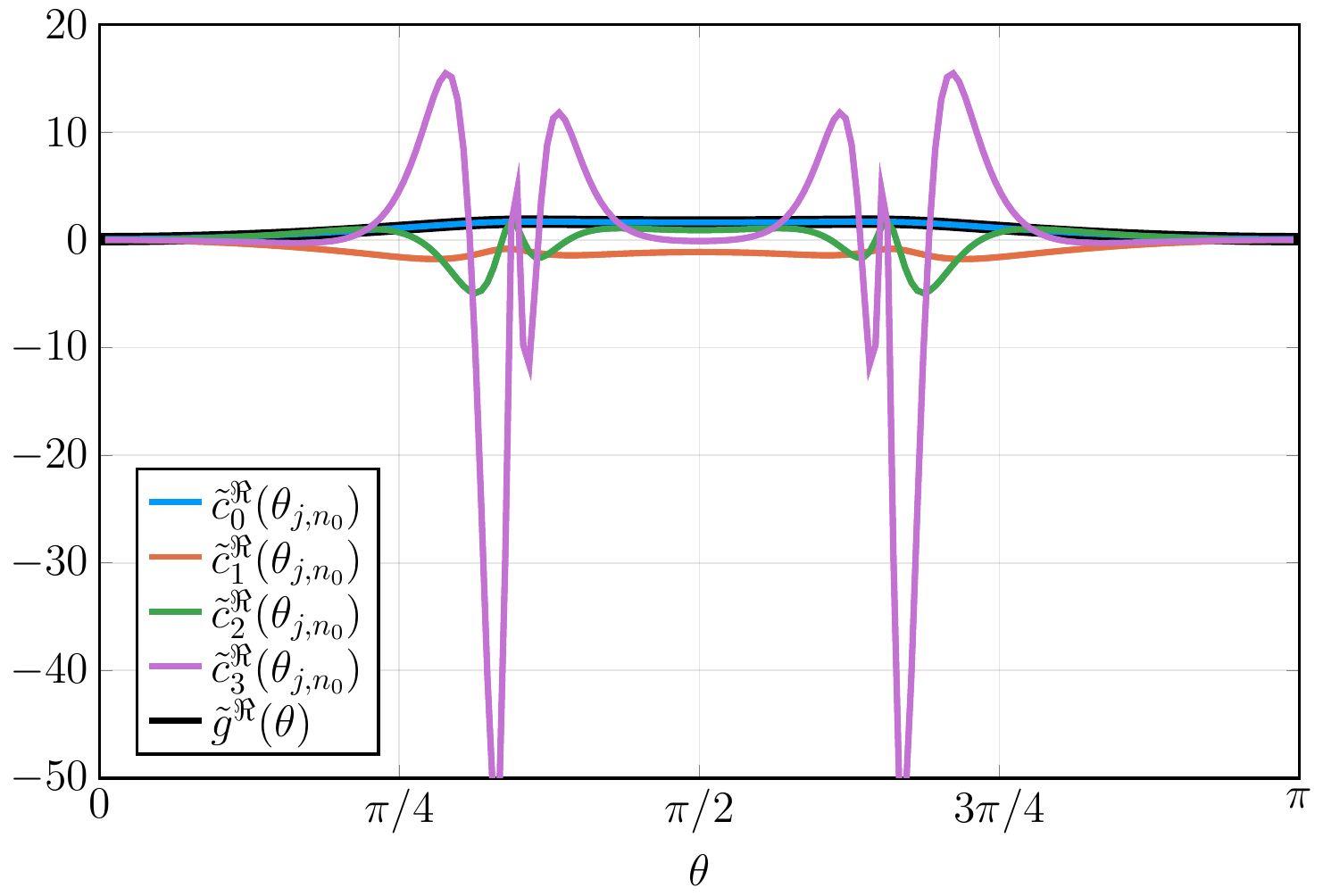}
\includegraphics[width=0.47\textwidth,valign=t]{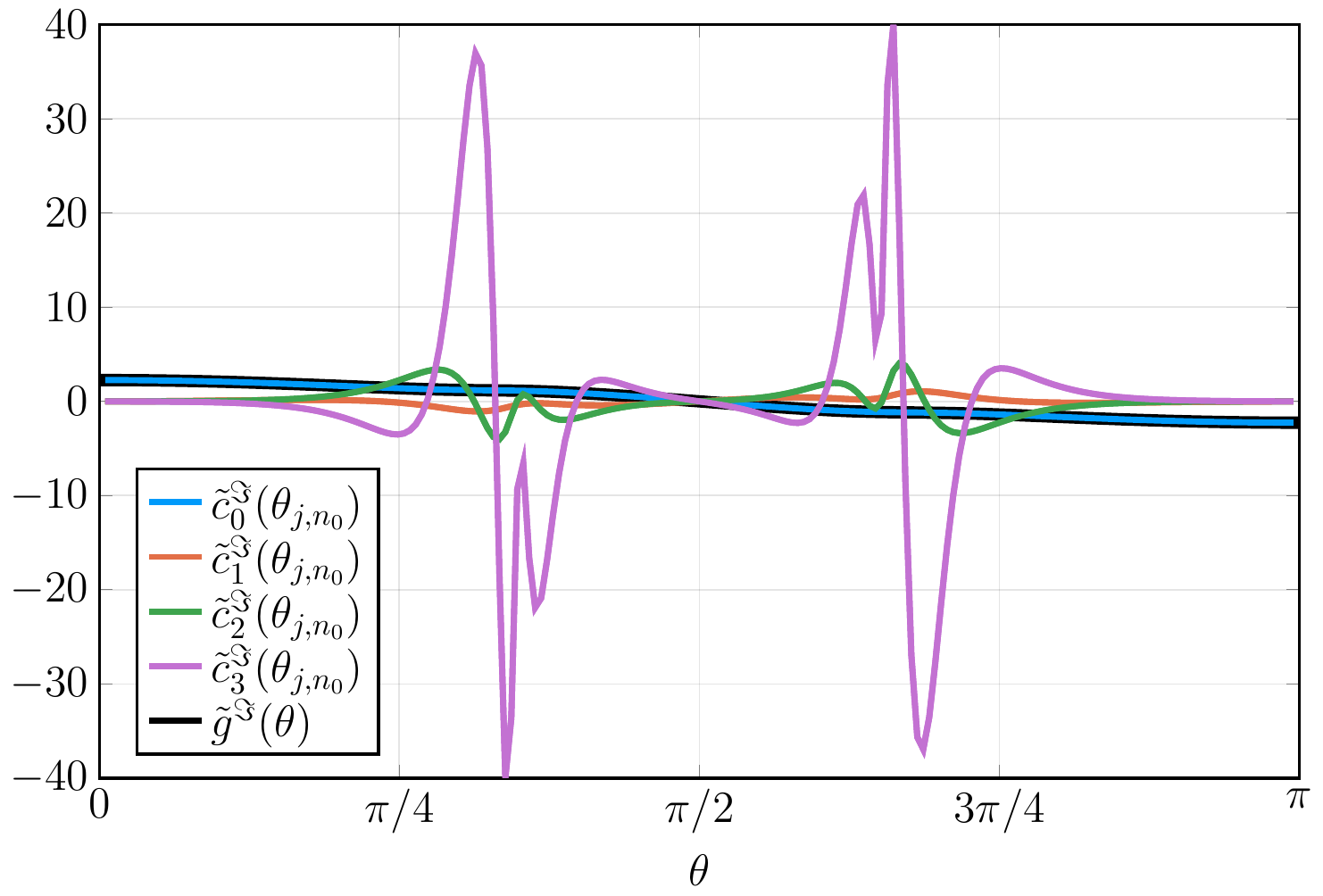}

\includegraphics[width=0.47\textwidth,valign=t]{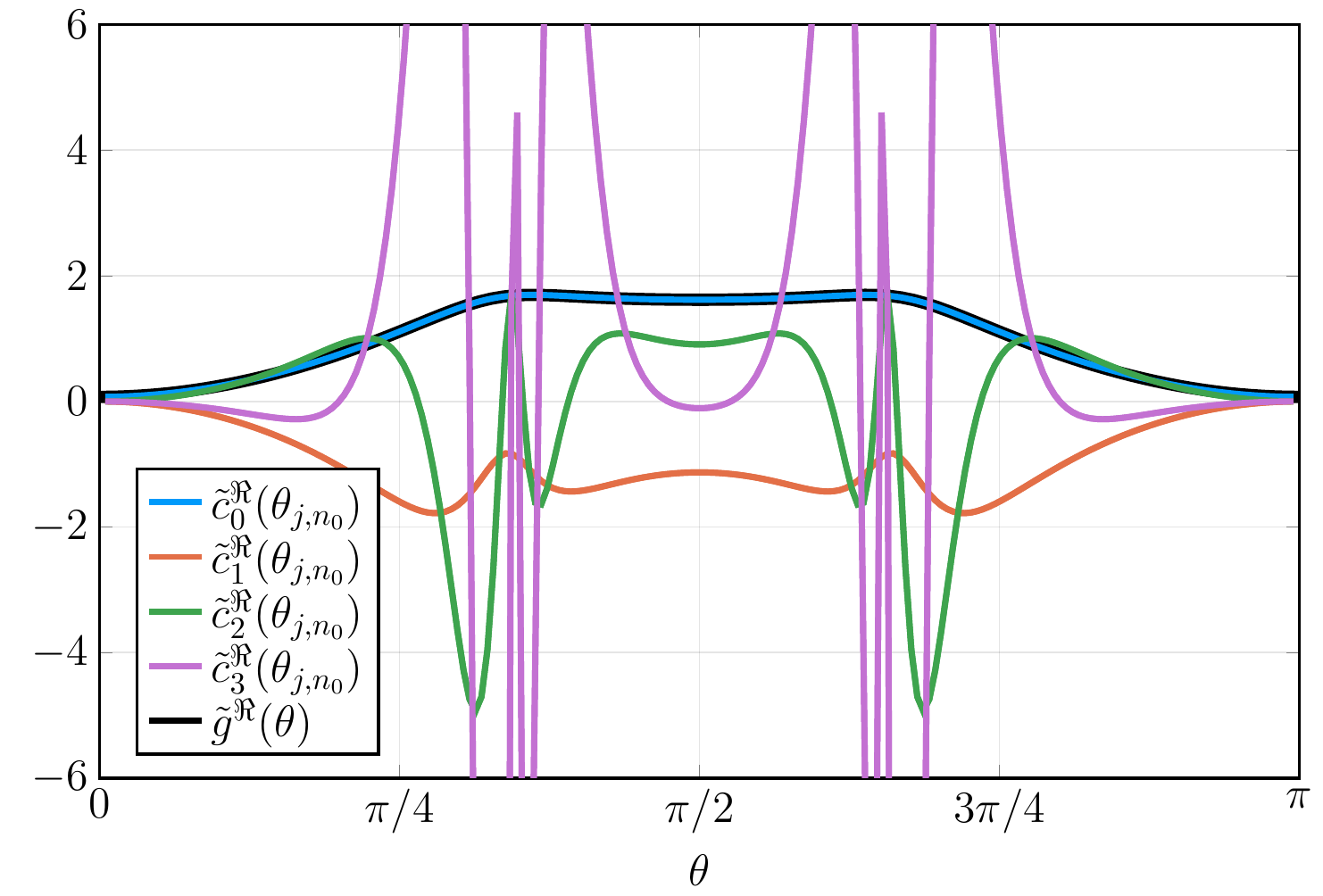}
\includegraphics[width=0.47\textwidth,valign=t]{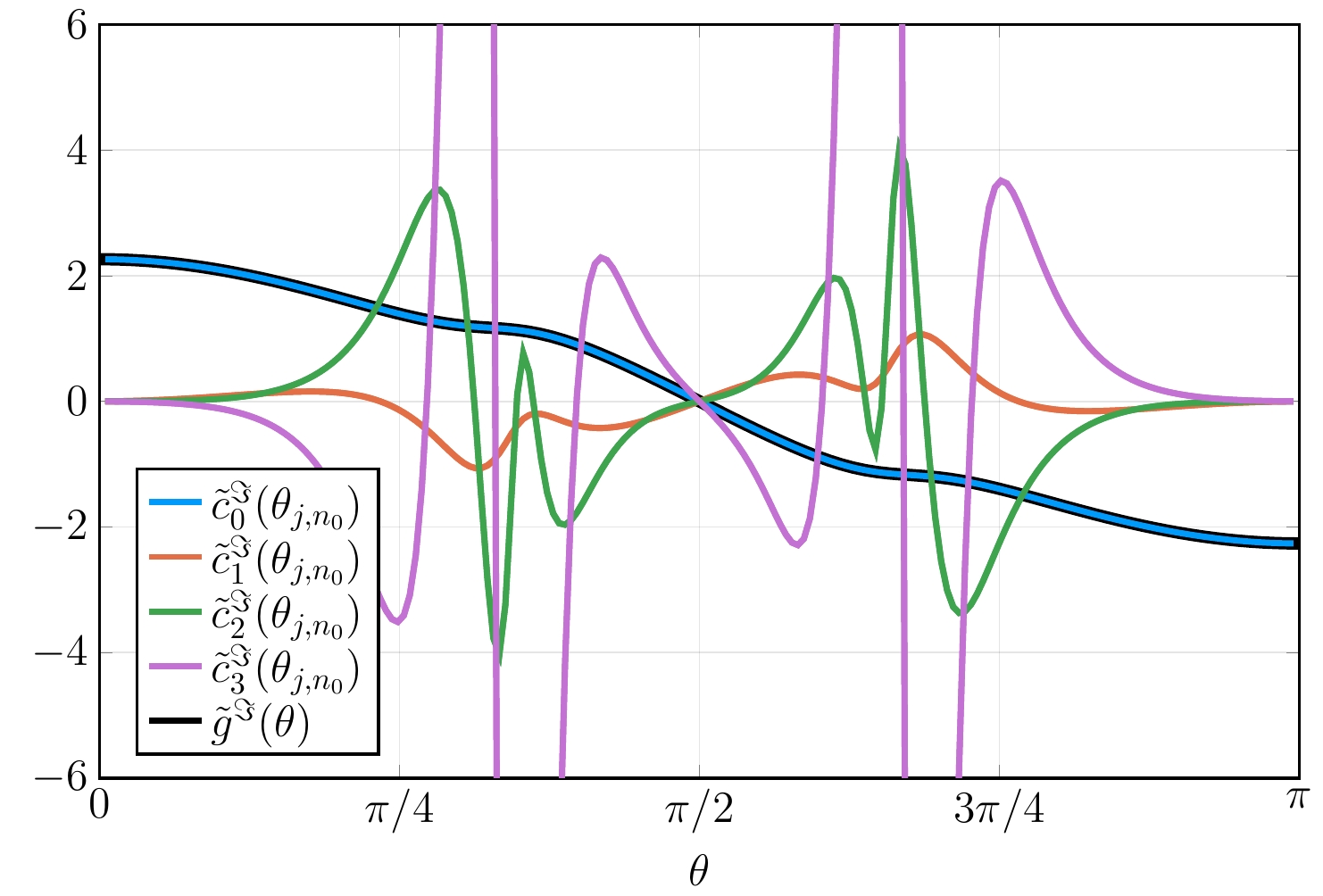}
\caption{[Example~\ref{exmp:8}: Symbol $f(\theta)=
 -\E^{\mathbf{i}\theta}+1+\E^{-\mathbf{i}\theta}+\E^{-2\mathbf{i}\theta}+\E^{-3\mathbf{i}\theta}$]
 Top left: The approximated expansion functions $\tilde{c}_k^\Re(\theta_{j,n_0})$, with $k=0,\ldots,\alpha$ and $(n_0,\alpha)=(100,3)$ (512 bits). 
 Top right: The corresponding $\tilde{c}_k^\Im(\theta_{j,n_0})$ as in the top left panel. Bottom: Detail of the top panels.}
\label{fig:exmp:8:expansion}
\end{figure}

\noindent In Table~\ref{tbl:exmp:8:gk} we present the first ten Fourier coefficients, approximated using Algorithm~\ref{algo:2} with the $\tilde{c}_0(\theta_{j,n_0})$ from Algorithm~\ref{algo:1} using $(n_0,\alpha)=(100,3)$ (512 bits). We note that every second real and imaginary part of $\hat{g}_k$ is zero to machine precision, and we can compactly express $\hat{g}_k$ by a real term times $\mathbf{i}^k$. Computations with $(n_0,\alpha)=(200,3)$ yield the same approximations of the Fourier coefficients as in the table.

\begin{table}[!ht]
\centering
\caption{[Example~\ref{exmp:8}: Symbol $f(\theta)=
 -\E^{\mathbf{i}\theta}+1+\E^{-\mathbf{i}\theta}+\E^{-2\mathbf{i}\theta}+\E^{-3\mathbf{i}\theta}$]
First ten computed ($\tilde{\hat{g}}_k^\Re$ and $\tilde{\hat{g}}_k^\Im$) Fourier coefficients of the unknown $g$. Approximations computed using $(n_0,\alpha)=(100,3)$ (512 bits).
 }
\label{tbl:exmp:8:gk}
\begin{tabular}{crrrr}
\toprule
$k$&\ \hfill$\tilde{\hat{g}}_k^\Re$\hfill\ &\ \hfill$\tilde{\hat{g}}_k^\Im$\hfill\ &\ \hfill$\tilde{\hat{g}}_k=\tilde{\hat{g}}_k^\Re+\mathbf{i}\tilde{\hat{g}}_k^\Im$\hfill\ \\
\midrule
0& $ 1.00000000$&$0$&$1.00000000\mathbf{i}^{k}$\\
1& $ 0$&$1.09011636$&$1.09011636\mathbf{i}^{k}$\\
2& $-0.43169755$&$0$&$0.43169755\mathbf{i}^{k}$\\
3& $ 0$&$-0.00623815$&$0.00623815\mathbf{i}^{k}$\\
4& $-0.07407497$&$0$&$-0.07407497\mathbf{i}^{k}$\\
5& $ 0$&$0.07509827$&$0.07509827\mathbf{i}^{k}$\\
6& $ 0.05451915$&$0$&$-0.05451915\mathbf{i}^{k}$\\
7& $ 0$&$-0.03011810$&$0.03011810\mathbf{i}^{k}$\\
8& $-0.00998665$&$0$&$-0.00998666\mathbf{i}^{k}$\\
9& $ 0$&$-0.00305026$&$-0.00305026\mathbf{i}^{k}$\\
 \bottomrule
\end{tabular}
\end{table}

\noindent In Table~\ref{tbl:exmp:8:g0} we present the numerical values in more detail of the approximated $\tilde{\hat{g}}_0=\tilde{\hat{g}}_0^\Re$ for different choices of $(n_0,\alpha)$ and floating point precision.
Bold digits are assumed to be true (i.e., $\hat{g}_0=1$). Underlined digits are the correct digits for the current $(n_0,\alpha)$. To get the true digits, i.e, zeros, then $(n_0,\alpha)$ has to be increased. The rest of the digits are either correct for the current $(n_0,\alpha)$ (computation with higher precision necessary to verify) or incorrect (computation done with too low precision).
We note that to get a ``converged'' computation with regards to precision, for $(n_0,\alpha)=(100,3)$ with 128 bits of precision, a 256 bits precision computation is required. However, we observe that for $(n_0,\alpha)=(100,3)$ we can only attain nine decimal digits of accuracy (if $\hat{g}_0=1$ is true), no matter what precision is used in the computations. Increasing to $n_0=200$ shows that indeed we get closer to $\hat{g}_0=1$; eleven decimal digits for 256 or 512 bits. Moreover, we note that 128 bits precision is not enough for computing the first Fourier coefficient more accurately than five decimal digits for $(n_0,\alpha)=(200,3)$.
\begin{table}[!ht]
\centering
\caption{[Example~\ref{exmp:8}: Symbol $f(\theta)=
 -\E^{\mathbf{i}\theta}+1+\E^{-\mathbf{i}\theta}+\E^{-2\mathbf{i}\theta}+\E^{-3\mathbf{i}\theta}$] Approximation of first Fourier coefficient $\hat{g}_0$ of the unknown eigenvalue symbol $g(\theta)$, computed with different $(n_0,\alpha)$ and floating point precisions (double precision is 53 bits).}
\label{tbl:exmp:8:g0}
\begin{tabular}{ccrlll}
\toprule
$n_0$&$\alpha$&\textsc{prec}&$\tilde{\hat{g}}_0=\tilde{\hat{g}}_0^\Re$\\
\midrule
100&3&128&$\scriptstyle\mathbf{1.000000000}\uline{96981752360}6915931263983957763$\\
100&3&256&$\scriptstyle\mathbf{1.000000000}\uline{969817523607333664429540902009}118778151114231558758571594613965961456$\\
100&3&512&\makecell{$\scriptstyle\mathbf{1.000000000}\uline{969817523607333664429540902009}091031390554423825856777436964297200837\ldots$\\$\scriptstyle\ldots80730231579597702253008990740101242875408766356306770197490694482736962297945$}\\
\midrule
200&3&128&$\scriptstyle\mathbf{1.00000}4156005612943843265403890568671694$\\
200&3&256&$\scriptstyle\mathbf{1.00000000000}\uline{10080815116696017593133725035485}35843579202835165183921982829075221$\\
200&3&512&\makecell{$\scriptstyle\mathbf{1.00000000000}\uline{10080815116696017593133725035485}97339656341897301003957102835037799\ldots$\\$\scriptstyle\ldots51486808562399236209004605736286999381231589886360905120549590011108832011936$}\\
\bottomrule
\end{tabular}
\end{table}

\noindent In Figure~\ref{fig:exmp:8:fourier} we show the approximated non-zero Fourier coefficients, real (left) and imaginary (right) for $(n_0,\alpha)=(100,3)$  (top) and $(200,3)$ (bottom) (512 bits).
Assuming that $\hat{g}_0=1$ we can estimate where the approximated coefficients are dominated by numerical noise. For $(n_0,\alpha)=(100,3)$, as seen in Table~\ref{tbl:exmp:8:g0}, we have approximately nine digits of accuracy for $\hat{g}_0$, and coefficients with modulus smaller than $10^{-8}$ can be considered to be noise. Indeed, when computing with $(n_0,\alpha)=(200,3)$ we have approximately eleven decimal digits of accuracy for $\hat{g}_0$, and coefficients with modulus smaller than $10^{-10}$ can probably be considered to be noise. The approximated Fourier coefficients $10^{-10}\leq|\tilde{\hat{g}}_k|\leq 10^{-8}$ of the bottom panels seem to follow the general pattern seen in the top panels. Further numerical investigations should be conducted where $(n_0,\alpha)$, and the floating point precision, are increased.
\begin{figure}[!ht]
\centering
\includegraphics[width=0.47\textwidth]{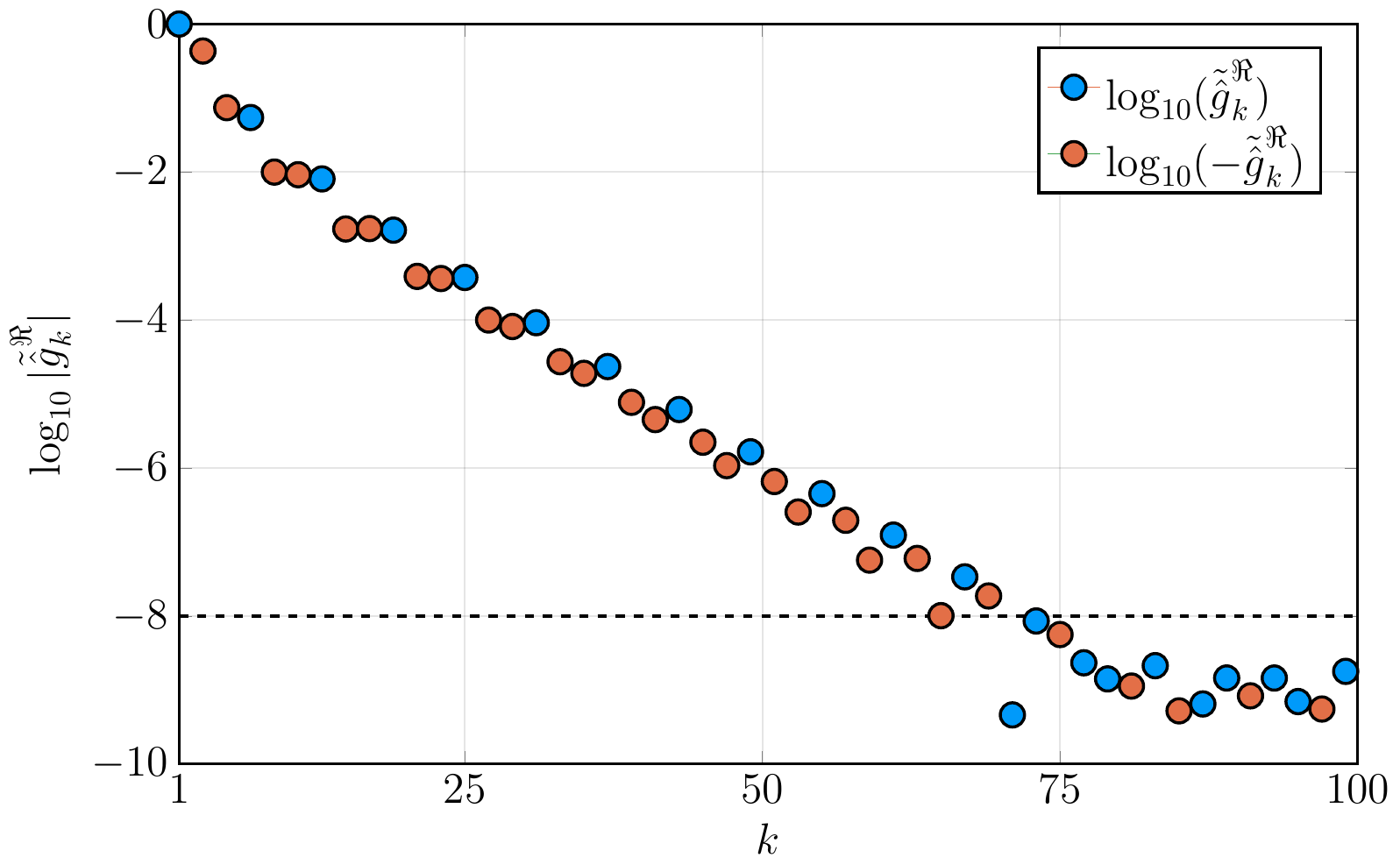}
\includegraphics[width=0.47\textwidth]{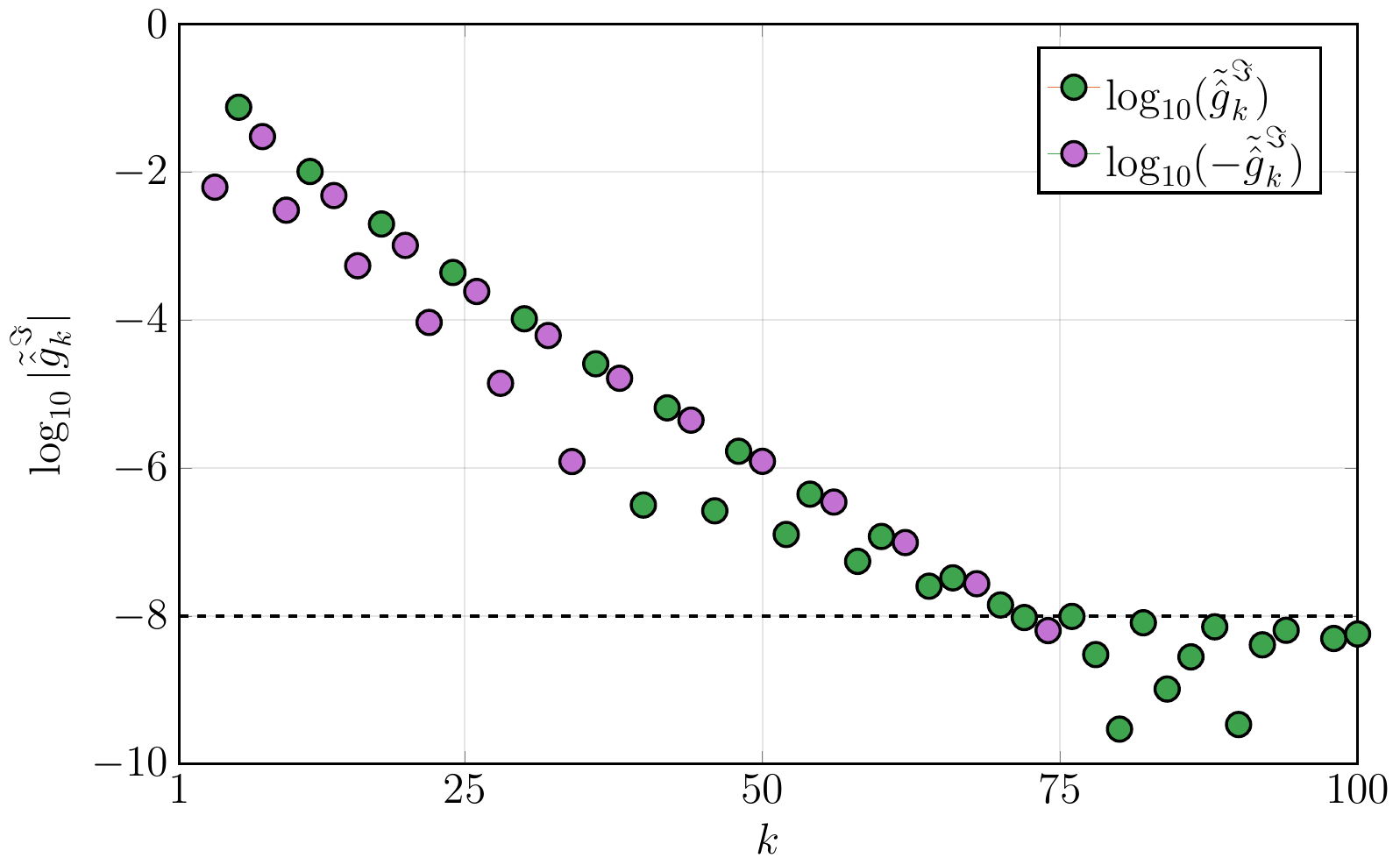}

\includegraphics[width=0.47\textwidth]{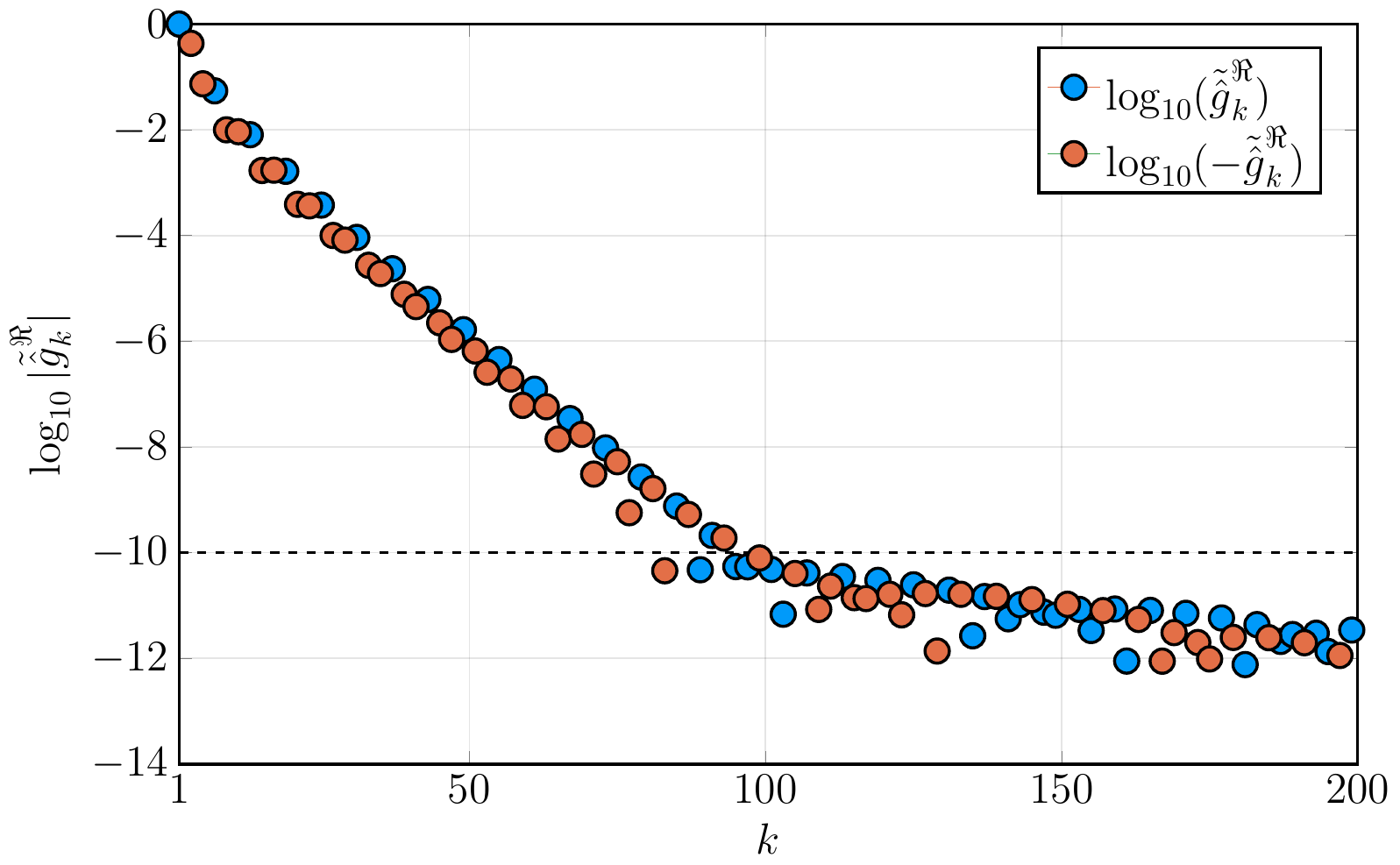}
\includegraphics[width=0.47\textwidth]{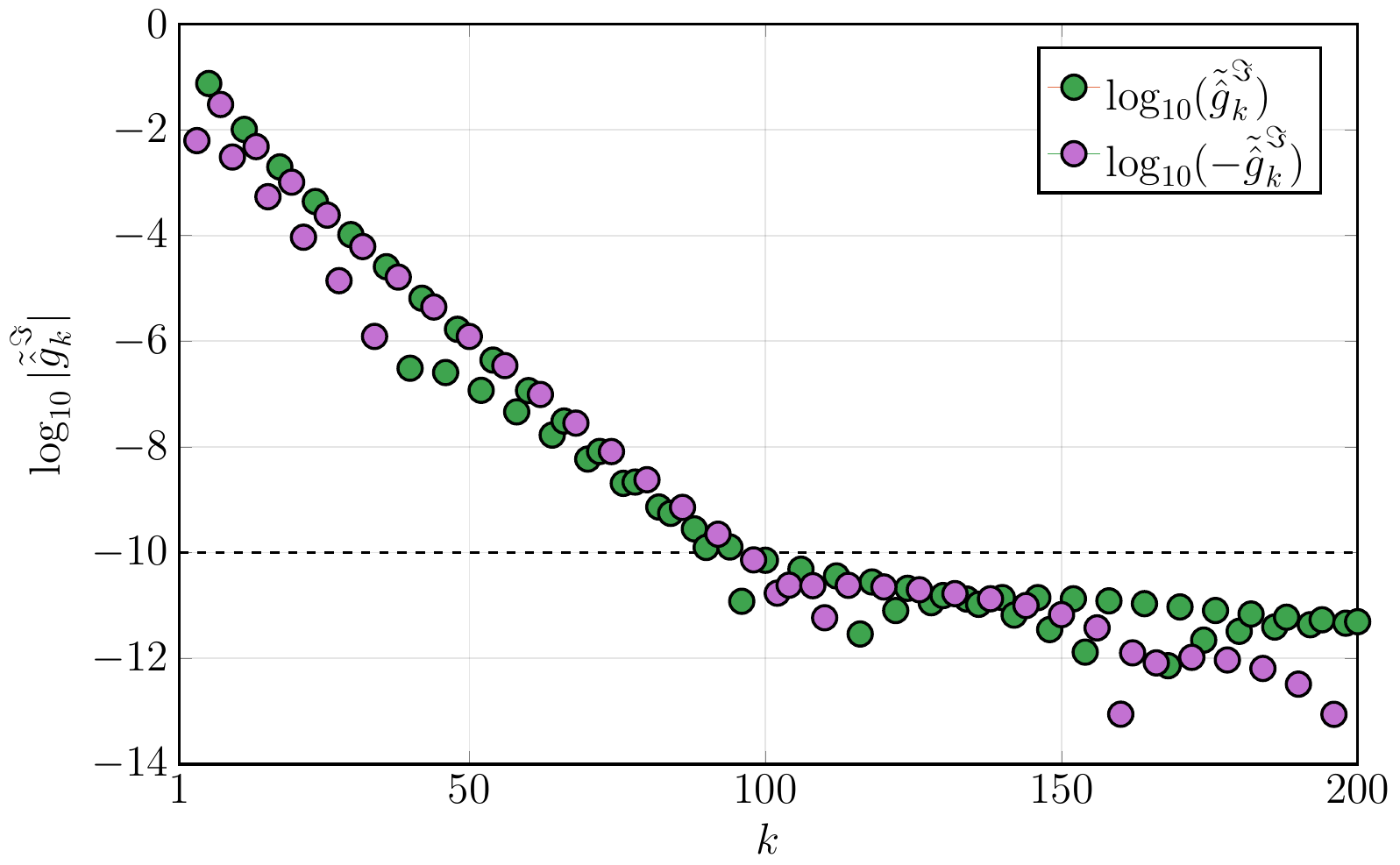}
\caption{[Example~\ref{exmp:8}: Symbol $f(\theta)=
 -\E^{\mathbf{i}\theta}+1+\E^{-\mathbf{i}\theta}+\E^{-2\mathbf{i}\theta}+\E^{-3\mathbf{i}\theta}$]
  The approximated non-zero Fourier coefficients $\tilde{\hat{g}}_k^\Re$ and $\tilde{\hat{g}}_k^\Im$, $k=1,\ldots, n_0$. Top: $(n_0,\alpha)=(100,3)$ (512 bits). Bottom: $(n_0,\alpha)=(200,3)$ (512 bits).}
\label{fig:exmp:8:fourier}
\end{figure}

\noindent Consider the generic symbol $b(z)=-z+1+z^{-1}+z^{-2}+z^{-3}$, where $z\in\mathbb{C}$. The symbol $f(\theta)=
 -\E^{\mathbf{i}\theta}+1+\E^{-\mathbf{i}\theta}+\E^{-2\mathbf{i}\theta}+\E^{-3\mathbf{i}\theta}$ is then given by $f(\theta)=b(\gamma(\theta))$ where $\gamma(\theta)=\E^{\mathbf{i}\theta}$ is a parametrization of the unit circle $\mathbb{T}=\{z\in\mathbb{C}:|z|=1\}$.
We are interested in another $\gamma(\theta)$ such that $g(\theta)=b(\gamma(\theta))$.
In Figure~\ref{fig:exmp:8:portrait} we have a visualization of the four main branches of such possible functions $\gamma_k(\theta)$, for $k=1,\ldots,4$, and $\theta\in[0,\pi]$.
Utilizing the approximations $\tilde{g}(\theta_{j,n_0})$ from Algorithm~\ref{algo:1}, we can draw the four curves $\gamma_k(\theta)$ at the grid points $\theta_{j,n_0}$. These points are given by the four roots of the polynomial $b(z)-\tilde{g}(\theta_{j,n_0})$ for $j=1,\ldots, n_0$.

In the left panel of Figure~\ref{fig:exmp:8:portrait} the four curves $\gamma_k(\theta)$ are shown, with arrows indicating the direction that a point on $\gamma_k(\theta)$ moves as $\theta\in[0,\pi]$ increases (given the computed $\tilde{g}(\theta_{j,n_0})$ of this example with $(n_0,\alpha)=(100,3)$  (512 bits)). For an increasing $\theta\in[-\pi,0]$ the directions are the opposite.

In the right panel of Figure~\ref{fig:exmp:8:portrait} the complex portrait of $b(z)$, using domain coloring~\cite{ludwig191}, is shown. Furthermore, the four $\gamma_k(\theta)$ are indicated in white. A future goal is to find an explicit $\gamma(\theta)$ such that $g(\theta)=b(\gamma(\theta))$, for $\theta\in[-\pi,\pi]$. A possible choice would be the Jordan curve defined by $\gamma_1$ and $\gamma_2$ in the left panel; e.g., $\gamma(\theta)=\gamma_1(\theta+\pi)$ for $\theta\in[-\pi,0]$ and $\gamma(\theta)=\gamma_2(\theta-\pi)$ for $\theta\in[0,\pi]$.
A related approach is \cite[Example 3.]{shapiro171} with symbol $b(z)=z^{-r}(1+az)^{r+s}$, where $r,s\in\mathbb{N}$ and $a\in\mathbb{R}\setminus \{0\}$. Then, $f(\theta)=b(\E^{\mathbf{i}\theta})$, $\gamma(\theta)=\sin(\omega\theta)/\sin((1-\omega)\theta)\E^{\mathbf{i}\theta}$, $\omega=r/(r+s)$, and $g(\theta)=b(\gamma(\theta))=\sin^{r+s}(\theta)/(\sin^r(\omega\theta)\sin^s((1-\omega)\theta))$. The case $(r,s,a)=(3,1,-1)$ is discussed, in the context of approximating the functions $c_k(\theta)$ using matrix-less methods, in~\cite[Examples 3 and 7.]{ekstrom193}.

\begin{figure}[!ht]
\centering
\includegraphics[width=0.39\textwidth]{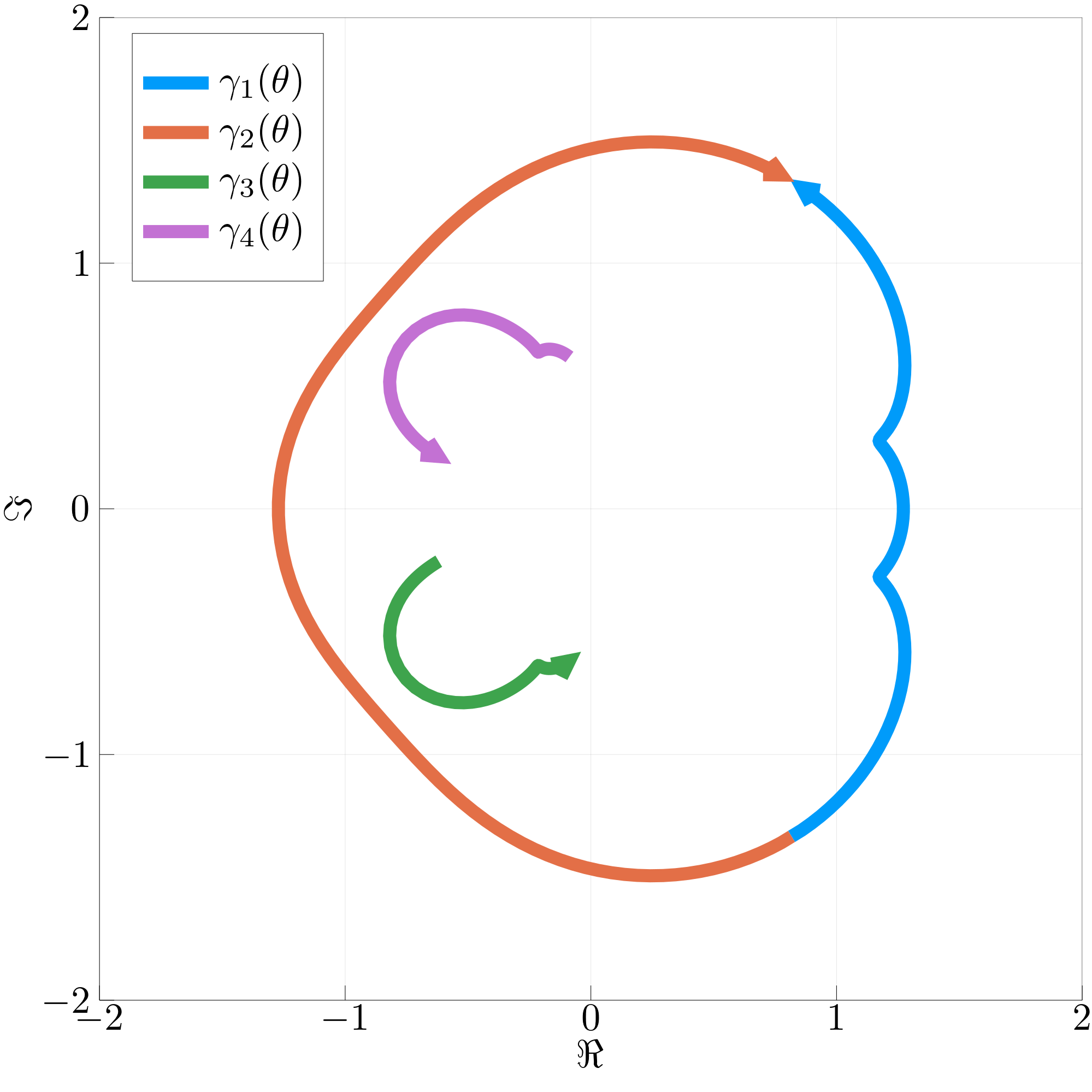}
\includegraphics[width=0.39\textwidth]{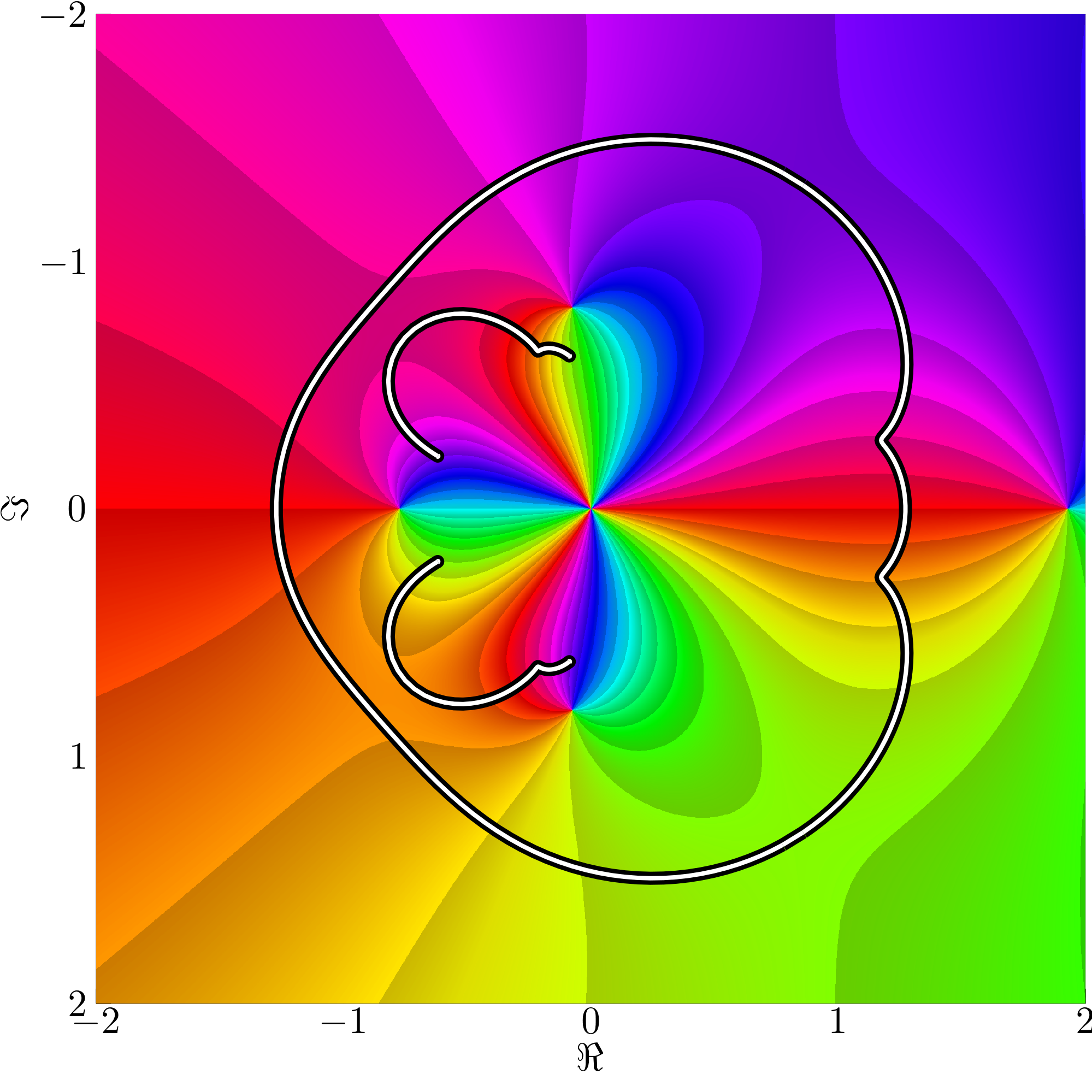}
\caption{[Example~\ref{exmp:8}: Symbol $f(\theta)=
 -\E^{\mathbf{i}\theta}+1+\E^{-\mathbf{i}\theta}+\E^{-2\mathbf{i}\theta}+\E^{-3\mathbf{i}\theta}$]
 We have in generic form the symbol $b(z)=-z+1+z^{-1}+z^{-2}+z^{-3}$ such that $f(\theta)=b(\E^{\mathbf{i}\theta})$. We search for a $\gamma(\theta)$ such that $g(\theta)=b(\gamma(\theta))$.
 Left: The four curves $\gamma_k(\theta)$ approximated by the roots of $b(z)-\tilde{g}(\theta_{j,n_0})$ for $j=1,\ldots,n_0$, with arrows indicating the direction that a point on $\gamma_k(\theta)$ moves as $\theta\in[0,\pi]$ increases.
 Right: The complex portrait, using domain coloring, of $b(z)$ with the four $\gamma_k(\theta)$ indicated in white.
 }
\label{fig:exmp:8:portrait}
\end{figure}
\end{exmp}

\section{Conclusions}
\label{sec:conclusions}
The working hypothesis in the current article concerns the existence of an asymptotic expansion, such that there exists of a function $g$ describing the eigenvalue distribution of the Toeplitz matrices $T_n(f)$ generated by a symbol $f$. The assumption is that $g$ is complex-valued, as opposed to \cite{ekstrom193} where $g$ is assumed to be real-valued.
We have shown in several numerical examples that we can recover an accurate approximation of the function $g$. 
This is done by a matrix-less method described in Algorithm~\ref{algo:1}, where no information of $f$ or $g$ is required, as long as the eigenvalues can be ordered in a consistent way, as $n$ varies. The input argument \texttt{eigfun} could for example encompass preconditioned matrices, matrices generated by matrix-valued symbols, and more complicated matrices, and should be further explored in the future.

Although not presented in the current paper, the approximations of $c_k^\Re$ and $c_k^\Im$ can be used to efficiently and accurately approximate the eigenvalues of $T_n(f)$, for any $n$, using the same interpolation--extrapolation procedure as in~\cite{ekstrom183}, but for the two expansions separately (real and imaginary parts of the spectrum).

In Algorithm~\ref{algo:2} we use the approximation $\tilde{c}_k^\Re(\theta_{j,n_0})$ and $\tilde{c}_k^\Im(\theta_{j,n_0})$ to approximate the Fourier coefficients $\hat{g}_k^\Re$ and $\hat{g}_k^\Im$, and thus reconstructing and approximation of $g$ by its Fourier series.
The presented algorithms can be valuables tool, using high precision computations, for the exploration of the spectral behavior of Toeplitz and Toeplitz-like matrices previously not easily understood; for example the Grcar matrix in Example~\ref{exmp:4} and~\ref{exmp:8}. 

For future research we propose the study of matrices more general than $T_n(f)$, the possibility of using the current results to compute asymptotic expansion of non-monotone real-valued symbols, and finding new explicit expressions for eigenvalue symbols $g$.

\bibliography{References}
\bibliographystyle{siam}
\end{document}